\documentclass[reqno]{amsart}
\usepackage{mathtools}
\usepackage{tkz-euclide}
\usepackage{xfrac}
\usepackage{hyperref}
\usepackage{amsmath}  
\usepackage{amsfonts}
\usepackage{amssymb} 
\usepackage{mathrsfs} 
\usepackage{graphicx}  
\usepackage{bm}
\usepackage{tikz} 
\usetikzlibrary{svg.path} 
\usepackage{amssymb}
\usepackage{pgfplots}
\usepackage{amsmath}
\usepackage[toc]{appendix}
\usetikzlibrary{calc} 
\usepackage{xcolor}  
\usetikzlibrary{arrows.meta} 
\usepackage{amsthm} 
\usepackage{float}
\usepackage{enumitem}
\usepackage{comment}

\usepackage[english]{babel}
\usepackage[font=footnotesize, labelfont=bf]{ caption}
\usepackage{ esint }
\usepackage{stackengine,scalerel,graphicx}
\stackMath

\definecolor{trueblue}{rgb}{0.0, 0.45, 0.81}
\definecolor{truegreen}{rgb}{0.13, 0.55, 0.13}

\newcommand{\NNN}{\color{black}}

\newcommand{\EEE}{\color{black}}

\theoremstyle{plain}
\begingroup
\newtheorem{theorem}{Theorem}[section]
\newtheorem{lemma}[theorem]{Lemma}

\newtheorem{remark}[theorem]{Remark}
\newtheorem{proposition}[theorem]{Proposition}

\endgroup

\usepackage{todonotes}

\theoremstyle{definition}
\begingroup
\newtheorem{definition}[theorem]{Definition}
\endgroup

\renewcommand{\tilde}{\widetilde}

\numberwithin{equation}{section}

\newcommand{\defas}{:=}

\usepackage[normalem]{ulem}

\begin{document}

\author[A. Donnarumma]{Antonio Flavio Donnarumma}
\address[Antonio Flavio Donnarumma]{Department of Mathematics, Friedrich-Alexander Universit\"at Erlangen-N\"urnberg. Cauerstr.~11,
D-91058 Erlangen, Germany}
\email{antonio.flavio.donnarumma@fau.de}

\title[Stochastic homogenization of second order singular perturbation models ]{$\boldsymbol{\Gamma}$-convergence and stochastic homogenization of second order singular perturbation models for phase transitions}
\subjclass[2010]{}
\maketitle
\begin{abstract}  
We study the effective behavior of random, heterogeneous, anisotropic, second order phase transitions energies that arise in the study of pattern formations in physical-chemical systems. Specifically, we study the asymptotic behavior, as $\varepsilon$ goes to zero, of random heterogeneous anisotropic functionals in which the second order perturbation competes not only with a double well potential but also with a possibly negative contribution given by the first order term. We prove that, under suitable  growth conditions and under a stationarity assumption, the functionals $\Gamma$-converge almost surely to a surface energy whose density is independent of the space variable. Furthermore, we show that the limit surface density can be described via a suitable cell formula and is deterministic when ergodicity is assumed.  
\end{abstract}

\section{Introduction}
Since the pioneering work of {\sc Modica and Mortola} \cite{modicamortola} (see also \cite{modica1987gradient}) the mathematical community has been interested in deriving sharp interface limits of diffuse interfaces in which the thickness of the transition layer is sent to zero. The prototypical physical system inspiring such kind of problems is the one formed by two immiscible fluids (oil and water for example) occupying a fixed container and behaving in such a way to minimize the surface area separating them. Mathematically such a system is described by modeling the container with an open bounded Lipschitz subset of $\mathbb{R}^3$, and by setting the configuration function $u$ to be equal to $1$ in the points where the first fluid is present and $-1$ otherwise. Then, the equilibrium configurations are the ones minimizing a surface functional of the form
\begin{equation}
\label{interface}
\mathcal{F}(u)=\sigma \mathcal{H}^{2}(J_u)\:\: \text{for}\:\: u\in BV(A;\{-1,1\}),  
\end{equation}
where $A$ is the container, $\sigma$ is a positive parameter called \emph{surface tension}, $BV(A;\{-1,1\})$ is the space formed by the functions of bounded variation defined on $A$ and taking values in $\{-1,1\}$, $J_u$ is the set of the discontinuity points of $u$, and $\mathcal{H}^{2}(J_u)$ is its 2-dimensional Hausdorff measure. In \cite{modica1987gradient} it was shown that, given a continuous double well potential $W$ with zeros exactly in $-1$ and $1$, the functionals 
\begin{equation}
\label{modicamortolaenergies}
    \mathcal{F}_{\varepsilon}(u)=\int_{A}\frac{W(u)}{\varepsilon}+\varepsilon \vert \nabla u \vert^2 \, \mathrm{d}x\:\: \text{for}\:\: u \in W^{1,2}(A)
\end{equation}
$\Gamma$-converge \cite{DalMaso:93}, with respect to the strong $L^1(A)$ topology, to an energy of the form \eqref{interface} with the surface tension $\sigma$ related to the potential by the following formula:
\begin{equation*}
    \sigma=2\int^{1}_{-1} \sqrt{W(u)}\, \mathrm{d}u.
\end{equation*}
One other important application of phase transition models is the study of pattern formations in physical and chemical systems, like the one occurring in biomembranes and oxygen vacancies in superconductors  \cite[Chapter XIII]{MACHLIN2007399}, \cite{ITsatskis_1994}. For example, one can suppose to have a monolayer $A$, at contact with water, formed by polar molecules having hydrophilic heads (water-attracting) and hydrophobic tails (water-repelling). It has been observed that the molecules arrange themselves in microdomains having different composition and local curvature. One of the first models linking the composition with the curvature of the biomembrane was derived by {\sc Andelman and Seul} in \cite{science1}: specifically, they derived a Modica-Mortola-type free energy, involving the second gradient of the phase variable, whose minimizers describe the equilibrium configurations of the bending biomembrane \NNN (see \cite[(1) and (4)]{science1})\EEE.
Thus, \NNN motivated by the applications mentioned above, \EEE in the last decades the work of Modica and Mortola has been extended in the mathematics literature in order to take into account more complex phenomena. For instance, there are models in which the gradient term is replaced or added to \NNN higher order terms (\EEE see \cite{brusca2024higherordersingularperturbationmodels, Chermisi, cicalesespadarozeppieri, Fonseca_2016, FonsecaMantegazza, GiHaPeZw23}) \EEE and models for heterogeneous, anisotropic and random media (see \cite{Ansini_Braides_ChiadòPiat_2003, cristoferi2023homogenization, cristoferi2020homogenization, marziani2022gammaconvergence, Morfe_2020}). \\
In this paper, we aim at combining some of the previously mentioned works by studying the macroscopic limit of second order phase transitions functionals for heterogeneous, anisotropic, random materials. More precisely, in the first part, we focus on the asymptotic behavior, as $\varepsilon$ goes to zero, of a collection of functionals of the form
\begin{equation}
\label{eq1.3}
\mathcal{E}_{\varepsilon}(u)=\frac{1}{\varepsilon}\int_A f_{\varepsilon}(x,u,\varepsilon\nabla u,\varepsilon^2 \nabla^2u)\, \mathrm{d}x \:\: \text{for}\:\: u\in W^{2,2}(A),
\end{equation}
for a suitable class of integrand functions $(f_{\varepsilon})_{\varepsilon}$.
\NNN We show that for such class of functionals the $\Gamma$-limit exists, up to extracting a subsequence, and that \EEE it is \EEE finite only on functions of bounded variation taking values in $\{-1,1\}$. In addition, such $\Gamma$-limit, denoted with $\mathcal{E}_0$, is a surface integral supported on the discontinuity points of the admissible configurations i.e., \EEE 
\begin{equation*}
    \mathcal{E}_0(u)=\int_{J_u\cap A}f_0\NNN(x,u^+,u^-,\nu_u)\EEE\, \mathrm{d}\mathcal{H}^{n-1}(x)\:\: \text{for}\:\: u \in BV(A;\{-1,1\}),
\end{equation*}
where \NNN $u^+$ and $u^-$ denote the traces on both sides of $J_u$ and \EEE the density $f_0$ can be derived from a suitable cell formula involving the initial energies. \\
In the second part, we focus on random functionals of the form
\begin{equation}
\label{roughrandomenergies}
\mathcal{E}_{\varepsilon}[\omega](u)=\frac{1}{\varepsilon}\int_A f\Big(\omega,\frac{x}{\varepsilon},u,\varepsilon\nabla u,\varepsilon^2 \nabla^2u\Big)\, \mathrm{d}x \:\: \text{for}\:\: u\in W^{2,2}(A),
\end{equation}
where $\omega$ is a random outcome of the sample space $\Omega$ modeling the random environment of the problem. Such functionals are obtained formally from the one in \eqref{eq1.3} by setting $f_{\varepsilon}(x,\cdot,\cdot,\cdot)=f(\frac{x}{\varepsilon},\cdot,\cdot,\cdot)$ and allowing dependence on the random outcome $\omega$. \NNN The energies \eqref{roughrandomenergies} \EEE can be interpreted as the heterogeneous, anisotropic and stochastic version of the functionals introduced by {\sc Chermisi, Dal Maso, Fonseca and Leoni} \cite{Chermisi} and {\sc Fonseca and Mantegazza} \cite{FonsecaMantegazza}, or as the extension of {\sc Marziani} 's energies in \cite{marziani2022gammaconvergence} when also a second gradient perturbation is involved. Indeed, the prototypical functionals $\mathcal{E}_{\varepsilon}[\omega]$ are of the form 
\begin{equation}
\label{andelmanseul}
\int_{A}a\Big(\omega,\frac{x}{\varepsilon},u,\varepsilon\nabla u, \varepsilon^2 \nabla^2 u\Big)\frac{W(u)}{\varepsilon}+ b\Big(\omega,\frac{x}{\varepsilon},u,\varepsilon\nabla u, \varepsilon^2 \nabla^2 u\Big)\varepsilon\big \vert \nabla u \big\vert^2+c\Big(\omega,\frac{x}{\varepsilon},u,\varepsilon\nabla u, \varepsilon^2 \nabla^2 u\Big)\varepsilon^3 \big\vert\nabla^2 u \big\vert^2\, \mathrm{d}x,
\end{equation}
with $a$, $b$, $c \in L^\infty(\Omega\times \mathbb{R}^n\times \mathbb{R}\times \mathbb{R}^n \times \mathbb{R}^{n \times n})$ such that $a(\omega,x,u,\xi,\zeta)\geq a_0$, $b(\omega,x,\xi,\zeta)\in [-q,+\infty)$ for some $q>0$ small enough, and $c(\omega,x,u,\xi,\zeta)\geq c_0$, for every $\omega \in \Omega$, $x\in \mathbb{R}^n$, $u \in \mathbb{R}$, $\xi \in \mathbb{R}^n$, $\zeta \in \mathbb{R}^{n \times n}$, and for some $a_0$, $c_0>0$. Notice that the energy introduced in \cite{Chermisi} is a special case of \eqref{andelmanseul} for $a=c=1$ and $b=\pm q$. \\ 
From a physical perspective, $u$ is a phase variable that takes values close to $1$ in the region where the concentration of the first phase is higher and close to $-1$ in the region where the concentration of the second phase is predominant. In the case of biomembranes \NNN (see for example \cite{science1}\EEE), the coefficient $b$ is related to the monolayer's stiffness and $c$ is a function depending on the surface tension and the bending modulus of the membrane in each point. The fact that $a$, $b$ and $c$ are not just parameters, but functions of $\omega$ and $x$, models the presence of random heterogeneous impurities on the biomembrane. Instead, the coefficients dependence on the gradient and the second gradient takes into account possible anisotropies and dependence of the physical parameters on the local curvature. We point out that the gradient perturbation term can contribute negatively to the total energy. Indeed, curvature instabilities and the creation of domain pattern in the biomembrane are expected when $b<0$  (see \cite{science1},\cite{MACHLIN2007399}). Finally, we point out that considering a double well potential $W$ with zeros in $-1$ and $1$, and thus scalar valued functions, is only for simplicity of notation. As a matter of fact all arguments in the paper can be repeated also when the zeros of the potential consist in two vectors $u_1$ and $u_2$ in $\mathbb{R}^m$, and when scalar valued functions are replaced with functions valued in $\mathbb{R}^m$, for some $m \in \mathbb{N}$.\\
\NNN We now explain \EEE more in details the structure of the paper, the results and the main ideas used in the proofs. The proof's strategy is similar to the one in \cite{marziani2022gammaconvergence} and it is based on the localization method of $\Gamma$-convergence \cite{DalMaso:93}. Specifically, $\Gamma$-convergence is not obtained by proving the $\Gamma-\liminf\limits$ inequality and $\Gamma-\limsup\limits$ inequality explicitly, but by showing that the family of energies is compact with respect to the notion of $\Gamma$-convergence, and that the extracted $\Gamma$-converging subsequence can be chosen independently of the integration domain $A$. Then, we show that the $\Gamma$-limits can be represented in an integral form. This is the content of the first part.
The main difficulties, in applying the localization method to our setting, are related to the fact that the energies densities are not necessarily non-negative and thus they do not define increasing set functions like in \cite{cagnetti2018gammaconvergence}, \cite[Chapters 14-16]{DalMaso:93}, or \cite{marziani2022gammaconvergence}. As a consequence, the proof of the fundamental estimate in Proposition \ref{fundamentalestimate} and the compactness and integral representation Theorem \ref{compintrepr} (which are the analogues of \cite[Proposition 4.1]{marziani2022gammaconvergence} and \cite[Theorem 4.2]{marziani2022gammaconvergence} respectively) become non standard. In order to overcome this problem, the key ingredient we use is \cite[Theorem 1.2]{Chermisi} which states roughly the following: for every open bounded set $A$ with $C^1$ boundary, there exists a sufficiently small $\varepsilon$ such that $\mathcal{E}_{\varepsilon}$ is positive for all configurations. In particular, the strategy is to select each time a sufficiently small $\varepsilon$ such that the energy is positive on the sets difference. In this way, we are able to prove a fundamental estimate, but only for a particular class of $C^1$ sets and for \NNN $\varepsilon$ smaller than some parameter \EEE depending on the sets triple. Nevertheless, this is sufficient for all results in this paper. 
For these reasons, consistently with \cite[Theorem 1.3]{Chermisi}, our $\Gamma$-convergence results hold for $C^1$ subsets and not for Lipschitz subsets like in the original work of Modica and Mortola. 
\\
In the second part, we move to a random environment and prove that, if the energies are of the form \eqref{roughrandomenergies} \NNN and satisfy a stationarity assumption, \EEE then all $\Gamma$-converging subsequences, almost surely, have limit functionals with the same density i.e.\ they converge to the same energy. \NNN In particular, under an additional assumption made only for the purpose of simplifying the notation (see Remark \ref{simplyremark}), \EEE the functionals \eqref{roughrandomenergies} homogenize in the sense that they $\Gamma$-convergence to a functional of the form 
\begin{equation}
\label{roughgammalimit}
\mathcal{E}_{\mathrm{hom}}[\omega](u)=\int_{J_u \cap A}f_{\mathrm{hom}}(\omega,\nu_u)\,\mathrm{d}\mathcal{H}^{n-1}(x)\:\: \text{for}\:\: u\in BV(A;\{-1,1\}),
\end{equation}
with respect to the strong $L^2(A)$ topology. We show also that the homogenized density $f_{\mathrm{hom}}$ can be computed through a cell formula. More precisely, we prove that
\begin{equation}
\label{introfhom}
    f_{\mathrm{hom}}(\omega, \nu)\EEE=\lim\limits_{r \to \infty}\frac{1}{r^{n-1}}\inf_u  \int_{Q^\nu_r(0)}f(\omega,x,u,\nabla u, \nabla^2u)\, \mathrm{d}x,
\end{equation}
where the infimum is taken on all $W^{2,2}$ functions defined on the oriented cube $Q_r^\nu(0)$ and attaining, as boundary datum, a suitable regularization of a piecewise constant function related to the normal vector $\nu$. We point out that $f_{\mathrm{hom}}$ is independent of $x$ \NNN as a consequence \EEE of the stationarity assumption on $f$.
Furthermore, if $f$ is also ergodic, we show that actually $f_{\mathrm{hom}}$ does not depend on $\omega$ i.e., the limit functional is deterministic. The proof strategy relies mostly on well-known techniques in the field of homogenization, except for taking into account two observations. The first one is that the notion of subadditive process, typically used in the variational approach to stochastic homogenization (see for example \cite[(7.5)]{marziani2022gammaconvergence}), is defined for positive energies evaluated at the scale $\varepsilon=1$, while our energies are in principle not bounded from below by zero. The key ingredient used to adapt the standard subadditive process to our setting is \cite[Proposition 3.9]{Chermisi}, which guarantees the positivity of the energies on sufficiently large rectangles. The second observation is that, due to the fact that the cell formula \ref{introfhom} is formulated on cubes and thus on sets which are not $C^1$, we need to adapt the previously mentioned $\Gamma$-convergence result to cubic domains. For this purpose, we use the fact that sets with $C^1$ boundaries are dense, in the sense of Definition \ref{setsdensity}, and that the energies of the sets difference of two cubes can be controlled with their surface uniformly in $\varepsilon$.    
We also remark that, in contrast to \cite{marziani2022gammaconvergence}, we also examine the case in which the random density is stationary with respect to a continuous group. \\
We mention that models involving two different scale parameters $\delta$ and $\varepsilon$ have been studied in the context of first order perturbation models and deterministic homogenization \cite{Ansini_Braides_ChiadòPiat_2003, supcristoferi2023homogenization, cristoferi2023homogenization,cristoferi2020homogenization} where one considers functionals of the type 
\begin{equation*}
 \mathcal{E}_{\delta,\varepsilon}(u)=\frac{1}{\varepsilon}\int_A f\Big(\frac{x}{\delta},u,\varepsilon\nabla u\Big)\, \mathrm{d}x,   
\end{equation*}
and the corresponding regimes $\varepsilon \gg \delta$ and $\delta\gg\varepsilon$.  The generalization to stochastic environments and second order models is still missing and it is beyond the scope of this work. \\
The paper, counting also this introduction, is divided in six sections. In Section \ref{sec: not} we expose the basic notation. In Section \ref{sec: setting of the problem and main results} we present the main notions and main results of this work which can be categorized in two parts: first a deterministic part in which we show compactness \NNN with \EEE respect to the notion of $\Gamma$-convergence and integral representation of the $\Gamma$-limit. Then, a stochastic part in which we state that there exists an event of probability 1 such that all subsequences converge to the same $\Gamma$-limit. 
The main results will be proved in Sections \ref{sec: gamma convergence} and \ref{sec: stochastichomogenization}. Finally, in Appendix \ref{sec: appendix} we collect some auxiliary results used in Sections \ref{sec: gamma convergence}, \ref{sec: stochastichomogenization}.

 \section{Notation}
 \label{sec: not}
We introduce basic notation. Let $n \in \mathbb{N}$. Given $x \in \mathbb{R}^n$ we denote by $|x|$ its Euclidean norm. For every $x,y \in \mathbb{R}^n$, $x \cdot y$ denotes the standard inner product on $\mathbb{R}^n$ between $x$ and $y$. For  $A,B \subset \mathbb{R}^n$ and $\lambda \in \mathbb{R}$,  we define
 \begin{equation*}
  A+B:=\{z \in \mathbb{R}^n: z=x+y, \ \ \   x   \in A\:\: \text{and}\:\: y \in B\}
 \end{equation*}
and 
 \begin{equation*}
 \lambda A:=\{ z \in \mathbb{R}^n: z=\lambda x,\:\: x \in A\}.
 \end{equation*}
By $A \setminus B$ we denote the difference between the two sets $A$ and $B$. With $A \triangle B$ we denote their symmetric difference. 
  We write $A \subset \subset B$ if $\overline{A} \subset B$, where $\overline{A}$ is the closure of $A$.  By $\mathbb{S}^{n-1} =  \{x \in \mathbb{R}^n: |x|=1\}$ we denote the unit sphere in $\mathbb{R}^n$. Then, we define its subsets
  \begin{equation*}
      \hat{\mathbb{S}}_{\pm}^{n-1}\defas \{x \in \mathbb{S}^{n-1}\colon \pm x_{i(x)}>0\},
  \end{equation*}
  where $i(x)$ is the largest index in $\{1,...,n\}$ such that $x_{i(x)}\neq 0$. Notice that $\mathbb{S}^{n-1}=\hat{\mathbb{S}}_{-}^{n-1}\cup \hat{\mathbb{S}}^{n-1}_{+}$.\\
Given $x \in \mathbb{R}^n$ and $\rho>0$ we indicate with $Q_{\rho}(x)$ the open cube, with center in $x$ and side length $\rho$, oriented according to the canonical orthonormal basis $\{e_1,...,e_n\}$, that is
 \begin{equation*}
 Q_{\rho}(x):=\Big\{ y \in \mathbb{R}^n: \max_{i=1,...,n}|y_i - x_i|<\frac{\rho}{2}\Big\}.
 \end{equation*}
 If $x=0$ and $\rho=1$ we will write $Q$ instead of $Q_1(0)$.
 For $\nu \in \mathbb{S}^{n-1}$ we fix an orthogonal matrix $R_{\nu}$ such that $R_{\nu}(e_n)=\nu$ and satisfying the following properties:
 \begin{itemize}
     \item the restrictions of the function $\nu \mapsto R_{\nu}$ to the sets $\hat{\mathbb{S}}_{\pm}^{n-1}$ are continuous,
     \item $R_{-\nu}Q=R_{\nu}Q$ for every $\nu \in \mathbb{S}^{n-1}$ and for every cube $Q$ with center in zero,
     \item $R_{\nu}\in O(n)\cap \mathbb{Q}^{n \times n}$ for every $\nu \in \mathbb{S}^{n-1}\cap \mathbb{Q}^n$.
 \end{itemize}
 An example of such a map $\nu \mapsto R_{\nu}$ satisfying these assumptions is provided in \cite[Example A.1 and Remark A.2]{cagnetti2018gammaconvergence}. For every $\nu \in \mathbb{S}^{n-1}\cap \mathbb{Q}^{n}$ we denote with $M_{\nu}$ the smallest integer, larger than 2, such that $M_{\nu}R_{\nu}$ belongs to $\mathbb{Z}^{n \times n}$.\\
 Then, we denote by $Q^{\nu}_{\rho}(x)$ the $n$ dimensional open cube, with center in $x$ and side length $\rho$, oriented according to the orthonormal basis $\{ R_{\nu}(e_1),...,\nu\}$, that is
 \begin{equation}\label{eq: Qnot}
     Q_{\rho}^{\nu}(x)=R_{\nu}Q_{\rho}(0)+x.
 \end{equation}
Similarly, with $Q_{\rho}^\prime$ we indicate the $(n-1)$-dimensional open cube, with center in $0$ and side length $\rho$. In this case, we will omit the subscript $\rho$ if $\rho=1$. 
For every $x \in \mathbb{R}^n$ and $\nu \in \mathbb{S}^{n-1}$ we indicate with $\Pi_x^{\nu}$ the hyperplane passing through $x$ and orthogonal to $\nu$, i.e.,
 \begin{equation*}
     \Pi_x^{\nu}\defas \{y \in \mathbb{R}^n \colon (y-x)\cdot \nu =0\}.
 \end{equation*}
 If $x=0$ we will just write $\Pi^{\nu}$ instead of $\Pi_0^{\nu}$.
 We proceed with further notation for sets and measures: 
 
 \begin{enumerate}

     \item[$(\mathrm{a})$] By  $\mathcal{A}$ we denote the family of all open, bounded subsets of $\mathbb{R}^n$, and by $\mathcal{A}_1$  the family of all  open, bounded subsets of $\mathbb{R}^n$ with $C^1$ boundary.
\item[$(\mathrm{b})$] By $\mathcal{L}^k$ and $\mathcal{H}^k$ we indicate, respectively, the $k$-dimensional Lebesgue and Hausdorff measure. 
\item[$(\mathrm{c})$] \NNN Given a topological space $X$, we indicate with $\mathcal{B}(X)$ its Borel $\sigma$-algebra. We denote with $\mathcal{B}^n$ the Borel $\sigma$-algebra of $\mathbb{R}^n$ and with $\mathcal{B}(\mathbb{S}^{n-1})$ the Borel $\sigma$-algebra of $\mathbb{S}^{n-1}$. \EEE
\item[$(\mathrm{d})$] Given $n$ measurable spaces $(X_1,\Sigma_1),...,(X_n,\Sigma_n)$, we denote with $\Sigma_1 \otimes ... \otimes \Sigma_n$ the product $\sigma$-algebra on $X_1 \times ... \times X_n$. 
 \end{enumerate}
 We denote with $BV(A;\{-1,1\})$ the set of all the functions $u\in L^1(A)$ having bounded variation and such that $u(x)\in \{-1,1\}$ for $\mathcal{L}^n$-a.e $x \in A$. For every $u \in BV(A;\{-1,1\})$, we indicate by $J_u$ the set of its approximate jump points,  see \cite{ambrosio2000fbv} for the definition.
\section{Setting of the problem and main results}
\label{sec: setting of the problem and main results}
In this section we introduce the setting and we present the main results. \NNN
\subsection{Setting in deterministic environments and preliminary results} \EEE
Let $W \colon \mathbb{R}\to [0,\infty)$ be a double well potential such that 
\begin{enumerate}
    \item[$(W_1)$] $W$ is continuous,
    \item[$(W_2)$] $W^{-1}(0)=\{-1,1\}$,
     \item[$(W_3)$] $W(s)\geq (\vert s \vert -1)^2$ for every $s \in \mathbb{R}$,
    \item[$(W_4)$] There exists a $c_0 \geq 1$ such that $W(s)\leq c_0 W(t)+c_0$ for every $t \in \mathbb{R}$ and every $s \in \mathbb{R}$ with $\vert s \vert \leq \vert t \vert$. 
\end{enumerate}
One example for $W$ satisfying $(W_1)$--$(W_4)$ is $W(s)=(s^2-1)^2$.
We point out that conditions $(W_3)$ and $(W_4)$ are additional technical assumptions needed for Theorem \ref{compintrepr} below, see also the $\Gamma$-convergence result \cite[Theorem 1.3]{Chermisi}. \\
We recall here a direct consequence of \cite[Proposition 3.9]{Chermisi} and \cite[Theorem 1.2]{Chermisi} and adapt them according to our notation. 
\begin{theorem}
\label{positivesubadditiveprocess}
Assume that $W$ satisfies $(W_1)$--$(W_3)$ and let $R$ be an open rectangle having all sides of at least length 1. Then, there exists a constant $q^*>0$, independent of $R$, such that for every $-\infty<q<\frac{q^*}{n}$ and for every $u \in W^{2,2}(R)$ it holds 
\begin{equation*}
\int_{R}W(u)-q\vert \nabla u \vert^2+\vert \nabla^2 u \vert^2 \, \mathrm{d}x \geq 0.
\end{equation*}
\end{theorem}
\begin{theorem}
\label{chermisitheorem}
Let $A \in \mathcal{A}_1$ and assume that $W$ satisfies $(W_1)$--$(W_3)$. Then there exists a constant $q^*>0$, independent of $A$, such that for every $-\infty<q<\frac{q^*}{n}$ there exists $\varepsilon_0=\varepsilon_0(A,q)$ such that
\begin{equation*}
    q\varepsilon^2\int_{A}\vert \nabla u \vert^2 \, \mathrm{d}x \leq \int_{A}W(u)\, \mathrm{d}x+\varepsilon^4\int_{A}\vert \nabla^2 u \vert^2 \, \mathrm{d}x
\end{equation*}
for every $\varepsilon \in (0,\varepsilon_0)$ and $u \in W^{2,2}(A)$.
\end{theorem}
Let $q^*$ such that Theorem \ref{positivesubadditiveprocess} and Theorem \ref{chermisitheorem} hold. Let $0<q< \frac{q^*}{n}$.
Let $0<c_1\leq c_2$ and $\mathcal{F}=\mathcal{F}(W,c_1,c_2,q)$ be the family of functions $f \colon \mathbb{R}^n \times \mathbb{R}\times \mathbb{R}^n \times \mathbb{R}^{n \times n} \to \mathbb{R}$ satisfying the following hypotheses:
\begin{enumerate}
    \item[$(f_1)$] (measurability) $f$ is Borel measurable on $\mathbb{R}^n \times \mathbb{R}\times \mathbb{R}^n \times \mathbb{R}^{n \times n}$, 
    \item[$\NNN (f_2)\EEE$] (lower bound) for every $x \in \mathbb{R}^n$, $u \in \mathbb{R}$, $\xi \in \mathbb{R}^n$, and $\zeta \in \mathbb{R}^{n \times n}$ it holds
    \begin{equation*}
       c_1(W(u)-q\vert \xi \vert^2+\vert \zeta \vert^2) \leq f(x,u,\xi,\zeta),
    \end{equation*}
    \item[$\NNN (f_3)\EEE$] (upper bound) for every $x \in \mathbb{R}^n$, $u \in \mathbb{R}$, $\xi \in \mathbb{R}^n$, and $\zeta \in \mathbb{R}^{n \times n}$ it holds
    \begin{equation*}
    f(x,u,\xi, \zeta)\leq c_2(W(u)+q\vert \xi \vert^2+\vert \zeta \vert^2).
    \end{equation*}
\end{enumerate} 
Consider a collection of functions $(f_\varepsilon)_{\varepsilon>0}$ in $\mathcal{F}$ and the corresponding functionals $\mathcal{E}_\varepsilon \colon L_{\mathrm{loc}}^2(\mathbb{R}^n)\times \mathcal{A} \to \mathbb{R}$ defined as
    \begin{equation}
    \label{defenergies}    
    \mathcal{E}_\varepsilon(u,A)=\begin{cases}
        \frac{1}{\varepsilon}\int_A f_\varepsilon(x,u,\varepsilon\nabla u,\varepsilon^2 \nabla^2u)\, \mathrm{d}x \quad \text{if}\:\: u_{\vert_A}\in W^{2,2}(A)\\
        +\infty \quad \quad \quad\quad \quad \quad \quad \quad \quad \quad\quad \:\:\text{otherwise}.
    \end{cases}
    \end{equation}
Given $\varepsilon >0$, let us denote by $\mathcal{M}^\pm_\varepsilon \colon L_{\mathrm{loc}}^2(\mathbb{R}^n)\times \mathcal{A}\to \mathbb{R}$ two functionals of the type introduced in \cite{Chermisi}, namely $\mathcal{M}^\pm_\varepsilon$ are defined as
\begin{equation}
\label{Chermisi}
 \mathcal{M}^\pm_\varepsilon(u,A)=\begin{cases}
        \int_A \frac{W(u)}{\varepsilon}\pm q\varepsilon\vert \nabla u\vert^2+\varepsilon^3\vert \nabla^2 u\vert^2\mathrm{d}x \quad \text{if}\:\: u_{\vert_A}\in W^{2,2}(A)\\
        +\infty \quad \quad \quad \quad \quad \quad \quad \quad\quad \quad \quad \quad \quad \:\:\text{otherwise}.
        \end{cases}
\end{equation}
In \cite[Theorem 1.3]{Chermisi} it has been proved that, for every $A\in \mathcal{A}_1$, the functionals $\mathcal{M}^\pm_\varepsilon(\cdot,A)$ $\Gamma$-converge, with respect to the strong $L^2(A)$ topology, to the surface functional
\begin{equation}
\label{FonsecaMantegazzalimit}
\mathcal{M}^\pm_{0}(u,A)=\begin{cases}
        \sigma^\pm\mathcal{H}^{n-1}(J_u \cap A) \quad \text{if}\:\: u_{\vert_A}\in BV(A;\{-1,1\})\\
        +\infty  \:\:\:\quad \quad \quad \quad \quad \quad \text{otherwise},
        \end{cases}
\end{equation}
where $\sigma^\pm$ is the surface tension defined in \cite[(1.4)]{Chermisi} ($m_N$ in the authors' notation) i.e.,
\begin{equation}
\label{chermisidensity}
   \sigma^\pm \defas \inf \Big\{\mathcal{M}^\pm_\varepsilon(u,Q_1(0)) \colon \varepsilon \in (0,1],\:\: u \in \mathcal{D}\Big\},    
\end{equation}
with 
\begin{equation*}
\begin{split}
    \mathcal{D}\defas \Big\{&u \in W_{\mathrm{loc}}^{2,2}(\mathbb{R}^n),\:\:  u=  1\:\: \text{near}\:\:  {x\cdot e_n} =  \frac{1}{2}, \:\: u=  -1\: \text{near}\:\:  {x\cdot e_n} =  -\frac{1}{2}, \\ &u(x+e_i)=u(x)\:\: \text{for all}\:\: x\in \mathbb{R}^n,\:\: i=1,...,n-1 \Big\}.
\end{split}    
\end{equation*}
Here with ``$u=1$ near $x \cdot e_n = \frac{1}{2}$" and ``$u=-1$ near $x \cdot e_n = -\frac{1}{2}$" we mean that there exists a neighborhood $N_{+}$ of $\{x \in Q_1(0) \colon x \cdot e_n = \frac{1}{2}\}$ and a neighborhood $N_{-}$ of $\{x \in Q_1(0) \colon x \cdot e_n = -\frac{1}{2}\}$ such that $u(x)=1$ for every $x \in N_{+}$ and $u(x)=-1$ for every $x \in N_{-}$.
In general, given a set $U$,  with ``$u=\Tilde{u}$ near $U$"  we will always intend that there exists a neighborhood $N$ of $ U$ such that $u=\Tilde{u}$ on $N$. We point out that the infimum problem \eqref{chermisidensity} is taken also over $0<\varepsilon\leq 1$ and not only over a space of functions. 
Notice that by virtue of $\NNN(f_2)$--$\NNN(f_3)\EEE$ it holds
\begin{equation}
\label{growthcondition}
c_1\mathcal{M}^-_\varepsilon(u,A)\leq \mathcal{E}_\varepsilon(u,A)\leq c_2 \mathcal{M}^+_\varepsilon(u,A).
\end{equation}
Thus it is natural to expect that $\mathcal{M}^-_{\varepsilon}$ and $\mathcal{M}^+_{\varepsilon}$ play a role in the analysis of the asymptotic behavior of $\mathcal{E}_{\varepsilon}$.    
If $\varepsilon=1$, we will write $\mathcal{E}$ instead of $\mathcal{E}_1$ and $\mathcal{M}^{\pm}$ instead of $\mathcal{M}_1^{\pm}$. 
Let $\eta \in C^2(\mathbb{R})$ be such that $\eta(t)=1$ if $t>\frac{1}{2}$ and $\eta(t)=-1$ if $t<-\frac{1}{2}$. For every $y \in \mathbb{R}^n$ we define ${u}^{\nu}_{x,\varepsilon}(y)\defas \eta\Big(\frac{(y-x)\cdot \nu}{\varepsilon}\Big)$. We point out, for arguments that will be used in the following, that $u^{\nu}_{x,\varepsilon}$ is \NNN constant \EEE along all directions orthogonal to $\nu$. To simplify the notation, if $\varepsilon=1$ we write ${u}^{\nu}_{x}$ instead of ${u}^{\nu}_{x,1}$. 
For every $x \in \mathbb{R}^n$ and $\nu \in \mathbb{S}^{n-1}$, we consider the oriented piecewise constant function $\overline{u}^\nu_x \colon \mathbb{R}^n \to \{-1,1\}$ defined as 
\begin{equation}
\label{jumpfunction}
    \overline{u}^\nu_x(y)=\begin{cases}
        \:\:\:1 \:\:  \text{if}\:\: (y-x)\cdot \nu \geq 0 \\
        -1 \:\: \text{if}\:\: (y-x)\cdot \nu < 0.
    \end{cases}
\end{equation}
\NNN Notice that the function ${u}^{\nu}_{x,\varepsilon}$ can be seen as regularizations of the jump function $\overline{u}^\nu_x$. 
We highlight to the reader that, frequently in the literature, the notation for $u^\nu_x$ and $\overline{u}_x^\nu$ are inverted. The choice of switching these two notations is motivated by the fact that we use the regularizations of \eqref{jumpfunction} more frequently than the function itself. \EEE \\
Given $A \in \mathcal{A}$ and $\tilde{u}\in W^{2,2}(A)$, we define
\begin{equation*}
    \mathcal{S}(\Tilde{u},A)\defas \{u \in W^{2,2}(A)\colon u=\Tilde{u}\:\: \text{near}\:\: \partial A\}.
\end{equation*}
\NNN Furthermore, \EEE for every $\varepsilon>0$, $A \in \mathcal{A}$, and $\Tilde{u}\in W^{2,2}(A)$ we define the infimum problem
\begin{equation}
\label{infimumproblem}
m_{\mathcal{E}_{\varepsilon}}(\Tilde{u},A)\defas \inf\{\mathcal{E}_{\varepsilon}(u,A)\colon u \in \mathcal{S}(\Tilde{u},A)\},
\end{equation}
and the following two densities:
\begin{equation}
\label{f´}
f^\prime(x,\nu)\defas \limsup\limits_{\rho \to 0}\liminf\limits_{\varepsilon \to 0}\frac{m_{\mathcal{E}_{\varepsilon}}(u^\nu_{x,\varepsilon},Q^\nu_{\rho}(x))}{\rho^{n-1}}
\end{equation}
and
\begin{equation}
\label{f´´}
f^{\prime \prime}(x,\nu)\defas \limsup\limits_{\rho \to 0}\limsup\limits_{\varepsilon\to 0}\frac{m_{\mathcal{E}_{\varepsilon}}(u^\nu_{x,\varepsilon},Q^\nu_{\rho}(x))}{\rho^{n-1}}.    
\end{equation}
Now, we proceed by reviewing the results in \cite{Chermisi} and highlighting the properties needed for the next \NNN sections. \EEE
Let $A \in \mathcal{A}$ be such that $A=A^\prime \times I$, with $A^\prime \subset \mathbb{R}^{n-1}$ open and bounded, and $I$ be an open bounded interval. Let $R_{\nu} \in O(n)$ as in Section \ref{sec: not} and let $A_{\nu}=R_{\nu}A$. Let $\mathcal{M}^\pm_{\varepsilon}$ be as in \eqref{Chermisi}. By applying Fubini's theorem \NNN and \EEE a change of variable we get 
\begin{equation}
\label{surfacecontrol}
    \mathcal{M}^\pm_\varepsilon({u}_{x,\varepsilon}^\nu,A_{\nu})\leq \int_{A^\prime}\int_{I/\varepsilon}W(\eta(t))+q\vert \eta(t)^\prime\vert^2+\vert \eta^{\prime \prime}(t)\vert^2 \, \mathrm{d}t \,\mathrm{d}x^\prime \leq \int_{A^\prime}\int_{\mathbb{R}}W(\eta(t))+q\vert \eta(t)^\prime\vert^2+\vert \eta^{\prime \prime}(t)\vert^2 \, \mathrm{d}t \, \mathrm{d}x^\prime
\end{equation}
\NNN Notice that, by virtue of $(W_2)$, the inner integral on the right hand-side of the last equation is finite. \EEE
Hence, by setting
\begin{equation}
\label{Ceta}
 C_{\eta}\defas \int_{\mathbb{R}}W(\eta(t))+q\vert \eta(t)^\prime\vert^2+\vert \eta^{\prime \prime}(t)\vert^2 \, \mathrm{d}t<\infty   
\end{equation}
\EEE
we get 
\begin{equation}
\label{controlsurface}
    \mathcal{M}^\pm_\varepsilon({u}^{\nu}_{x,\varepsilon},A_{\nu})\leq C_{\eta}\mathcal{L}^{n-1}(A^\prime).
\end{equation}
For every $\nu \in \mathbb{S}^{n-1}$ consider the orthonormal basis $e^{\prime}_1,...,\nu$, where $e^{\prime}_i\defas R_{\nu}e_i$ for every $i \in \{1,...,n-1\}$. For every $x\in \mathbb{R}^n$ and $\rho>\varepsilon>0$ consider the infimum problem
\begin{equation}
\begin{split}
\label{FMinfimumproblem}
    \Tilde{m}_{\mathcal{M}^\pm_{\varepsilon}}({u}^\nu_{x,\varepsilon}, Q^\nu_{\rho}(x))\defas \inf \Big\{\mathcal{M}^\pm_{\varepsilon}(u,Q^\nu_{\rho}(x))\colon u \in \mathcal{D}(u^{\nu}_{x,\varepsilon},Q^{\nu}_{\rho}(x))\Big\},
\end{split}
\end{equation}
where 
\begin{equation*}
\begin{split}
    \mathcal{D}(u^\nu_{x,\varepsilon},Q^{\nu}_{\rho}(x))\defas \Big\{&u \in W_{\mathrm{loc}}^{2,2}(\mathbb{R}^n),\:\:  u=  u^{\nu}_{x,\varepsilon}\:\: \text{near}\:\:  {(y-x)\cdot \nu} =  \frac{\rho}{2} \:\: \text{and near}\:\:  {(y-x)\cdot \nu} =  -\frac{\rho}{2}, \\ &u(x+\NNN \rho e^{\prime}_i\EEE)=u(x)\:\: \text{for all}\:\: x\in \mathbb{R}^n,\:\: i=1,...,n-1 \Big\}.
\end{split}    
\end{equation*}
\NNN We point out that $\mathcal{D}(u^\nu_{0,\varepsilon},Q_1^\nu(0))=\mathcal{D}$ if $\nu=e_n$ and $\varepsilon<1$. \EEE In the next \NNN Remark \EEE we state that \eqref{chermisidensity} can be viewed also as limit of the infimum problems defined in \eqref{FMinfimumproblem}.
\begin{remark}
\label{claimsection6}
Let $\sigma^+$ and $\sigma^-$ be the infima defined in \eqref{chermisidensity} respectively for $\mathcal{M}^+_{\varepsilon}$ and $\mathcal{M}^-_{\varepsilon}$.
Then for every $x \in \mathbb{R}^n$, $\nu \in \mathbb{S}^{n-1}$, and $\rho>0$ it holds
\begin{equation}
\lim\limits_{\varepsilon \to 0}\Tilde{m}_{\mathcal{M}^\pm_{\varepsilon}}({u}^\nu_{x,\varepsilon}, Q^\nu_{\rho}(x))=\sigma^\pm \rho^{n-1}.
\end{equation}
The proof is based on the idea that $u^{\nu}_{x,\varepsilon}$ coincides with $\overline{u}^{\nu}_{x}$ in a neighborhood of $\NNN \{ y \in Q_{\rho}^\nu(x): {(y-x)\cdot \nu} =  \pm\frac{\rho}{2}\}\EEE$ if $\varepsilon<\rho$. In particular, with \cite[Remark 4.3]{Chermisi}, in \eqref{chermisidensity} one can replace $Q_1(0)$ with general cubes of side length equal to 1. Then, with a change of variable one can show that, for every $\varepsilon<\rho$, all $u_{\varepsilon}\in \mathcal{D}(u^\nu_{x,\varepsilon},Q^{\nu}_{\rho}(x))$, up to a rescaling by $\rho$, are admissible for the infimum problem \eqref{chermisidensity} on rotated cubes and that they satisfy
\begin{equation*}
    \sigma^\pm \rho^{n-1}\leq \mathcal{M}^{\pm}_{\varepsilon}(u_{\varepsilon},Q^{\nu}_{\rho}(x)),
\end{equation*}
and thus
\begin{equation*}
    \sigma^\pm \rho^{n-1}\leq \liminf\limits_{\varepsilon \to 0}\Tilde{m}_{\mathcal{M}^\pm_{\varepsilon}}({u}^\nu_{x,\varepsilon}, Q^\nu_{\rho}(x)).
\end{equation*} 
In order to obtain the other desired inequality, i.e., 
\begin{equation*}
\limsup\limits_{\varepsilon \to 0}\Tilde{m}_{\mathcal{M}^\pm_{\varepsilon}}({u}^\nu_{x,\varepsilon}, Q^\nu_{\rho}(x))\leq  \sigma^\pm \rho^{n-1},
\end{equation*}
the strategy is to apply
\cite[Proof of Theorem 1.3, Substep 2A]{Chermisi} getting in this way a recovery sequence $(u_{\varepsilon})_{\varepsilon}$ for $\overline{u}^\nu_{x}$ on $Q^\nu_{\rho}(x)$ such that, for every $\varepsilon$ small enough, $u_{\varepsilon}$ is admissible for the infimum problem \eqref{FMinfimumproblem}.
\end{remark}

\EEE

\subsection{Compactness and integral representation} In this subsection, we give our main result on the compactness (with respect to $\Gamma$-convergence) for functionals of type \eqref{defenergies} and the integral representation of their limit. Subsequently we also state the main properties of the limit densities: in particular, they coincide with $f^\prime$ and $f^{\prime \prime}$ along the $\Gamma$-converging subsequence, and they are bounded from above and from below by two quantities depending only on $c_1$, $c_2$, $n$, and $q$.
\begin{theorem}[Compactness and integral representation]
\label{compintrepr}
Let $(\mathcal{E}_\varepsilon)_\varepsilon$ be as in \eqref{defenergies}. Then, there exists a subsequence $(\mathcal{E}_{\varepsilon_j})_j$ of $(\mathcal{E}_\varepsilon)_\varepsilon$ and a functional $\mathcal{E}_{0}(\cdot,A)$ such that $\mathcal{E}_{\varepsilon_j}(\cdot,A)$ $\Gamma$-converges to $\mathcal{E}_0(\cdot,A)$ with respect to the strong $L^2(A)$ topology, for every $A \in \mathcal{A}_1$. Moreover, for every $u \in L_{\mathrm{loc}}^2(\mathbb{R}^n)$ and $A \in \mathcal{A}$, $\mathcal{E}_{0}$ is given by
\begin{equation}
\label{energycompactness}
\mathcal{E}_{0}(u,A)=\begin{cases}
\int_{J_u \cap A}f_{0}\NNN(x,u^+,u^-,\nu_u)\EEE\, \mathrm{d}\mathcal{H}^{n-1}(x)\quad \text{if}\:\: u_{\vert_A}\in BV(A;\{-1,1\})\\
+\infty \quad \quad \quad \quad \quad \quad \quad \quad \quad \quad \quad \quad\quad \:\, \text{otherwise},
\end{cases}
\end{equation}
with $\NNN f_{0}\colon \mathbb{R}^n \times \{-1,1\}\times \{-1,1\}\times \mathbb{S}^{n-1}\to [0,\infty)\EEE$  defined \NNN by \EEE 
\begin{equation}
\label{cellformula}
\NNN f_{0}(x,1,-1,\nu)=\limsup\limits_{\rho \to 0}\frac{1}{\rho^{n-1}}m_{\mathcal{E}_{0}}(\overline{u}_x^\nu,Q_{\rho}^\nu(x))\EEE
\end{equation}
\NNN and
\begin{equation}
\label{cellformula2} 
f_{0}(x,-1,1,\nu)=\limsup\limits_{\rho \to 0}\frac{1}{\rho^{n-1}}m_{\mathcal{E}_{0}}(-\overline{u}_x^\nu,Q_{\rho}^\nu(x)),
\end{equation}
\EEE
for every $x \in \mathbb{R}^n$ and $\nu \in \mathbb{S}^{n-1}$, where 
\begin{equation}
\label{infimumproblem infty}
m_{\mathcal{E}_0}(\overline{u},A)\defas \inf\{\mathcal{E}_{0}(u,A)\colon u \in BV(A;\{-1,1\})\:\: \text{and}\:\: u=\overline{u}\:\: \text{near}\:\: \partial A\}
\end{equation}
for $A \in \mathcal{A}$ and $\overline{u}\in BV(A;\{-1,1\})$. 
\end{theorem}
\NNN
\begin{remark}
\label{simplyremark}
It can be shown that $f_0(x,1,-1,\nu)=f_0(x,-1,1,-\nu)$, for every $x \in \mathbb{R}^n$ and for every $\nu \in \mathbb{S}^{n-1}$. Hence, in order to uniquely determine $f_0$, it is sufficient to study the dependence on $x$ and $\nu$ of $f_0(x,1,-1,\nu)$. For this reason, only with the purpose of simplifying the notation, in the following we assume $f_0(x,\nu)\defas f_0(x,-1,1,\nu)=f_0(x,1,-1,\nu)$. \EEE This is the case for example when $f_{\varepsilon}$ satisfy $f_{\varepsilon}(x,u,\xi,\zeta)=f_{\varepsilon}(x,-u,-\xi,-\zeta)$ for every $\varepsilon>0$. \EEE   
\end{remark}
\begin{proposition}[Properties of $f_0$]
\label{propositioncontroldensity}
Let $(\mathcal{E}_{\varepsilon})_{\varepsilon}$ be as in \eqref{defenergies} and let $(\varepsilon_j)_j$ and $f_0$ be as in Theorem \ref{compintrepr} \NNN under the assumption of Remark \ref{simplyremark}. \EEE Let $f^\prime$ and $f^{\prime \prime}$ be defined as in \eqref{f´} and \eqref{f´´}, with $\varepsilon$ replaced by $\varepsilon_j$, let $\sigma^\pm$ be as in \eqref{chermisidensity}. Then, $f_0$ is Borel measurable and, for every $x \in \mathbb{R}^n$ and $\nu \in \mathbb{S}^{n-1}$, it holds
\begin{equation}
\label{firstequationprop3.4}
    f^\prime(x,\nu)=f^{\prime \prime}(x,\nu)=f_0(x,\nu)
\end{equation}
and
\begin{equation}
\label{secondequationprop3.4}
    c_1\sigma^- \leq f_0(x,\nu)\leq c_2 \sigma^+,
\end{equation}  
\end{proposition}
Let us assume that $f_{\varepsilon}$ has the form  $f_{\varepsilon}(x,u,\xi,\zeta)= f(\frac{x}{\varepsilon},u,\xi,\zeta)$, for some $f \in \mathcal{F}$. Then further information can be added to the previous results. In particular, \NNN if \EEE a certain limit exists and does not depend on $x$ (see \eqref{homlim})\EEE, then actually the $\Gamma$-limit exists, i.e., there exists a unique limit density which does not depend on the extracted $\Gamma$-converging subsequence of energies. In Theorem \ref{stochhomformula} below, we show that this assumption is always verified if $f$ is stationary in the sense of Definition \ref{stationarity}. Specifically, we have the following.
\begin{theorem}[Deterministic homogenization]
\label{deterministichomogenization}
Let $f \in \mathcal{F}$ and let $m_{\mathcal{E}}(u^\nu_{rx},Q_{rx}^\nu(rx))$ be as in \eqref{infimumproblem} with $\varepsilon=1$, $\Tilde{u}=u^\nu_{rx}$ and $A=Q_{rx}^\nu(rx)$. Assume that for every $x \in \mathbb{R}^n$, $\nu \in \mathbb{S}^{n-1}$ the limit
\begin{equation}
\label{homlim}
f_{\mathrm{hom}}(\nu)\defas\lim\limits_{r \to \infty}\frac{m_{\mathcal{E}}(u^\nu_{rx},Q_{rx}^\nu(rx))}{r^{n-1}},
\end{equation}
exists and is independent of $x$. Then, \NNN under the assumption of Remark \ref{simplyremark}, \EEE for every $A \in \mathcal{A}_1$ the functionals $\mathcal{E}_{\varepsilon}$ defined in \eqref{defenergies} with $f_{\varepsilon}(x,u,\xi,\zeta)\defas f(\frac{x}{\varepsilon},u,\xi,\zeta)$, $\Gamma$-converge, with respect to the strong $L^2(A)$ topology, to the functional $\mathcal{E}_{\mathrm{hom}}\colon L_{\mathrm{loc}}^2(\mathbb{R}^n)\times \mathcal{A}\to [0,\infty]$ given by
\begin{equation}
\label{dethomfunctional}
\mathcal{E}_{\mathrm{hom}}(u,A)=\begin{cases}
\int_{J_u \cap A}f_{\mathrm{hom}}(\nu_u(x))\,\mathrm{d}\mathcal{H}^{n-1}(x) & u_{|_A}\in BV(A;\{-1,1\}),\\
+\infty &  \text{otherwise}.
\end{cases}
\end{equation}
\end{theorem}
\begin{proof}
The proof follows the same lines of the one in \cite[Theorem 3.3]{marziani2022gammaconvergence} (see also \cite[Theorem 3.8]{cagnetti2018gammaconvergence}).  Let $f^\prime$, $f^{\prime \prime}$ as in \eqref{f´} and \eqref{f´´} respectively. By virtue of Theorem \ref{compintrepr}, Proposition \ref{propositioncontroldensity}, and the Urysohn's lemma \cite[Proposition 8.3]{DalMaso:93}, it is sufficient to show 
\begin{equation}
\label{fhom=f´=f´´}
    f_{\mathrm{hom}}(\nu)=f^\prime(x,\nu)=f^{\prime \prime}(x,\nu),
\end{equation}
for every $x \in \mathbb{R}^n$ and $\nu \in \mathbb{S}^{n-1}$. Let $x\in \mathbb{R}^n$, $\nu \in \mathbb{S}^{n-1}$ and $\rho>\varepsilon>0$ be fixed. Given $u \in \mathcal{S}(u_{x,\varepsilon}^\nu,Q_{\rho}^\nu(x))$, we can define $u_{\varepsilon}\colon Q_{\rho/\varepsilon}^\nu(x/\varepsilon)\to \mathbb{R}$ as $u_{\varepsilon}(y)\defas u(\varepsilon y)$. Clearly, $u_{\varepsilon}\in \mathcal{S}(u^\nu_x,Q_{\rho/\varepsilon}^\nu(x/\varepsilon))$. With a change of variable we get 
\begin{equation*}
\mathcal{E}_{\varepsilon}(u,Q_{\rho}^\nu(x))=\varepsilon^{n-1}\int_{Q_{\rho/\varepsilon}^\nu(x/\varepsilon)}f(y,u_{\varepsilon},\nabla u_{\varepsilon},\nabla^2 u_{\varepsilon})\,\mathrm{d}y.
\end{equation*}
Hence, by setting $r\defas \frac{\rho}{\varepsilon}$ we obtain
\begin{equation*}
    m_{\mathcal{E}_{\varepsilon}}(u_{x,\varepsilon}^\nu,Q_\rho^\nu(x))=\frac{\rho^{n-1}}{r^{n-1}} m_{\mathcal{E}}(u_{rx/\rho}^\nu,Q_{r}^\nu(rx/\rho)).
\end{equation*}
Passing to the limit as $\varepsilon \to 0$ (i.e., $r \to \infty$), and using \eqref{homlim} with $x/\rho$ instead of $x$, the thesis follows.
\end{proof}
Theorem \ref{compintrepr} and Proposition \ref{propositioncontroldensity} are proved in Section \ref{sec: gamma convergence}.
\subsection{Setting in random environments and stochastic homogenization}
We now move to a probabilistic setting. Here and in the following, we denote with $(\Omega,\mathcal{I},\mathbb{P})$ a complete probability space.
We consider a family of functions $(f_{\varepsilon})_{\varepsilon>0}$, with $f_\varepsilon \colon \Omega \times \mathbb{R}^n \times \mathbb{R}\times \mathbb{R}^n \times \mathbb{R}^{n\times n}\to \mathbb{R}$ defined as 
\begin{equation*}
    f_\varepsilon(\omega,x,u,\xi,\zeta)=f\Big(\omega,\frac{x}{\varepsilon},u,\xi,\zeta\Big),
\end{equation*}
for some function $f \colon \Omega \times \mathbb{R}^n \times \mathbb{R}\times \mathbb{R}^n \times \mathbb{R}^{n\times n}\to \mathbb{R}$ satisfying
\begin{enumerate}
    \item[$(F_1)$] $f$ is $\mathcal{I}\times {\mathcal{B}}^n \times \mathcal{B}\times \mathcal{B}^n \times \mathcal{B}^{n \times n}$ measurable,
    \item[$(F_2)$] for every $\omega \in \Omega$, $f(\omega,\cdot,\cdot,\cdot,\cdot) \in \mathcal{F}$. 
\end{enumerate}
In particular, we call a function $f$ satisfying these properties a \emph{random density}.
The previously mentioned concepts of stationarity and ergodicity for random densities are defined using groups of $\mathbb{P}$-preserving transformations, as we will explain in the following definitions. 
\begin{definition}{(Group of $\mathbb{P}$-preserving transformations)} 
\label{def: pre}
 An $n$-dimensional group of \emph{$\mathbb{P}$-preserving transformations} on $(\Omega,\mathcal{I}, \mathbb{P})$ is a family $(\tau_z)_{z \in \mathbb{Z}^n}$ (resp.\  $(\tau_z)_{z \in \mathbb{R}^n}$) of transformations $\tau_z \colon \Omega \to \Omega$ satisfying the following properties:
\begin{enumerate}
    \item[$(\mathrm{\tau}_1)$] (measurability) $\tau_{z}$ is $\mathcal{I}$-measurable for every $z \in \mathbb{Z}^n$ (resp.\  for every $z \in \mathbb{R}^n)$, 
    \item[$(\NNN\mathrm{\tau}_2\EEE)$] (group property) 
    $(\tau_{z})_{z \in \mathbb{Z}^n}$ (resp.\  $(\tau_{z})_{z \in \mathbb{R}^n}$) is a group with respect to the composition operation, stable with the sum in $\mathbb{R}^n$, that is $\tau_{x} \circ \tau_{y} = \tau_{y} \circ \tau_{x} = \tau_{x+y}$ for every $x,y \in \mathbb{Z}^n $ (resp.\  for every $x,y \in \mathbb{R}^n)$, and $\tau_0 \colon \Omega \to \Omega$ is the identity, 
    \item[$(\NNN \mathrm{\tau}_3\EEE)$] (invariance) $(\tau_{z})_{z \in \mathbb{Z}^n}$ (resp.\  $(\tau_{z})_{z \in \mathbb{R}^n}$) preserves probability, that is $\mathbb{P}(\tau_z(E))=\mathbb{P}(E)$ for every $z \in \mathbb{Z}^n$ (resp.\ for every $z \in \mathbb{R}^n)$ and for every $E \in \mathcal{I}$. 
\end{enumerate}
 In addition, if it  also holds that
 \begin{enumerate}
\item[$(\NNN\mathrm{\tau}_4\EEE)$] given $E \in \mathcal{I}$, $\mathbb{P}(E \triangle \tau_z(E))=0$ for every $z\in \mathbb{Z}^n$ (resp.\  for every $z\in \mathbb{R}^n$) implies $\mathbb{P}(E)=0$ or $\mathbb{P}(E)=1$, 
 \end{enumerate}
 we say that $(\tau_z)_{z \in \mathbb{Z}^n}$ (resp.\  $(\tau_z)_{z \in \mathbb{R}^n})$ is \emph{ergodic}.
\end{definition}
We emphasize that $\NNN (\mathrm{\tau}_4)\EEE$   will be needed only to show that the homogenized $\Gamma$-limit is deterministic whereas all other results in this paper hold also without this property.    
\begin{definition}(Stationarity)
\label{stationarity}
We say that a random density $f$ is \emph{stationary} with respect to a ($n$-dimensional) group of $\mathbb{P}$-preserving transformations $(\tau_z)_{z \in \mathbb{Z}^{n}}$ (resp.\  $(\tau_z)_{z \in \mathbb{R}^{n}}$) on $(\Omega,\mathcal{I},\mathbb{P})$ if      
\begin{equation*}
\begin{split}
f(\tau_{z}\omega,x,u,\xi,\zeta)=f(\omega,x+z,u,\xi,\zeta), 
\end{split}
\end{equation*}
 for every $\omega \in \Omega$, $x \in \mathbb{R}^n$, $u \in \mathbb{R}$, $\xi \in \mathbb{R}^n$, $\zeta \in \mathbb{R}^{n \times n}$, and $z \in \mathbb{Z}^n$ (resp.\  $z \in\mathbb{R}^n$).
\end{definition}
Let $f$ be a random density and let $\mathcal{E}_{\varepsilon}\colon \Omega \times L_{\mathrm{loc}}^2(\mathbb{R}^n)\times \mathcal{A}\to \NNN \mathbb{R}\EEE$ be the corresponding energies defined by
\begin{equation}
    \label{defrandomenergies}    
    \mathcal{E}_\varepsilon[\omega](u,A)=\begin{cases}
        \frac{1}{\varepsilon}\int_A f\Big(\omega,\frac{x}{\varepsilon},u,\varepsilon\nabla u,\varepsilon^2 \nabla^2u \Big)\, \mathrm{d}x \quad \text{if}\:\: u_{\vert_A}\in W^{2,2}(A)\\
        +\infty \quad \quad \quad \quad \quad\quad \quad \quad \quad \quad \quad \quad\quad \text{otherwise}.
    \end{cases}
    \end{equation}
If $\varepsilon=1$ we again write $\mathcal{E}$ instead of $\mathcal{E}_1$. For every $\omega \in \Omega$, $A\in \mathcal{A}$, and $\Tilde{u}\in W^{2,2}(A)$, we consider the stochastic infimum problem defined as 
\begin{equation}
\label{randominfimumproblem}
m_{\mathcal{E}[\omega]}(\Tilde{u},A)\defas\inf\{\mathcal{E}[\omega](u,A)\colon u \in \mathcal{S}(\Tilde{u},A)\}.
\end{equation}
Notice that \eqref{randominfimumproblem} is obtained from \eqref{infimumproblem} by replacing $\mathcal{E}_{\varepsilon}$ with $\mathcal{E}[\omega]$.
The following theorem on existence of limits of asymptotic cell formulas will be used in the statement of the stochastic homogenization result.
\begin{theorem}[Homogenization formula]
\label{stochhomformula}
 Let $f$ be a stationary random density with respect to a group $\{\tau_z\}_{z \in \mathbb{Z}^n}$ (resp.\  $\{\tau_z\}_{z \in \mathbb{R}^n}$) of $\mathbb{P}$-preserving transformations on $(\Omega,\mathcal{I}, \mathbb{P})$, and let $(\mathcal{E}_{\varepsilon})_{\varepsilon}$ be the corresponding random energies defined according to \eqref{defrandomenergies}. In addition, for every $\omega \in \Omega$ let $m_{\mathcal{E}[\omega]}$ be as in Definition \ref{randominfimumproblem}.
 Then, there exists an event $\Omega' \in \mathcal{I}$, with $\mathbb{P}(\Omega')=1$, and a function $f_{\mathrm{hom}} \colon \Omega \times  \mathbb{S}^{n-1} \to [0,\infty]$, independent of $x$, such that for every $\omega \in \Omega'$, $x \in \mathbb{R}^n$,  $\nu \in \mathbb{S}^{n-1}$ and for every function $r\colon (0,\infty)\to (0,\infty)$, with $r(t)\geq t$ for every $t>0$, it holds
\begin{equation}
\label{eqhom}
f_{\mathrm{hom}}(\omega,\nu)=\lim\limits_{t \to \infty}\frac{ m_{\mathcal{E}[\omega]}(u^\nu_{tx},Q_{r(t)}^{\nu}(tx))}{{r(t)}^{n-1}}.
\end{equation}
 Moreover, if $(\tau_z)_{z \in \mathbb{Z}^{n}}$ (resp.\  $(\tau_z)_{z \in \mathbb{R}^{n}}$) is ergodic, then $f_{\mathrm{hom}}$ does not depend on $\omega$, and we have 
\begin{equation}
\label{eqerg}
f_{\mathrm{hom}}(\nu)=\lim\limits_{t \to \infty} \frac{1}{{r(t)}^{n-1}}\int_{\Omega}{  m_{\mathcal{E}[\omega]}(u^\nu_{0},Q_{r(t)}^{\nu}(0))\,\mathrm{d}\mathbb{P}(\omega)}.
\end{equation}
\end{theorem}
We now come to the main statement. 
\begin{theorem}[Stochastic homogenization]
\label{stochastichomogenization}
Let $f$ be a stationary random density with respect to a group $(\tau_z)_{z \in \mathbb{Z}^n}$ (resp.\  $(\tau_z)_{z \in \mathbb{R}^n}$) of $\mathbb{P}$-preserving transformations on $(\Omega,\mathcal{I},\mathbb{P})$. Let $\mathcal{E}_{\varepsilon}$ be as in \eqref{defrandomenergies}, let $\Omega' \in \mathcal{I}$ (with $\mathbb{P}(\Omega')=1$), $f_{\mathrm{hom}}$ as in Theorem \ref{stochhomformula}, and let $\mathcal{E}_{\mathrm{hom}}\colon \Omega \times L_{\mathrm{loc}}^2(\mathbb{R}^n)\times \mathcal{A} \to [0,\infty]$ be the surface functional defined by
\begin{equation}
\label{homogenisedenergy}
\mathcal{E}_{\mathrm{hom}}[\omega](u,A)=\begin{cases}
\int_{J_u \cap A}f_{\mathrm{hom}}(\omega,\nu_u(x))\,\mathrm{d}\mathcal{H}^{n-1}(x) & u_{|_A}\in BV(A;\{-1,1\}),\\
+\infty &  \text{otherwise}
\end{cases}
\end{equation}
for $\omega \in \Omega$ and $A \in \mathcal{A}$. 
Then, \NNN under the assumption of Remark \ref{simplyremark},\EEE 
\begin{equation}\label{gamma}
\mathcal{E}_{\varepsilon}[\omega](\cdot, A)\:\: \Gamma\text{-converge}\:\:\text{to}\:\: \mathcal{E}_{\mathrm{hom}}\NNN [\omega]\EEE(\cdot,A) \ \ \ \ \text{with respect to the strong $L^2(A)$ topology}, 
\end{equation}
for every $\omega \in \Omega'$ and every $A  \in  \mathcal{A}_1 $.
Furthermore, if $(\tau_z)_{z \in \mathbb{Z}^n}$ (resp.\  $(\tau_z)_{z \in \mathbb{R}^n}$) is ergodic, then  $\mathcal{E}_{\mathrm{hom}}$ is a deterministic functional,  i.e.\  it does not depend on $\omega$.
\end{theorem}
The results in this subsection are proved in Section \ref{sec: stochastichomogenization}.
\section{$\Gamma$-convergence}
\label{sec: gamma convergence}
For convenience of the reader, we divide this section in \NNN two \EEE parts. The \NNN first \EEE subsection is devoted to the statement and the proof of the fundamental estimate (Proposition \ref{fundamentalestimate}) and its application to Theorem \ref{compintrepr}. 
In the \NNN second \EEE subsection, we study and characterize the surface integrand of the $\Gamma$-limit, establishing its relationship with the original densities $(f_{\varepsilon})_{\varepsilon}$. Finally we prove Proposition \ref{propositioncontroldensity}.

\subsection{$\Gamma$-convergence and integral representation}
We start this subsection by proving a fundamental estimate for $(\mathcal{E}_{\varepsilon})_{\varepsilon}$. Fundamental estimates have been proved to be standard tools in the field of homogenization. For instance, together with abstract $\Gamma$-convergence compactness results, they ensure the existence of $\Gamma$-converging subsequences independent of the chosen subdomain. We refer to the proof of Theorem \ref{compintrepr} for further details. \\
\NNN
In the following, for every $u \in L_{\mathrm{loc}}^2(\mathbb{R}^n)$ and $A \in \mathcal{A}$, we set
\begin{equation}
\label{gammaliminfsup}
\mathcal{E}^{\prime}_0(u,A)\defas \Gamma-\liminf\limits_{\varepsilon \to 0}\mathcal{E}_{\varepsilon}(u,A) \:\: \text{and}\:\: \mathcal{E}^{\prime \prime}_0(u,A)\defas \Gamma-\limsup\limits_{\varepsilon \to 0}\mathcal{E}_{\varepsilon}(u,A),  
\end{equation}
meant with respect to the strong  $L^2(A)$ topology. We notice that, because of \eqref{growthcondition} and \cite[Theorem 1.3]{Chermisi}, for every $A \in \mathcal{A}_1$ it holds
\begin{equation}
\label{eq4.30}
c_1 \mathcal{M}^-_0(u,A)\leq \mathcal{E}^\prime_0(u,A)\leq \mathcal{E}_0^{\prime \prime}(u,A)\leq c_2 \mathcal{M}^+_0(u,A),
\end{equation}
where $\mathcal{M}^{\pm}_0$ are the functionals defined in \eqref{FonsecaMantegazzalimit}. 
Finally, we recall that given a set function $F \colon \mathcal{A}\to \mathbb{R}^+$, its inner regular envelope, denoted by $F_{-}$, is defined as 
\begin{equation*}
    F_{-}(A)\defas \sup\{ F(U)\colon U \subset \subset A, \:\: U \in \mathcal{A}\},
\end{equation*}
for every $A \in \mathcal{A}$.
\EEE
\NNN
Before proceeding with the fundamental estimate, we make a remark on \cite[Theorem 1.2]{Chermisi}.
\begin{remark}
\label{remarkonTheorem1.2}
Let $A^\prime \in \mathcal{A}_1$. Because of \cite[Theorem 1.2]{Chermisi} there exists an $\varepsilon_0(A^\prime,q)$ such that
\begin{equation*}    
\mathcal{M}_{\varepsilon}^-(u,A^\prime)\geq 0 \:\: \text{for every}\:\: u \in L_{\mathrm{loc}}^2(\mathbb{R}^n)\:\: \text{whenever}\:\: 0<\varepsilon <\varepsilon_0(A^\prime,q).
\end{equation*}   
 \EEE
\EEE In this remark we want to show that, for a certain family of sets, $\varepsilon_0$ can be chosen uniformly. More precisely, there exists a family of sets $(O_{\delta})_{\delta \in (0,\delta_0)}\subset \mathcal{A}_1$, with $\delta_0>0$, such that $O_{\delta_2}\subset \subset O_{\delta_1}\subset \subset A^\prime$ if $\delta_1<\delta_2$, and an $\overline{\varepsilon} \in (0,\varepsilon_0(A^\prime,q))$ such that \EEE 
\begin{equation}
\label{remarkeps1}
\inf_{\delta \in (0,\delta_0)}\mathcal{M}_{\varepsilon}^-(u,\EEE O_{\delta} \EEE)\geq 0 \:\: \text{for every}\:\: u \in L_{\mathrm{loc}}^2(\mathbb{R}^n)\:\: \text{whenever}\:\:0<\varepsilon<\overline{\varepsilon}. \EEE 
\end{equation}
Sets of the type $O_{\delta}$ will be used as intermediate sets for the proof of Proposition \ref{fundamentalestimate} as they will guarantee the possibility to take only one $\varepsilon$ small enough making the energy positive on all the sets involved in the proof. \EEE \\
For every $x \in \partial A^\prime$ and $r>0$, let us denote with $C_r(x)$ the right circular cylinder centered at $x$ and with height $2r$, radius $r$ and axis parallel to $\nu_{\partial A^\prime}(x)$. \EEE As given by \cite[Section 3.2]{Chermisi}, $\varepsilon_0(q,A^\prime)$ coincides with some $\frac{r}{2}$ small enough to have
\begin{equation*}
 C_r(x)\cap \partial A^\prime \:\: \text{is the graph of a $C^1$-function defined on the basis of the cylinder},   
\end{equation*}
and
\begin{equation*}
\sup_{x \in \partial A^\prime}\mathcal{H}^{n-1}(\{ \xi \in \mathbb{S}^{n-1}\colon \xi \in T_{\partial A^\prime}(y), \:\: y \in C_r(x)\cap \partial A^\prime\})<\eta,    
\end{equation*}
where $\eta>0$ depends only on $q$ and $n$, and $T_{\partial A^\prime}(y)$ denotes the tangent space to $\partial A^\prime$ in $y\in \partial A^\prime$.\\
Let us fix $O \subset \subset A^\prime$ with $O$ open, bounded and with $C^2$-boundary. \EEE For every $\delta>0$ we define 
\begin{equation*}
\EEE O_{\delta}\defas \{ x \in O\colon \mathrm{dist}(x,\partial O)>\delta\}. 
\end{equation*}
Notice that the sets $O_{\delta}$ approximate $O$, as $\delta$ becomes smaller and smaller, in the sense of \cite[Theorem 4.1]{Doktor1976} and \cite[Proposition 3, (12)-(13) and Lemma 2]{HenaoMora} (see also \cite[Theorem 16.25.2]{MR0350769}). In particular, we notice that there exists a $\delta_0=\delta_0(O)>0$ small enough such that, for every $\delta<\delta_0$, it holds $O_{\delta}\in \mathcal{A}_1$ and $T_{\partial O_{\delta}}( x-\delta\nu_{\partial O}(x))=T_{\partial O}(x)$ for every $x \in \partial O$. In addition, up to taking $\delta_0>0$ even smaller but still depending on $O$, we can find an $0< \overline{r}=\overline{r}(q,O)< r$ depending only on $q$ and $O$, such that \EEE
\begin{equation*}
 C_{\overline{r}\EEE}(x)\cap \partial O_{\delta} \:\: \text{is the graph of a $C^1$-function defined on the basis of the cylinder, for every}\:\: \delta \in (0,\delta_0),   
\end{equation*}
and
\begin{equation*}
\sup_{x \in \partial \EEE O_{\delta}\EEE}\mathcal{H}^{n-1}(\{ \xi \in \mathbb{S}^{n-1}\colon \xi \in T_{\partial \EEE O_{\delta}\EEE}(y), \:\: y \in C_{\overline{r}\EEE}(x)\cap \partial \EEE O_{\delta}\EEE\})<\eta.
\end{equation*}
In particular, by setting $\overline{\varepsilon}=\frac{\overline{r}}{2}$ \EEE we have \eqref{remarkeps1}.
Similarly, if $A \in \mathcal{A}_1$ with $A \subset \subset A^\prime$ and $\mathrm{dist}(A,\mathbb{R}^n \setminus A^\prime)=d>0$, it holds 
\begin{equation}
\label{remarkeps2}
\mathcal{M}_{\varepsilon}^-(u,A^\prime \setminus \overline{A})\geq 0 \:\: \text{for every}\:\: u \in L_{\mathrm{loc}}^2(\mathbb{R}^n)\:\: \text{whenever}\:\:0<\varepsilon<\min\Big\{\varepsilon_0(A,q),\varepsilon_0(A^\prime,q),\frac{d}{8}\Big\}.
\end{equation}
\end{remark}

\begin{proposition}[Fundamental estimate]
\label{fundamentalestimate}
Let $(\mathcal{E}_\varepsilon)_\varepsilon$ be as in \eqref{defenergies}. Let $A$, $A^{\prime} \in \mathcal{A}_1$ with $A \subset \subset A^{\prime}$. Let $B \in \mathcal{A}_1$ be such that $A^{\prime}\setminus \overline{A} \subset B$. Let $S \defas A^\prime \setminus \overline{A}$ and let $(u_{\varepsilon})_{\varepsilon} \subset W^{2,2}(A^\prime)$ and $(v_{\varepsilon})_{\varepsilon}\subset W^{2,2}(B)$. Then there exists a $\varepsilon_0>0$, depending on $q$, $A$, $A^{\prime}$ and $B$, and functions $ w_{\varepsilon} \in W^{2,2}(A \cup B)$, such that for every $\varepsilon\leq \varepsilon_0$ it holds
\begin{equation*}
\NNN w_\varepsilon=u_{\varepsilon} \:\: \text{a.e in}\:\: A ,\:\: w_\varepsilon=v_{\varepsilon} \:\: \text{a.e in}\:\: B \setminus \overline{A}^\prime,     
\end{equation*}
and
\begin{equation}
\label{fundamentalestimate1}
\NNN \mathcal{E}_{ \varepsilon}(w_\varepsilon ,A \cup B)\leq 
\mathcal{E}_{ \varepsilon}( u_{\varepsilon},A^\prime)+\mathcal{E}_{\varepsilon}( v_{\varepsilon},B)+\omega_\varepsilon(u_{\varepsilon}, v_{\varepsilon},A,A^\prime),
\end{equation}
where $\NNN \omega_\varepsilon \colon (L_{\mathrm{loc}}^2(\mathbb{R}^n))^2 \times \mathcal{A}_1^2 \to [0,\infty)\EEE$ satisfies, for every measurable set $\Tilde{S} \subset S$ such that $\NNN u_{\varepsilon}=v_{\varepsilon}\EEE$ on $\Tilde{S}$, 
\begin{align}
\label{omegaestimate}
\NNN\omega_\varepsilon( u_{\varepsilon} , v_{\varepsilon},A,A^\prime)\EEE&\leq C\mathcal{L}^{n}(S\setminus \Tilde{S})+\NNN C\int_{S \setminus \tilde{S}}\vert u_{\varepsilon}-v_{\varepsilon}\vert^2 \, \mathrm{d}x+ C\varepsilon (\mathcal{M}^-_{ \varepsilon}(u_{\varepsilon},S)+\mathcal{M}^-_{ \varepsilon}( v_{\varepsilon},S)), \EEE
\end{align}
 for some constant $C=C(c_1,c_2,W,q)>0$, and
 \begin{equation}
 \label{limomega}
\NNN \lim\limits_{\varepsilon \to 0}\omega_\varepsilon(u_{\varepsilon},v_{\varepsilon},A,A^\prime)=0,
\end{equation}
whenever $\sup_{\varepsilon}\mathcal{E}_{\varepsilon}(u_{\varepsilon},S)+\mathcal{E}_{\varepsilon}(v_{\varepsilon},S)<\infty$ and $u_{\varepsilon}-v_{\varepsilon}$ converges to zero in $L^2(S)$.
\end{proposition}
\begin{proof}
\emph{Step 1}: Let $A$, $A^\prime$, $B \in \mathcal{A}_1$ and $\NNN (u_{\varepsilon})_{\varepsilon}\EEE$,  $\NNN (v_{\varepsilon})_{\varepsilon}\EEE$ be as in the statement. In this step we show a preliminary inequality. Let $\NNN \varepsilon>0\EEE$ \EEE and $O$ be open, bounded, with $C^2$ boundary and such that $A\subset \subset O \subset \subset A^\prime$. \EEE Let $d=\mathrm{dist}(A;\mathbb{R}^n \setminus \NNN O \EEE)$ \EEE and 
\EEE
$\NNN N_{\varepsilon}=\lfloor 2\varepsilon^{-1} \rfloor \EEE$. \EEE Let us define intermediate sets $\EEE O_1,...,O_{N_{\varepsilon}+1} \EEE \in \mathcal{A}_1$ with
\begin{equation*}
   \EEE A \subset \subset O_1 \subset \subset ... \subset \subset O_{N_{
   \varepsilon
   }+1}\subset \subset O \subset \subset A^\prime ,
\end{equation*}
\NNN and such that $\EEE O_i\defas O_{\delta_i}$, for some $\delta_i \in (0,\delta_0)$ where 
\begin{equation*}
    \EEE O_{\delta_i}=\{ x \in O \colon \mathrm{dist}(x,\partial O)>\delta_i\},
\end{equation*}
and $\EEE \delta_0=\delta_0(O) \EEE$ is the parameter defined in Remark \ref{remarkonTheorem1.2}. 
\NNN  Notice that for every $i \in \{1,...,N_{\varepsilon}+1\}$, we have $\EEE O_i \cap B \in \mathcal{A}_1$. \EEE
Furthermore, \EEE up to assuming without loss of generality $\delta_0<d$, \EEE we choose $\EEE (O_i)_i$ to satisfy also
\begin{equation*}
\mathrm{dist}(\EEE O_i,\mathbb{R}^n \setminus O_{i+1})\geq \EEE \frac{\delta_0}{\NNN N_{\varepsilon}+2\EEE},\:\: \NNN \text{for every}\:\: i=1,...,N_{\varepsilon}.    
\end{equation*}
For each $i=1,...,\NNN N_{\varepsilon}\EEE$ let $\varphi_i$ be a \NNN smooth function with compact support such that $\varphi_i$ takes values in $[0,1]$, $\varphi_i=1$ on $\EEE O_i$, $\EEE \mathrm{supp}\, \varphi_i \subset O_{i+1}$, \EEE
\begin{equation}
\label{control1cutoff}
\NNN L_{\varepsilon}\EEE \defas \max_{i=1,...,\NNN N_{\varepsilon}\EEE}\vert \nabla \varphi_i \vert \leq \frac{2(\NNN N_{\varepsilon}\EEE+2)}{\EEE \delta_0}  
\end{equation}
and 
\begin{equation}
\label{control2cutoff}
\NNN M_{\varepsilon}\EEE \defas  \max_{i=1,...,\NNN N_{\varepsilon}\EEE}\vert \nabla^2 \varphi_i \vert \leq 4\Big(\frac{\NNN N_{\varepsilon}\EEE+2}{\EEE \delta_0}\Big)^2 .
\end{equation}
For $i=1,...,\NNN N_{\varepsilon}\EEE$ we define $\NNN w_{\varepsilon}^i \in W^{2,2}(A \cup B)\EEE$ by
\begin{equation}
\label{interpolation}
    \NNN w_{\varepsilon}^i \defas \varphi_i u_{\varepsilon} + (1-\varphi_i)v_{\varepsilon}.\EEE
\end{equation} 
Because of \eqref{growthcondition}, \cite[Theorem 1.2]{Chermisi} and \eqref{remarkeps1}--\eqref{remarkeps2} in Remark \ref{remarkonTheorem1.2}, \NNN there exists a $\varepsilon_{0}=\varepsilon_{0}(q,A,A^{\prime},B,O)$ such that for every $\varepsilon \in (0,\varepsilon_0]$ it holds
\begin{equation}
\label{correctionfundamentalestimate}
\mathcal{E}_{\varepsilon}(z,U)\geq 0\:\: \text{for every}\:\: z \in W^{2,2}(U)\:\: \text{with}\:\:U\in \{A^{\prime}\setminus \EEE \overline{O}_{i},\:\: O_i \cap B \colon \:i=1,...,N_{\varepsilon}\}.
\end{equation}
\EEE
We set $S_i \defas \EEE O_{i+1}\setminus \overline{O}_i$. Thus, we obtain 
\begin{align}
\label{eq3.5}
&\NNN\mathcal{E}_{\varepsilon}(w_{\varepsilon}^i, A \cup B)=\mathcal{E}_{\varepsilon}(w_{\varepsilon}^i,\EEE O_{i+1}\cup B \EEE) =\mathcal{E}_{\varepsilon}(u_{\varepsilon},\EEE O_i \EEE)+\mathcal{E}_{\varepsilon}(v_{\varepsilon}, B\setminus \EEE \overline{O}_{i+1}\EEE)+\mathcal{E}_{\varepsilon}(w_{\varepsilon}^i,S_i)  \nonumber\nonumber\\
& \NNN \leq \mathcal{E}_{\varepsilon}(u_{\varepsilon},\EEE O_i \EEE)+\mathcal{E}_{\varepsilon}(u_{\varepsilon},A^{\prime}\setminus \EEE \overline{O}_i \EEE)+\mathcal{E}_{\varepsilon}(v_{\varepsilon}, \EEE B\setminus \overline{O}_{i+1}\EEE)+\mathcal{E}_{\varepsilon}(v_{\varepsilon},\EEE {O}_{i+1}\cap B \EEE)+\mathcal{E}_{\varepsilon}(w_{\varepsilon}^i,S_i) \nonumber \\
& \NNN =\mathcal{E}_{\varepsilon}(u_{\varepsilon},A^\prime)+\mathcal{E}_{\varepsilon}(v_{\varepsilon},B)+\mathcal{E}_{\varepsilon}(w_{\varepsilon}^i,S_i),
\end{align}
where we used \eqref{correctionfundamentalestimate}, the fact that $\NNN \mathcal{E}_{\varepsilon}(z,\cdot)\EEE$ is an additive set function, and that \NNN for every \NNN $z \in W^{2,2}_{\mathrm{loc}}(\mathbb{R}^n)$ the integral defining $\NNN \mathcal{E}_{\varepsilon}(z,\cdot)\EEE$ vanishes on sets with zero Lebesgue measure. 
\EEE
\\
\emph{Step 2}: In this step we define \NNN the functions $w_\varepsilon$ and $\omega_{\varepsilon}$ of the statement, and we show \eqref{fundamentalestimate1}. We focus on studying $ \mathcal{E}_{\varepsilon}(w_{\varepsilon}^i,S_i)$\EEE. Because of \eqref{growthcondition} it holds
\begin{equation}
\label{eq3.6}
\NNN \mathcal{E}_{\varepsilon}(w_{\varepsilon}^i,S_i)\leq c_2 \int_{S_i}\frac{W(w_{\varepsilon}^i)}{\varepsilon}+q\varepsilon\vert \nabla w_{\varepsilon}^i  \vert^2+\varepsilon^3\vert \nabla^2  w_{\varepsilon}^i \vert^2 \, \mathrm{d}x. 
\end{equation}
Using the convexity of the function $\vert \cdot \vert^2$, and that $\vert \varphi_i \vert \leq 1$, we deduce that there exists a universal constant $C>0$ such that 
\begin{equation}
\label{convexity0}
\NNN \vert \nabla w_{\varepsilon}^i \vert^2 \leq C(\vert \nabla \varphi_i \vert^2 \vert u_{\varepsilon}-v_{\varepsilon}\vert^2+\vert \nabla u_{\varepsilon}\vert^2+\vert \nabla v_{\varepsilon}\vert^2)  \EEE 
\end{equation}
and
\begin{equation}
\label{convexity}
\NNN \vert \nabla^2 w_{\varepsilon}^i \vert^2 \leq C[\vert \nabla^2 \varphi_i \vert^2 \vert u_{\varepsilon}-v_{\varepsilon} \vert^2 + \vert \nabla \varphi_i \vert^2(\vert \nabla u_{\varepsilon}\vert^2+\vert \nabla v_{\varepsilon}\vert^2)+\vert \nabla^2 u_{\varepsilon} \vert^2+\vert \nabla^2 v_{\varepsilon} \vert^2]\EEE.
\end{equation}
Combining \eqref{control1cutoff}, \eqref{control2cutoff}, \eqref{eq3.6}, \eqref{convexity0} and \eqref{convexity}, we get that there exists a constant $\Tilde{C}=\Tilde{C}(q,c_1,c_2)$ such that
\begin{align}
\label{eq3.10}
\NNN \mathcal{E}_{\varepsilon}(w_{\varepsilon}^i,S_i)&\NNN\leq c_2 \int_{S_i}\frac{W(w_{\varepsilon}^i )}{\varepsilon}\, \mathrm{d}x+\Tilde{C}\varepsilon L^2_{\varepsilon}\int_{S_i}\vert u_{\varepsilon}-v_{\varepsilon} \vert^2 \, \mathrm{d}x+\Tilde{C}\varepsilon\int_{S_i}\vert \nabla u_{\varepsilon}\vert^2+\vert \nabla v_{\varepsilon}\vert^2 \, \mathrm{d}x \nonumber \\&\NNN+\Tilde{C}\varepsilon^3 \int_{S_i}L^2_{\varepsilon}(\vert \nabla u_{\varepsilon}\vert^2+\vert \nabla v_{\varepsilon}\vert^2)+M^2_{\varepsilon} \vert u_{\varepsilon}-v_{\varepsilon} \vert^2+\vert \nabla^2 u_{\varepsilon} \vert^2+\vert \nabla^2 v_{\varepsilon} \vert^2\, \mathrm{d}x. \EEE    
\end{align}
\NNN Using that $(S_i)_i$ are disjoint and $S_i\subset S$, we can sum and average in \eqref{eq3.5} and then apply \eqref{eq3.10}. In this way, we find \EEE that there exists an $\NNN i^* \in {1,...,N_{\varepsilon}}\EEE$ such that 
\begin{equation}
\label{i^*}
\NNN\mathcal{E}_{\varepsilon}(w_{\varepsilon}^{i^*}, A \cup B)\leq \mathcal{E}_{\varepsilon}(u_{\varepsilon},A^\prime)+\mathcal{E}_{\varepsilon}(v_{\varepsilon},B)+\frac{1}{N_{\varepsilon}}\sum^{N_{\varepsilon}}_{i=1}\mathcal{E}_{\varepsilon}(w_{\varepsilon}^i,S_i)\EEE
\end{equation}
with
\begin{align}
\label{eq3.12}    
& \NNN\frac{1}{N_{\varepsilon}}\sum^{N_{\varepsilon}}_{i=1}\mathcal{E}_{\varepsilon}(w_{\varepsilon}^i,S_i)\leq \Tilde{C}\Big(\frac{\varepsilon L^2_{{\varepsilon}}}{N_{{\varepsilon}}}+\frac{M^2_{\varepsilon}\varepsilon^3}{N_{\varepsilon}}\Big)\int_{S}\vert u_{\varepsilon}-v_{\varepsilon} \vert^2 \, \mathrm{d}x+\Tilde{C}\Big(\frac{\varepsilon}{N_{\varepsilon}}+\frac{L^2_{\varepsilon} \varepsilon^3}{N_{\varepsilon}}\Big)\int_{S}\vert \nabla u_{\varepsilon}\vert^2+\vert \nabla v_{\varepsilon}\vert^2\, \mathrm{d}x
\\& \NNN +\frac{\tilde{C}\varepsilon^3}{N_{\varepsilon}}\int_{S}\vert \nabla^2 u_{\varepsilon} \vert^2+\vert \nabla^2 v_{\varepsilon} \vert^2\, \mathrm{d}x 
+\frac{c_2}{N_{\varepsilon}}\sum^{N_{\varepsilon}}_{i=1}\int_{S_i}\frac{W(w_{\varepsilon}^i)}{\varepsilon}\, \mathrm{d}x. 
\end{align}
\NNN Now, taking $\NNN \varepsilon\in (0,1)$ \EEE and using that $\NNN \frac{2-\varepsilon}{\varepsilon}\leq N_{\varepsilon}\leq \frac{2}{\varepsilon}$, \eqref{control1cutoff} and \eqref{control2cutoff}, we get that there exist a constant 
$\NNN C=C(q,c_1,c_2, \EEE \delta_0 \EEE)>4\EEE$ such that
\begin{equation}
\label{eq3.13}    
\NNN \frac{\Tilde{C}\varepsilon L^2_{\varepsilon}}{N_{\varepsilon}}\leq \frac{C}{2}, \:\: \frac{\Tilde{C}\varepsilon}{N_{\varepsilon}}\leq C\varepsilon^2, \:\: \frac{\Tilde{C}L^2_{\varepsilon}\varepsilon^3}{N_{\varepsilon}}\leq \frac{C}{2}\varepsilon^2\:\:, \frac{\Tilde{C}M^2_{\varepsilon}\varepsilon^3}{N_{\varepsilon}}\leq \frac{C}{2},\:\: \frac{\tilde{C}c_2}{\varepsilon N_{\varepsilon}}\leq C,
\end{equation}
which implies
\begin{equation*}
\NNN \mathcal{E}_{\varepsilon}(w_{\varepsilon}^{i^*},A \cup B)\leq \mathcal{E}_{\varepsilon}(u_{\varepsilon},A^\prime)+\mathcal{E}_{\varepsilon}(v_{\varepsilon},B)+\omega_{\varepsilon}(u_{\varepsilon},v_{\varepsilon},A,A^\prime), \EEE
\end{equation*}
where
\begin{equation}
\begin{split}
\label{defomega}
\NNN\omega_{\varepsilon}(u_{\varepsilon},v_{\varepsilon},A,A^\prime)=&\NNN C\Big[\int_{S}\vert u_{\varepsilon}-v_{\varepsilon}\vert^2\, \mathrm{d}x+\varepsilon^2\int_{S}\vert \nabla u_{\varepsilon}\vert^2+\vert \nabla v_{\varepsilon}\vert^2\, \mathrm{d}x \\
&\NNN+ \varepsilon^4\int_{S}\vert \nabla^2 u_{\varepsilon} \vert^2+\vert \nabla^2 v_{\varepsilon} \vert^2\, \mathrm{d}x+\sum^{N_{\varepsilon}}_{i=1}\int_{S_i}W(w_{\varepsilon}^i)\, \mathrm{d}x\Big].    
\end{split}
\end{equation}
\NNN Denoting $ w_\varepsilon \defas w_{\varepsilon}^{i^*}$, \eqref{fundamentalestimate1} follows. 
\EEE
\\
\emph{Step 3}: In this step we prove \eqref{omegaestimate}. Let $\Tilde{S}\subset S$ measurable and such that $\NNN u_{\varepsilon}=v_{\varepsilon}\EEE$ a.e in $\Tilde{S}$. Then it clearly holds
\begin{equation*}
\NNN \int_{\Tilde{S}}\vert u_{\varepsilon}-v_{\varepsilon}\vert^2\, \mathrm{d}x=0. 
\end{equation*}
Hence, in view of \eqref{defomega}, in order to prove \eqref{omegaestimate} we just have to show that there exists a constant $C>0$ such that, \NNN for $\varepsilon$ small enough, \EEE it holds 
\NNN
\begin{equation}
\label{claimstep3}
\begin{split}
\varepsilon^2\int_{S}\vert \nabla u_{\varepsilon}\vert^2+\vert \nabla v_{\varepsilon}\vert^2\, \mathrm{d}x
&+ {\varepsilon^4}\int_{S}\vert \nabla^2 u_{\varepsilon} \vert^2+\vert \nabla^2 v_{\varepsilon} \vert^2\, \mathrm{d}x+\sum^{N_{\varepsilon}}_{i=1}\int_{S_i}W(w_{\varepsilon}^i)\, \mathrm{d}x \\   
& \leq C\mathcal{L}^{n}(S\setminus \Tilde{S})+C\varepsilon(\mathcal{M}^-_{ \varepsilon}(u_{ \varepsilon},S)+\mathcal{M}^-_{ \varepsilon}( v_{\varepsilon},S)).    
\end{split}
\end{equation}
\EEE We start by noticing that since $u_{\varepsilon}=v_{\varepsilon}$ on $\tilde{S}\subset S$ and the $S_i$ are disjoint and contained in $S$, we have
\begin{equation}
\label{termwithStilde}
\sum^{N_{\varepsilon}}_{i=1}\int_{S_i \cap \tilde{S}}W(w^i_{\varepsilon})\, \mathrm{d}x = \sum^{N_{\varepsilon}}_{i=1}\int_{S_i \cap \tilde{S}}W(u_{\varepsilon})\, \mathrm{d}x\leq \int_{S}W(u_{\varepsilon})\, \mathrm{d}x.
\end{equation}
We first want to show
\begin{equation}
\label{controlStilde}   
\int_{S}W(u_{\varepsilon})\, \mathrm{d}x \leq C\varepsilon \mathcal{M}^-_{\varepsilon}(u_{\varepsilon},S)
\end{equation}
and
\begin{equation}
\label{eq3.15bis}
\NNN \sum^{N_{\varepsilon}}_{i=1}\int_{S_i \setminus \Tilde{S}}W(w_{\varepsilon}^i)\, \mathrm{d}x \leq  C\mathcal{L}^{n}(S \setminus \Tilde{S})+{C\varepsilon}(\mathcal{M}^-_{\varepsilon}(u_{\varepsilon},S)+\mathcal{M}^-_{\varepsilon}(v_{\varepsilon},S)). \EEE
\end{equation}
\EEE
Because of the convexity of $\vert \cdot \vert$, $(W_4)$, and $\varphi_i \in [0,1]$, it must hold \NNN$\max\{c_0W(u_{\varepsilon})+c_0,c_0W(v_{\varepsilon})+c_0\}\geq W(w_{\varepsilon}^i)$\EEE. Consequently, we have
\begin{equation}
\label{eq3.15}
\begin{split}
\NNN\int_{S_i \setminus \Tilde{S}}W(w_{\varepsilon}^i)\, \mathrm{d}x\leq  2c_0\mathcal{L}^{n}(S_i \setminus \Tilde{S})+c_0\int_{S_i\setminus \Tilde{S}}(W(u_{\varepsilon})+W(v_{\varepsilon}))\, \mathrm{d}x.
\end{split}
\end{equation}
By virtue of \cite[Theorem 1.2]{Chermisi}, since $q<\frac{q^*}{n}$, there exist $\alpha=\alpha(q) \in (0,1)$ and $\NNN \varepsilon^\prime_0=\varepsilon^\prime_0(q,A,A^\prime)>0$ such that 
\begin{equation}
\label{eq3.18bis}
    \int_{S}\frac{\alpha W(z)}{\varepsilon}-q\varepsilon\vert \nabla z \vert^2+\varepsilon^3\alpha\vert \nabla^2 z \vert^2\, \mathrm{d}x>0
\end{equation}
for every $\NNN \varepsilon\leq \varepsilon^\prime_0$ and for every $z \in W^{2,2}(S)$. Then, from \eqref{eq3.18bis}, it easily follows that, up to taking \NNN \NNN$\varepsilon_0\leq\varepsilon^\prime_0$ (from which $\varepsilon\leq \varepsilon^\prime_0$ follows), it holds
\begin{equation}
\label{Stilde2}
    \EEE \int_{S}W(u_{\varepsilon})\, \mathrm{d}x \leq \frac{\varepsilon}{(1-\alpha)}\mathcal{M}^-_{\varepsilon}(u_{\varepsilon},S)\leq C\varepsilon \mathcal{M}^-_{\varepsilon}(u_{\varepsilon},S)
\end{equation}
\EEE and \EEE
\begin{equation}
\label{eq3.19bis}
\begin{split}
\int_{S\setminus \Tilde{S}}(W(u_{\varepsilon})+W(v_{\varepsilon}))\, \mathrm{d}x&\NNN\leq \frac{\varepsilon}{(1-\alpha)}(\mathcal{M}^-_{\varepsilon}(u_{\varepsilon},S)+\mathcal{M}^-_{\varepsilon}(v_{\varepsilon},S))\\
&\NNN\leq C\varepsilon(\mathcal{M}^-_{\varepsilon}(u_{\varepsilon},S)+\mathcal{M}^-_{\varepsilon}(v_{\varepsilon},S))  
\end{split}
\end{equation}
\NNN for some suitable constant $C=C(q)>0$, where we used \eqref{Chermisi}. \EEE
Finally, \EEE \eqref{controlStilde} follows from \eqref{Stilde2} while \eqref{eq3.15bis}  follows from \eqref{eq3.19bis}, by summing over $\NNN i=1,...,N_{\varepsilon}\EEE$ in \eqref{eq3.15}, and recalling that all the $S_i$ are disjoint and contained in $S$.\\
Because of \eqref{eq3.18bis}, it also holds
\begin{equation}
\label{gradienttermgoingtozero}
\NNN  \varepsilon^2\int_{S}\NNN \vert \nabla u_{\varepsilon}\vert^2+\vert\nabla v_{\varepsilon}\vert^2  \,\mathrm{d}x\leq \frac{2\alpha \varepsilon}{q(1-\alpha)}(\mathcal{M}^-_{\varepsilon}(u_{\varepsilon},S)+\mathcal{M}^-_{\varepsilon}(v_{\varepsilon},S)).\EEE
\end{equation}
\NNN Similarly, it can be proved
\begin{equation}
\label{hessiantermgoingtozero}
    {\varepsilon^4}\int_{S} \vert \nabla^2 u_{\varepsilon}\vert^2+\vert \nabla^2 v_{\varepsilon}\vert^2 \, \mathrm{d}x\leq \frac{\varepsilon}{(1-\alpha)}(\mathcal{M}^-_{\varepsilon}(u_{\varepsilon},S)+\mathcal{M}^-_{\varepsilon}(v_{\varepsilon},S)),
\end{equation}
that together with \EEE \eqref{termwithStilde}-\eqref{eq3.15bis} \EEE and \eqref{gradienttermgoingtozero}, implies \eqref{claimstep3}.\EEE 
\\
\emph{Step 4}: In this step we prove \eqref{limomega}. Let $(u_\varepsilon)_\varepsilon$, $(v_\varepsilon)_\varepsilon$ be as in the statement \NNN and let \EEE $w$ \NNN be \EEE their limit in $L^2(S)$. Then we have 
\begin{equation*}
    \lim\limits_{j \to \infty}\int_{S}\NNN\vert u_{\varepsilon}-v_{\varepsilon}\EEE\vert^2\,\mathrm{d}x=0.
\end{equation*}
\NNN In addition, because $\sup_{\varepsilon>0}(\mathcal{E}_{\varepsilon}(u_{\varepsilon},S)+\mathcal{E}_{\varepsilon}(v_{\varepsilon},S))<\infty$, \eqref{growthcondition},  \eqref{gradienttermgoingtozero} and \eqref{hessiantermgoingtozero}, we have
\begin{equation*}
\lim\limits_{\varepsilon \to 0}\varepsilon^2\int_{S}\NNN \vert \nabla u_{\varepsilon}\vert^2+\vert\nabla v_{\varepsilon}\vert^2  \,\mathrm{d}x+{\varepsilon^4}\int_{S} \vert \nabla^2 u_{\varepsilon}\vert^2+\vert \nabla^2 v_{\varepsilon}\vert^2 \, \mathrm{d}x=0.    
\end{equation*}
\EEE
So, recalling \eqref{defomega}, we have just to show that
\begin{equation}
\label{claimfundamentalestimate}
\NNN \lim\limits_{\varepsilon \to 0}\sum^{N_{\varepsilon}}_{i=1}\int_{S_i}W( w^i_{\varepsilon})\, \mathrm{d}x=0. \EEE  
\end{equation}
\EEE
For every $\NNN \varepsilon>0$ consider $\NNN \tilde{w}_\varepsilon \EEE\colon S \to \mathbb{R}$ defined as 
\begin{equation*}
    \NNN\tilde{w}_\varepsilon(x)\EEE=\begin{cases}
        \NNN w^i_{\varepsilon}(x)\EEE\:\: \text{if}\:\: x \in S_i\:\:\text{for}\:\: i \in \{1,...,\NNN N_{\varepsilon\EEE}\}, \\
        w(x)\:\:\:\:\text{otherwise}.
    \end{cases}
\end{equation*}
Then, clearly $\NNN \tilde{w}_\varepsilon\EEE$ converges to $w$ in $L^2(S)$ as $\NNN \varepsilon \to 0$. Notice that, by virtue of \cite[Theorem 1.1]{Chermisi} \EEE applied to the sequence $(u_{\varepsilon})_{\varepsilon}$ (which also converges to $w$)\EEE, we have that $w \in BV(S;\{-1,1\})$ and so, because of \NNN $(W_1)$ and $(W_2)$\EEE, up to subsequences $W(\NNN \tilde{w}_\varepsilon\EEE)\to 0$ a.e. in $S$.
\NNN In addition, using $(W_4)$ and arguing like in the proof \eqref{eq3.15}, we can show that
\begin{equation}
\label{newpotentialinequality}
 \NNN \int_{S \cap \{\vert \tilde{w}_\varepsilon \vert >1\}}W(\tilde{w}_\varepsilon)\, \mathrm{d}x \leq c_0\mathcal{L}^n(S \cap \{\vert \tilde{w}_\varepsilon\vert>1\})+c_0\int_{S}W(u_{\varepsilon})+W(v_{\varepsilon})\,\mathrm{d}x.\EEE
\end{equation}
Finally, arguing like in the proof of \eqref{eq3.19bis}, using \eqref{newpotentialinequality} and the reverse Fatou's lemma, we have
\begin{gather*}
\limsup\limits_{\varepsilon \to 0}\sum^{N_{\varepsilon}}_{i=1}\int_{S_i}W(w^i_{\varepsilon})\, \mathrm{d}x  \leq  \limsup\limits_{\varepsilon \to 0} \int_{S}W(\tilde{w}_\varepsilon)\,\mathrm{d}x\leq \\ \NNN \limsup\limits_{\varepsilon \to 0} \int_{S\cap \{\vert \tilde{w}_\varepsilon \vert \leq  1 \}}W(\tilde{w}_\varepsilon)\,\mathrm{d}x+ C\sup_{\delta>0}(\mathcal{M}^-_\delta(u_\delta,S)+ \mathcal{M}^-_\delta(v_\delta,S))\limsup\limits_{\varepsilon\to 0}\varepsilon\leq  \int_{S}\limsup\limits_{\varepsilon \to 0} W(\Tilde{w}_\varepsilon)\, \mathrm{d}x=0
\end{gather*}
which implies \eqref{claimfundamentalestimate}. \EEE
\end{proof}
\EEE
Now, let us provide a definition that will be used in the proof of Theorem \ref{compintrepr}.
\begin{definition}
\label{setsdensity}
Let $\mathcal{P}^\prime$ be a family of subsets of $\mathbb{R}^n$. A subfamily $\mathcal{P}\subset \mathcal{P}^\prime$ is said to be dense in $\mathcal{P}^\prime$ if for every $A$, $C \in \mathcal{P}^\prime$, with $A \subset\subset C$, there exists a $B \in \mathcal{P}$ such that $A \subset B \subset C$.    
\end{definition}
\EEE
\begin{proof}[Proof of Theorem \ref{compintrepr}] 
The proof follows the localization method of $\Gamma$-convergence. \\
\emph{Step 1:} Since $L_{\mathrm{loc}}^2(\mathbb{R}^n)$ has a countable basis and \cite[Theorem 8.5]{DalMaso:93}, for every $A \in \mathcal{A}_1$, the family of functionals $(\mathcal{E}_{\varepsilon}(\cdot,A))_{\varepsilon}$ admits a $\Gamma$-converging subsequence which in principle depends on $A$.  Let $\mathcal{D}_1$ be a countable dense subset of $\mathcal{A}_1$, intended in the sense of Definition \ref{setsdensity}. 
Indeed, one way to construct $\mathcal{D}_1$ is the following: let $(R_k)_k$ be an enumeration of all the poly-rectangles with rational vertices. It can be easily checked that such family is dense in $\mathcal{A}$. Now, for every $k \in \mathbb{N}$, we select a sequence $(A^k_n)_n \subset \mathcal{A}_1$ such that $A^k_n \subset \subset R_k$ and $A_n^k \uparrow R_k$ as $n \to \infty$. Then, one can set $\mathcal{D}_1 \defas (A^k_n)_{k,n}$.
Up to a diagonal argument, we can find a subsequence $(\varepsilon_j)_j$ such that $\mathcal{E}_{\varepsilon_j}(\cdot,A)$ $\Gamma$-converges with respect to the strong $L^2(A)$ topology for every $A \in \mathcal{D}_1$. In particular, along the subsequence $(\varepsilon_j)_j$, it holds $\mathcal{E}^{\prime}_0(u,A)=\mathcal{E}^{\prime \prime}_0(u,A)$ for every $u \in L_{\mathrm{loc}}^2(\mathbb{R}^n)$ and $A \in \mathcal{D}_1$. Furthermore, we recall that, because of \cite[Proposition 6.8]{DalMaso:93}, the functionals $\mathcal{E}^{\prime}(\cdot,A)$ and $\mathcal{E}^{\prime \prime}(\cdot,A)$ are lower semicontinuous with respect to the strong $L^2(A)$  topology for every $A \in \mathcal{A}$. 
Combining this with Lemma \ref{increasingsetfunctionlemma} in the appendix and \cite[Theorem 15.18 $(b)$ and $(e)$]{DalMaso:93}, we get that the inner regular envelopes of $\mathcal{E}_0^\prime$ and $\mathcal{E}^{\prime \prime}_0$ coincide on $L_{\mathrm{loc}}^2(\mathbb{R}^n)\times \mathcal{A}_1$. Thus, we define a functional $\mathcal{E}_0 \colon L^2_{\mathrm{loc}}(\mathbb{R}^n)\times \mathcal{A}_1 \to [0,\infty]$ by setting it to be equal to the inner regular envelope of $\mathcal{E}^{\prime}_{0}\colon L_{\mathrm{loc}}^2(\mathbb{R}^n)\times \mathcal{A}_1 \to [0,\infty]$ i.e. 
\begin{equation}
    \mathcal{E}_0(u,A)\defas \sup\{\mathcal{E}^\prime_0(u,U)\colon U \in \mathcal{A}_1 \:\: \text{and}\:\: U \subset \subset A\}\:\: \text{for every}\:\: u \in L_{\mathrm{loc}}^2(\mathbb{R}^n)\:\: \text{and}\:\: A \in \mathcal{A}_1.
\end{equation} 
\emph{Step 2:} In this step we consider $\mathcal{E}_0^{\prime}$ and $\mathcal{E}^{\prime \prime}_0$ evaluated only along the sequence $(\varepsilon_j)_j$ of \emph{Step 1}. Because of \emph{Step 1} and Lemma \ref{increasingsetfunctionlemma}, we have $\mathcal{E}_0(u,U)\leq \mathcal{E}^\prime_0(u,U)\leq \mathcal{E}^{\prime \prime}_0(u,U)$ for every $u \in L_{\mathrm{loc}}^2(\mathbb{R}^n)$ and $U \in \mathcal{A}_1$. So to prove that $\mathcal{E}_{\varepsilon_j}(\cdot,U)$ $\Gamma$-converges to $\mathcal{E}_0(\cdot,U)$ for every $U \in \mathcal{A}_1$, it is sufficient to prove that $\mathcal{E}_0(u,U)\geq \mathcal{E}_0^{\prime \prime}(u,U)$ for every $U \in \mathcal{A}_1$.
\NNN Notice that, by virtue of \eqref{FonsecaMantegazzalimit} and \eqref{eq4.30}, and since $\mathcal{M}^-_0$ is inner regular, it holds $\mathcal{E}_0(u,U)\geq c_1\mathcal{M}_0^-(u,U)=\infty$ if $u \notin BV(U;\{-1,1\})$\EEE. Hence, we only have to prove
\begin{equation}
\label{eq4.32}
\mathcal{E}_0(u,U)\geq \mathcal{E}_0^{\prime \prime}(u,U)\:\: \text{if}\:\: u \in BV(U;\{-1,1\}).   
\end{equation}
 Let \NNN $A$, $A^{\prime}$, $B \in \mathcal{A}_1$, $A \subset \subset A^{\prime}$ \EEE and $A^{\prime}\setminus \overline{A} \subset B$.
Notice that using  Lemma \ref{increasingsetfunctionlemma}, the fundamental estimate in Proposition \ref{fundamentalestimate}, \eqref{eq4.30}, and arguing like in the proof of \cite[Proposition 18.3]{DalMaso:93} we can show
\begin{equation}
\label{eq4.31}
\mathcal{E}^{\prime \prime}_0(u,A \cup B)\leq \mathcal{E}^{\prime \prime}_0(u,A^\prime)+c_2\mathcal{M}^{+}_0(u,B).    
\end{equation}
Indeed, let $(u_j)_j \subset W^{2,2}(A^{\prime})$ and $(v_j)_j \subset W^{2,2}(B)$ such that $u_j \to u$ in $L^2(A^{\prime})$, $v_j \to u$ in $L^2(B)$, 
\begin{equation*}
    \lim\limits_{j \to \infty}\mathcal{E}_{\varepsilon_j}(u_j,A^{\prime})=\mathcal{E}^{\prime \prime}_0(u,A^{\prime})\:\: \text{and}\:\:  \lim\limits_{j \to \infty}\mathcal{M}^{+}_{\varepsilon_j}(v_j,B)=\mathcal{M}^+_0(u,B).
\end{equation*}
Thanks to \eqref{fundamentalestimate1}, \eqref{limomega} in Proposition \ref{fundamentalestimate} and $\NNN(f_2)\EEE$, \NNN up to taking $\varepsilon_j$ smaller than some $\varepsilon_0$, depending on $A$, $A^{\prime}$ and $B$, we can construct a sequence of functions $(w_{j})_j \subset W^{2,2}(A\cup B)$ such that $w_{j}$ converges to $u$ in $L^2(A\cup B)$ as $j \to \infty$,\EEE
\begin{equation}
\label{bonusinequality1}
\limsup\limits_{j \to \infty}\mathcal{E}_{\varepsilon_{j}}(w_{j},A\cup B)\leq \lim\limits_{j \to \infty}\mathcal{E}_{\varepsilon_{j}}(u_{j},A^{\prime})+c_2 \lim\limits_{j \to \infty}\mathcal{M}^{+}_{\varepsilon_{j}}(v_{j},B)=\mathcal{E}^{\prime \prime}_0(u,A^{\prime})+c_2\mathcal{M}_0^{+}(u,B),   
\end{equation}
and 
\begin{equation*}
\mathcal{E}^{\prime \prime}_0(u,A\cup B)\leq \limsup\limits_{j \to \infty}\mathcal{E}_{\varepsilon_{j}}(w_{j},A\cup B).   
\end{equation*}
The last equation combined with \eqref{bonusinequality1} gives \eqref{eq4.31}. 
\\
 Let $U\in \mathcal{A}_1$, $u \in BV(U;\{-1,1\})$, $\delta>0$ and $K$ be a compact set such that $\mathcal{M}^{+}_0(u, U\setminus K)\leq \delta$. Let us choose $A$, $A^{\prime}\in \NNN\mathcal{A}_1\EEE$ such that $K \subset \subset A \subset \subset A^{\prime}\subset \subset U$.  Thus,  by setting $B \defas U\setminus K$, we get from \eqref{eq4.31}
\begin{equation}
\label{eq4.33}
\mathcal{E}^{\prime \prime}_0(u,U)\leq  \mathcal{E}^{\prime \prime}_0(u,A^{\prime})+c_2\mathcal{M}^{+}_0(u,U \setminus K)\leq \mathcal{E}_0(u,U)+c_2\delta,  
\end{equation}
where we used that $A^{\prime} \subset \subset U$ and the fact that $\mathcal{E}_0$ coincides with the inner regular envelope of $\mathcal{E}^{\prime \prime}_{0}$ on $\mathcal{A}_1$. By sending $\delta \to 0$ in \eqref{eq4.33}, \eqref{eq4.32} follows. 
\\
 \emph{Step 3}: In the previous steps, we proved that there exists a subsequence $(\varepsilon_j)_j$ such that $\mathcal{E}_{\varepsilon_j}(\cdot,A)$ $\Gamma$-converges to a functional $\mathcal{E}_{0}(\cdot,A)$ with respect to the strong $L^2(A)$  topology, for every $A \in \mathcal{A}_1$. In this step we want to show that $\mathcal{E}_0$ admits an integral representation and that the limit integrand is given by \eqref{cellformula} \NNN and \eqref{cellformula2}\EEE. To this purpose it is sufficient to show that the hypotheses of the representation theorem \cite[Theorem 3]{Bouchitt2002AGM} hold.  \\
We point out that in \cite[Theorem 3]{Bouchitt2002AGM} the domain of the energies is a space of functions defined on a fixed open and bounded set. However, one can easily generalize the authors' result to energies defined on $L_{\mathrm{loc}}^1(\mathbb{R}^n)\times \mathcal{A}$.
For every $u \in L_{\mathrm{loc}}^2(\mathbb{R}^n)$ we extend $\mathcal{E}_0(u,\cdot)$ to the whole $\mathcal{A}$ by setting
 \begin{equation}
 \label{extendedinnerregularity}
 \mathcal{E}_0(u,A)\defas\sup\{\mathcal{E}_0(u,U): U\in \mathcal{A}_1 \:\: \text{and}\:\: U \subset \subset A\}\:\: \text{for every}\:\: A \in \mathcal{A}.    
 \end{equation}
Then, for every $A \in \mathcal{A}$ we extend $\mathcal{E}_0(\cdot,A)$ to $L_{\mathrm{loc}}^1(\mathbb{R}^n)$ by taking
\begin{equation*}
    \mathcal{E}_0(u,A)= \begin{cases}
     \infty \:\: \text{if}\:\: u \in L_{\mathrm{loc}}^1(\mathbb{R}^n)\setminus L_{\mathrm{loc}}^2(\mathbb{R}^n)\:\: \text{and}\:\: A \neq \emptyset, \\
     0 \:\:\:\:\, \text{if}\:\: u \in L_{\mathrm{loc}}^1(\mathbb{R}^n)\setminus L_{\mathrm{loc}}^2(\mathbb{R}^n)\:\: \text{and}\:\: A = \emptyset.
    \end{cases}
\end{equation*}
We do the same for $\mathcal{M}_0^\pm$. 
We claim that $\mathcal{E}_0$ satisfies the assumptions of \cite[Theorem 3]{Bouchitt2002AGM}, i.e.,
\begin{itemize}
    \item[$i)$] \NNN $\frac{1}{C}\mathcal{H}^{n-1}(J_u \cap A)\leq \mathcal{E}_0(u,A)\leq C\mathcal{H}^{n-1}(J_u \cap A)$, for some $C\geq 1$, for every $A \in \mathcal{A}$ and $u\in BV(A;\{-1,1\})$;\EEE
    \item[$ii)$] $\mathcal{E}_0(u,A)=\mathcal{E}_0(v,A)$ whenever $u=v$ $\mathcal{L}^n$-a.e. on $A \in \mathcal{A}$;
    \item[$iii)$] $\mathcal{E}_0(\cdot,A)$ is lower semicontinuous with respect to the strong $L^1(A)$ topology;
    \item[$iv)$] $\mathcal{E}_0(u,\cdot)$ is the restriction to $\mathcal{A}$ of a Radon measure.
\end{itemize}
\NNN Property $i)$ is a consequence of \eqref{eq4.30}, \eqref{extendedinnerregularity}, and the fact that $\mathcal{M}_0^\pm$ are inner regular. \EEE
We notice that property \NNN $ii)$ \EEE follows using Lemma \ref{increasingsetfunctionlemma}, \cite[Remark 15.25]{DalMaso:93} and arguing like in the proof of \cite[Proposition 16.15]{DalMaso:93}). Because of \cite[Proposition 6.8]{DalMaso:93}, \cite[Remark 15.10]{DalMaso:93}, and the fact that the supremum of lower semicontinuous functions is still lower semicontinuous, we have that for every $A \in \mathcal{A}$ the function $u \in L^2(A)\mapsto \mathcal{E}_0(u,A)$ is lower semicontinuous with respect to the strong $L^2(A)$ topology. Combining this with \eqref{eq4.30}, which implies that $\mathcal{E}_0(u,A)<\infty$ if and only if $u \in BV(A;\{-1,1\})$, \NNN $iii)$ \EEE follows. Indeed on $L^{\infty}$ functions defined on bounded domains, the notion of $L^1$ convergence and $L^2$ convergence coincide. 
Thus we are left with proving \NNN $iv)$\EEE. Notice that, by construction, $\mathcal{E}_0$ defines an increasing set function. In addition, from \cite[Remark 15.10]{DalMaso:93}), and the density of $\mathcal{A}_1$ in $\mathcal{A}$, it follows that $\mathcal{E}_0$ is also inner regular. Arguing like in Lemma \ref{increasingsetfunctionlemma} below, and using the additivity of $\mathcal{E}_{\varepsilon}$, it can be shown that the set function $\mathcal{E}_0(u,\cdot)$ is superadditive on $\mathcal{A}_1$, and thus, with a density argument, on the whole $\mathcal{A}$. \\
\EEE We are left to prove \EEE that $\mathcal{E}_0(u,\cdot)$ is a subadditive set function for every $u \in L_{\mathrm{loc}}^2(\mathbb{R}^n)$. \EEE Indeed, by combining this and the previous properties with the De Giorgi-Letta criterium \cite[Theorem 14.23]{DalMaso:93}, we will get that $\mathcal{E}_0(u,\cdot)$ is the restriction to $\mathcal{A}$ of a Borel measure and so $iv)$ holds. \EEE \\
Let $U$, $V \in \mathcal{A}$. Because of \eqref{eq4.30}, we can restrict the proof to the case $u \in BV(U\cup V;\{-1,1\})$. 
Let $O$, $A$, $A^\prime$, $C \in \mathcal{A}_1$ with $O \subset \subset A \subset \subset A^\prime \subset \subset U $ and $C \subset \subset  V  $.
We choose $C$ also in such a way that $\mathcal{H}^{n-1}(J_u \cap \partial C)=0$. 
Let $B\in \mathcal{A}_1$, such that $(A^\prime \setminus \overline{A})\cup C \subset \subset B \subset \subset (U \setminus \overline{O})\cup V$. Since $\mathcal{E}_{\varepsilon}$ $\Gamma$-converges to $\mathcal{E}_0$ on sets in $\mathcal{A}_1$, there exist sequences $(u_{\varepsilon})_{\varepsilon}\subset W^{2,2}(A^\prime)$ and $(v_{\varepsilon})_{\varepsilon}\subset W^{2,2}(B)$ such that $u_{\varepsilon}\to u$ in $L^2(A^\prime)$, $v_{\varepsilon}\to u$ in $L^2(B)$,
\begin{equation}
\label{recoverysequences}
    \lim\limits_{\varepsilon \to 0}\mathcal{E}_{\varepsilon}(u_{\varepsilon},A^\prime)=\mathcal{E}_0(u,A^\prime),\:\: \text{and}\:\:  \lim\limits_{\varepsilon\to 0}\mathcal{E}_{\varepsilon}(v_{\varepsilon},B)=\mathcal{E}_0(u,B).
\end{equation}
\EEE At the end of the proof we will show that \EEE
\begin{equation}
\label{subadditivityimproved}
\mathcal{E}_0(u,B)=\mathcal{E}_0(u,B \setminus \overline{C})+\mathcal{E}_0(u,{C}). 
\end{equation}
\EEE With this, by applying \EEE Proposition \ref{fundamentalestimate},\eqref{eq4.30}, \eqref{recoverysequences} and \eqref{subadditivityimproved}, we get
\begin{equation}
\begin{split}
\mathcal{E}_0(u,A\cup B)\leq & \mathcal{E}_0(u,A^\prime)+\mathcal{E}_0(u,B)=\mathcal{E}_0(u,A^\prime)+\mathcal{E}_0(u,B \setminus \overline{C})+\mathcal{E}_0(u,C)\\
\leq
&\mathcal{E}_0(u,A^\prime)+c_2\mathcal{M}^+_0(u,(U\setminus \overline{O})\cup (V \setminus \overline{C}))+\mathcal{E}_0(u,C).      
\end{split}
\end{equation}
By sending $O \uparrow U$, $C \uparrow V$, using \cite[Lemma 14.20]{DalMaso:93} and the inner regularity of $\mathcal{E}_0(u,\cdot)$, we get $\mathcal{E}_0(u,U\cup V)\leq \mathcal{E}_0(u,U)+\mathcal{E}_0(u,V)$ and so the subadditivity of $\mathcal{E}_0(u,\cdot)$. \EEE \\
\EEE Thus, we are left with showing \eqref{subadditivityimproved}. \EEE We start by observing that 
\begin{equation}
\label{additionalstep}
    \mathcal{E}_0(u,B\setminus \partial C)=\mathcal{E}_0(u,B \setminus \overline{C})+\mathcal{E}_0(u,{C}).    
\end{equation}
In fact, one inequality follows from the superadditivity of $\mathcal{E}_0(u,\cdot)$. \EEE The other inequality is a consequence of the following argument. Let $\EEE W \subset \subset B \setminus \overline{C}$, with $W\in \mathcal{A}_1$. Let $(\tilde{u}_{\varepsilon})_{\varepsilon}$ be a recovery sequence for $u$ on $W$ and $(\tilde{v}_{\varepsilon})_{\varepsilon}$ be a recovery sequence for $u$ on $C$. We define $\tilde{w}_{\varepsilon}\in W^{2,2}(W\cup C)$ as
\begin{equation*}
    \tilde{w}_{\varepsilon}(x)=\begin{cases}
        \tilde{u}_{\varepsilon}(x)\:\: \text{if}\:\: x \in W \\
        \tilde{v}_{\varepsilon}(x)\:\:\: \text{if}\:\: x \in C.
    \end{cases}
\end{equation*}
Then, from $\mathcal{E}_{\varepsilon}(\tilde{w}_{\varepsilon}, W \cup C)=\mathcal{E}_{\varepsilon}(\tilde{u}_{\varepsilon},W)+\mathcal{E}_{\varepsilon}(\tilde{v}_{\varepsilon},C)$ we obtain, by sending $\varepsilon \to 0$, using that $W \cup C \in \mathcal{A}_1$ and that $\mathcal{E}_0(u,\cdot)$ is an increasing set function, 
\begin{equation*}
\mathcal{E}_0(u,W\cup C)\leq \mathcal{E}_0(u,W)+\mathcal{E}_0(u,C)\leq \mathcal{E}_0(u,B \setminus \overline{C})+\mathcal{E}_0(u,C).    
\end{equation*}
By sending $W \uparrow B \setminus \overline{C}$, using the inner regularity of $\mathcal{E}_0(u,\cdot)$ and \cite[Lemma 14.20]{DalMaso:93}, \eqref{additionalstep} follows. \EEE
Thus, in order to obtain \eqref{subadditivityimproved}, it is sufficient to show that
\begin{equation}
\label{eq*}
\mathcal{E}_0(u,B \setminus \partial C)=\mathcal{E}_0(u,B).
\end{equation}
\EEE To this end, let $ D \subset \subset E \subset \subset E^\prime \subset \subset B \setminus \partial C$, with $D$, $E$, $E^\prime\in \mathcal{A}_1$. We notice that $E^\prime \setminus \overline{E} \subset B \setminus \overline{D}$. 
Applying Proposition \ref{fundamentalestimate}, on the sets $E$, $E^\prime$ and $B \setminus \overline{D}$ (in place of $A$, $A^\prime$ and $B$ respectively), \eqref{growthcondition}, and using standard arguments involving recovery sequences for $u$ combined with fact that $\mathcal{H}^{n-1}(J_u \cap \partial C)=0$, \EEE one can prove
\begin{equation*}
\begin{split}
    \mathcal{E}_0(u,B\setminus \partial C)\leq \mathcal{E}_0(u,B)=\mathcal{E}_0(u,E \cup (B\setminus \overline{D}))&\leq \mathcal{E}_0(u,E^\prime)+\mathcal{E}_0(u,B\setminus \overline{D})  \\&\leq \mathcal{E}_0(u,B \setminus \partial C)+ c_2\mathcal{M}^+_0(u,B \setminus \overline{D})\\
    &=\mathcal{E}_0(u,B \setminus \partial C)+\EEE c_2\mathcal{M}^+_0(u,(B \setminus \partial C)\setminus \overline{D})\EEE,    
\end{split}
\end{equation*}
Then by arbitrariness of \EEE $D \subset \subset B \setminus \partial C$ \EEE we obtain \eqref{eq*}.

\end{proof}

\subsection{Properties and identification of the $\Gamma$-limit}
This subsection is dedicated to the study of the limit density $f_0$ and to the proof of Proposition \ref{propositioncontroldensity}.  

For $\rho$, $\delta$, $\varepsilon>0$ with $\rho>\delta>2\varepsilon$ we set
\begin{equation*}
    m_{\mathcal{E}_{\varepsilon}}^\delta({u}_{x,\varepsilon}^\nu, Q_{\rho}^\nu(x))\defas \inf\{\mathcal{E}_\varepsilon(u,Q_{\rho}^\nu(x))\colon u \in \mathcal{S}^\delta({u}_{x,\varepsilon}^\nu, Q_{\rho}^\nu(x))\},
\end{equation*}
where 
\begin{equation*}
\mathcal{S}^\delta({u}_{x,\varepsilon}^\nu, Q_{\rho}^\nu(x)) \defas \{u \in W^{2,2}(Q_{\rho}^\nu(x))\colon u={u}_{x,\varepsilon}^\nu\:\: \text{on}\:\: Q_{\rho}^\nu(x)\setminus \overline{Q}_{\rho-\delta}^{\nu}(x)\}.   
\end{equation*}
Then, for every $x \in  \mathbb{R}^n$ and $\mathbb{S}^{n-1}$, we define
\begin{equation}
\label{f_rho^prime}
f_{\rho}^\prime(x,\nu)\defas\inf_{\delta>0}\liminf\limits_{\varepsilon \to 0}m^\delta_{\mathcal{E}_{\varepsilon}}({u}_{x,\varepsilon}^\nu,Q_{\rho}^\nu(x))=\lim\limits_{\delta \to 0}\liminf\limits_{\varepsilon \to 0}m^\delta_{\mathcal{E}_{\varepsilon}}({u}_{x,\varepsilon}^\nu,Q_{\rho}^\nu(x)),
\end{equation}
and
\begin{equation}
\label{f_rho^second}
f_{\rho}^{\prime \prime}(x,\nu)\defas\inf_{\delta>0}\limsup\limits_{\varepsilon \to 0}m^\delta_{\mathcal{E}_{\varepsilon}}({u}_{x,\varepsilon}^\nu,Q_{\rho}^\nu(x))=\lim\limits_{\delta \to 0}\limsup\limits_{\varepsilon \to 0}m^\delta_{\mathcal{E}_{\varepsilon}}({u}_{x,\varepsilon}^\nu,Q_{\rho}^\nu(x)).
\end{equation}
\begin{proposition}
\label{prop6.1}
Let $f_{\rho}^\prime$, $f_{\rho}^{\prime \prime}$ as in \eqref{f_rho^prime} and \eqref{f_rho^second} respectively. Then the following holds:
\begin{enumerate}
    \item[$(i)$] The restrictions of $f_{\rho}^\prime, f_{\rho}^{\prime \prime}$ to the sets $\mathbb{R}^n \times \hat{\mathbb{S}}_+^{n-1}$ and $\mathbb{R}^n \times \hat{\mathbb{S}}_-^{n-1}$ are upper semicontinuous. 
    \item[($ii$)] For every $x \in \mathbb{R}^n$ and $\nu \in \mathbb{S}^{n-1}$, the functions $\rho \to f_{\rho}^\prime(x,\nu)-c_2 C_{\eta}\rho^{n-1}$ and $\rho \to f_{\rho}^{\prime \prime}(x,\nu)-c_2 C_{\eta}\rho^{n-1}$ are not increasing on $(0,\infty)$, \NNN where $C_{\eta}$ is the positive constant defined in \eqref{Ceta}\EEE. 
    \item[$(iii)$] For every $x \in \mathbb{R}^n$ and $\nu \in \mathbb{S}^{n-1}$ we have
    \begin{equation}
    \label{eq6.3}
    f^\prime(x,\nu)= \limsup\limits_{\rho \to 0}\frac{1}{\rho^{n-1}}f_{\rho}^\prime(x,\nu)
    \end{equation}
    and 
    \begin{equation}
    \label{eq6.4}
    f^{\prime \prime}(x,\nu)= \limsup\limits_{\rho \to 0}\frac{1}{\rho^{n-1}}f_{\rho}^{\prime \prime}(x,\nu).
    \end{equation}
\end{enumerate}
\end{proposition}

\begin{proof} 
The proof is an adaptation of the one in \cite[Lemma 6.1]{marziani2022gammaconvergence} up to some modifications. These modifications are mostly due to the fact that the statement of our fundamental estimate requires the sets to have $C^1$-boundaries. As a matter of fact, the only information regarding the explicit form of the energies that the author makes use of are the fundamental estimate and an analogue version of \eqref{controlsurface}.  We give a sketch of the proof only for $(i)$ and for $f_{\rho}^\prime$, since $(ii)$--$(iii)$ can be obtained following exactly the one of \cite[Lemma 6.1 $(ii)$--$(iii)$] {marziani2022gammaconvergence} up to replace \cite[(2.10)]{marziani2022gammaconvergence} with \eqref{controlsurface}. 
Let $\rho>0$ and let $(x_j,\nu_j)\subset \mathbb{R}^n \times \hat{\mathbb{S}}_{\pm}^{n-1}$ be a sequence converging to $(x,\nu)\in \mathbb{R}^{n}\times \hat{\mathbb{S}}^{n-1}$. We claim that
\begin{equation}
\label{eq6.5}
\limsup\limits_{j \to \infty}f^{\prime}_{\rho}(x_j,\nu_j)\leq f^{\prime}_{\rho}(x,\nu).    
\end{equation}
Let $\gamma>0$. Up to passing to a subsequence (not relabeled) we can assume that the $\liminf\limits$ in \eqref{f_rho^prime} is a limit.
Following the same argument of \cite[Lemma 6.1 $(i)$]{marziani2022gammaconvergence}, which mostly only exploits the continuity of the restrictions of the map $\nu \to R_{\nu}$ on $\hat{\mathbb{S}}_{\pm}^{n-1}$, we get that for $\delta>0$ sufficiently small there exist a $j_0=j_0(\delta)$ such that
\begin{equation*}
Q_{\rho-5\delta}^{\nu_j}(x_j)\subset Q_{\rho-4\delta}^\nu(x)\subset Q^{\nu}_{\rho-2\delta}(x)\subset Q_{\rho-\delta}^{\nu_j}(x_j)
\end{equation*}
for every $j \geq j_0$, and a sequence $(u_\varepsilon)_\varepsilon \subset \mathcal{S}^{3 \delta}({u}^\nu_{x,\varepsilon},Q_{\rho}^\nu(x))$ satisfying
\begin{equation}
\label{eq6.6}
\mathcal{E}_\varepsilon(u_\varepsilon,Q_{\rho}^\nu(x))\leq m^{3\delta}_{\mathcal{E}_{\varepsilon}}({u}^{\nu}_{x,\varepsilon},Q_{\rho}^\nu(x))+\gamma\leq c_2C_{\eta}\rho^{n-1}+\gamma
\end{equation}
and
\begin{equation}
\label{eq6.7}
u_\varepsilon={u}_{x,\varepsilon}^\nu={u}_{x_j,\varepsilon}^{\nu_j}\:\: \text{on}\:\:  (Q_{\rho-2\delta}^\nu(x) \setminus \overline{Q}_{\rho-3\delta}^\nu(x))\setminus \overline{R}_{\varepsilon,j},
\end{equation}
where 
\begin{equation*}
R_{\varepsilon,j}\defas R_{\nu}\Big((Q_{\rho-2\delta}^\prime \setminus \overline{Q}_{\rho-3\delta}^\prime)\times (-h_{\varepsilon,j},h_{\varepsilon,j})\Big)+x,
\end{equation*}
with $h_{\varepsilon,j}\defas \varepsilon +\vert x-x_j\vert +\frac{\sqrt{n}}{2}\rho \vert \nu -\nu_j\vert$. Let $U_j^\prime$, $U_j \in \mathcal{A}_1$ such that $Q_{\rho-5\delta}^{\nu_j}(x_j)\subset \subset U_j \subset \subset Q^{\nu}_{\rho-2\delta}(x) \subset \subset U^\prime_j \subset \subset Q_{\rho-\delta}^{\nu_j}(x_j)$.
Fix a $j \geq j_0$. Now, we can apply Proposition \ref{fundamentalestimate} with $A$, $A^\prime \in \mathcal{A}_1$ such that $Q_{\rho-3\delta}^\nu(x)\subset \subset U_j \subset \subset A  \subset \subset A^{\prime}\subset \subset  Q_{\rho-2\delta}^\nu(x)\subset \subset U_j^\prime$, $B=U_j^\prime\setminus \overline{U}_j$, $\NNN v_{\varepsilon}={u}_{x_j,\varepsilon}^{\nu_j}\EEE$, obtaining, \NNN up to taking $\varepsilon$ smaller than some $\varepsilon_0$  depending on $A$, $A^\prime$ and $B$, functions $w_\varepsilon \in W^{2,2}(U_j^\prime)$ \EEE such that
\begin{align}
\label{eq6.8}
&\mathcal{E}_\varepsilon(w_\varepsilon,U_j^\prime)
   \leq  \mathcal{E}_\varepsilon(u_\varepsilon,A^{\prime})+\mathcal{E}_\varepsilon({u}_{x_j,\varepsilon}^{\nu_j},B)+ \omega_\varepsilon(u_\varepsilon,{u}_{x_j,\varepsilon}^{\nu_j},A,A^\prime) \nonumber\\
& =\mathcal{E}_{\varepsilon}(u_{\varepsilon},Q_{\rho}^{\nu}(x))-\mathcal{E}_{\varepsilon}(u^{\nu}_{x,\varepsilon},Q^{\nu}_{\rho}(x)\setminus \overline{A}^{\prime})+\mathcal{E}_\varepsilon({u}_{x_j,\varepsilon}^{\nu_j},B)+\omega_\varepsilon(u_\varepsilon,{u}_{x_j,\varepsilon}^{\nu_j},A,A^\prime) \nonumber \\
&\leq \mathcal{E}_{\varepsilon}(u_{\varepsilon},Q_{\rho}^{\nu}(x))+c_{2}C_{\eta}[(\rho-2\delta)^{n-1}-(\rho-3\delta)^{n-1})+\rho^{n-1}-(\rho-5\delta)^{n-1})]+\omega_\varepsilon(u_\varepsilon,{u}_{x_j,\varepsilon}^{\nu_j},A,A^\prime), 
\end{align}
and $w_\varepsilon={u}_{x_j,\varepsilon}^{\nu_j}$ on $ B \setminus \overline{A}^{\prime}=U_j^{\prime}\setminus \overline{A}^{\prime}$. In \eqref{eq6.8} we also used that $-\mathcal{E}_{\varepsilon}(z,A)\leq -c_1 \mathcal{M}^{-}_{\varepsilon}(z,A)\leq c_2 \mathcal{M}_{\varepsilon}^{+}(z,A)$ for every $A \in \mathcal{A}$ and $z \in W^{2,2}(A)$. We extend then $w_{\varepsilon}$ on the whole $Q_{\rho}^{\nu_j}(x_j)$ by setting it equal to $u^{\nu_j}_{x,\varepsilon_j}$ in $Q_{\rho}^{\nu_j}(x_j) \setminus \overline{U}_j^\prime$. As a consequence, $w_\varepsilon$ is admissible for $m^{\delta}_{\mathcal{E}_{\varepsilon}}({u}_{x_j,\varepsilon}^{\nu_j},Q_{\rho}^{\nu_j}(x_j))$ and we have, by virtue of $\NNN(f_3)\EEE$ and \eqref{controlsurface},
    \begin{equation}
    \label{eq5.9}
    \mathcal{E}_{\varepsilon}(w_{\varepsilon},Q_{\rho}^{\nu_j}(x_j))\leq \mathcal{E}_{\varepsilon}(w_{\varepsilon},U_j^\prime)+c_2C_{\eta}(\rho^{n-1}-(\rho-5\delta)^{n-1}).    
    \end{equation}
Since $R_{\varepsilon,j}\subset Q_{\rho-2\delta}^\nu \setminus \overline{Q}^\nu_{\rho-3\delta}$, we have $u_{\varepsilon}=u^{\nu}_{x,\varepsilon}$ and $\vert u_{\varepsilon}-u^\nu_{x_j,\varepsilon}\vert \leq 2$ on $R_{\varepsilon,j}$. \NNN In addition, because of \eqref{controlsurface}, we have
    \begin{align*}
     &\sup_j\sup_\varepsilon (\mathcal{M}^-_\varepsilon(u_\varepsilon,A^\prime \setminus \overline{A})+\mathcal{M}^-_\varepsilon({u}^\nu_{x_j,\varepsilon},A^\prime \setminus \overline{A})) \\ 
     &\leq \EEE \sup_j \EEE \sup_\varepsilon (\mathcal{M}^+_\varepsilon(u^{\nu}_{x_j,\varepsilon}, Q^{\nu}_{\rho}(x) )+\mathcal{M}^+_\varepsilon({u}^\nu_{x_j,\varepsilon}, Q^{\nu}_{\rho}(x) )<2C_{\eta}\rho^{n-1}<C,   
    \end{align*}
\NNN for some $C$ large enough. \EEE 
Combining this fact with \eqref{eq6.7} and \eqref{omegaestimate} in Proposition \ref{fundamentalestimate}, with $\Tilde{S}=(A^{\prime}\setminus \overline{A})\setminus \overline{R}_{\varepsilon,j}$, we get it holds
\begin{equation}
\label{eq6.9}
\omega_\varepsilon(u_\varepsilon,{u}_{x_j,\varepsilon}^{\nu_j},A,A^\prime)\leq C\mathcal{L}^n(R_{\varepsilon,j})\NNN+C\varepsilon.\EEE  
\end{equation}
Notice now that 
\begin{equation}
\label{eq6.10}
    \mathcal{L}^n(R_{\varepsilon,j})\leq C\delta \rho^{n-2}h_{\varepsilon,j}.
\end{equation}
Combining \eqref{eq6.6} and \eqref{eq6.8}--\eqref{eq6.10} we get 
\begin{align*}
m^{\delta}_{\mathcal{E}_{\varepsilon}}({u}_{x_j,\varepsilon}^{\nu_j},Q_{\rho}^{\nu_j}(x_j))&\leq m^{3\delta}_{\mathcal{E}_{\varepsilon}}({u}^{\nu}_{x,\varepsilon},Q_{\rho}^\nu(x))+\gamma +c_{2}C_{\eta}((\rho-2\delta)^{n-1}-(\rho-3\delta)^{n-1})\\
&+c_{2}C_{\eta}(\rho^{n-1}-(\rho-5\delta)^{n-1}))+ C\delta \rho^{n-2}h_{\varepsilon,j}+\NNN C\varepsilon\EEE,    
\end{align*}
and by sending, in this order, $\varepsilon \to 0$, then $\delta \to 0$ , $j \to \infty$, and $\gamma \to 0$, we obtain
\begin{equation*}
    \limsup\limits_{j \to \infty} f_{\rho}^\prime(x_j,\nu_j) \leq f_{\rho}^\prime(x,\nu)
\end{equation*}
from which \eqref{eq6.5} follows.   
\end{proof}
\begin{lemma}[Equivalence of infimum problems with different boundary conditions]
\label{lemma5.2}
Define for every $x \in \mathbb{R}^n$, $\nu \in \mathbb{S}^{n-1}$ and $\varepsilon$, $\rho>0$, $m_{\mathcal{M}^\pm_{\varepsilon}}(u^{\nu}_{x,\varepsilon},Q^{\nu}_{\rho}(x))$ as in \eqref{infimumproblem} up to replace $\mathcal{E}_{\varepsilon}$ with $\mathcal{M}^{\pm}_{\varepsilon}$, $\Tilde{u}$ with $u^{\nu}_{x,\varepsilon}$, and $A$ with $Q^\nu_{\rho}(x)$.
Then, it holds
\begin{equation}
    \lim\limits_{\varepsilon \to 0}m_{\mathcal{M}^\pm_{\varepsilon}}(u^\nu_{x,\varepsilon},Q^\nu_{\rho}(x))=\sigma^\pm \rho^{n-1}.
\end{equation}
\end{lemma}
\begin{proof} 
Let $\Tilde{m}_{\mathcal{M}^{\pm}_{\varepsilon}}(u^{\nu}_{x,\varepsilon},Q^{\nu}_{\rho}(x))$ be as in \eqref{FMinfimumproblem}. Clearly by virtue of Proposition \ref{claimsection6} and the properties of the infimum it holds
\begin{equation}
\sigma^\pm\rho^{n-1}=\lim\limits_{\varepsilon\to 0}\Tilde{m}_{\mathcal{M}^\pm_{\varepsilon}}(u^\nu_{x,\varepsilon},Q^\nu_{\rho}(x))\leq \liminf\limits_{\varepsilon \to 0}m_{\mathcal{M}^\pm_{\varepsilon}}(u^\nu_{x,\varepsilon},Q^\nu_{\rho}(x)).  
\end{equation}
Thus, we have just to show
\begin{equation}
\label{eq5.15}
\limsup\limits_{\varepsilon \to 0}m_{\mathcal{M}^\pm_{\varepsilon}}(u^\nu_{x,\varepsilon},Q^\nu_{\rho}(x))\leq \sigma^\pm \rho^{n-1}. 
\end{equation}
Up to extracting a subsequence (not relabeled) we can assume that the $\limsup\limits$ in \eqref{eq5.15} is a limit. Let $\delta \in (0,1)$ and $A$, $A^\prime$, $U \in \mathcal{A}_1$ such that $Q_{(1-\delta)\rho}^\nu(x)\subset \subset A \subset \subset A^\prime \subset \subset U \subset \subset Q_{\rho}^\nu(x)$ and $B \defas U \setminus \overline{A}$. Thanks to \cite[Theorem 1.3]{Chermisi}, there exists a recovery sequence $(u^\pm_{\varepsilon})_{\varepsilon}\subset W^{2,2}(A^{\prime})$ such that $u^\pm_{\varepsilon} \to \overline{u}_x^\nu$ in $L^2(A^{\prime})$ and
\begin{equation}
\label{eq5.16}
    \limsup\limits_{\varepsilon \to 0}\mathcal{M}_{\varepsilon}^\pm(u^\pm_{\varepsilon},A^{\prime})\leq \mathcal{M}_0^\pm(\overline{u}_x^\nu,A^{\prime})\leq \sigma^\pm \rho^{n-1}.
\end{equation}
We now apply Proposition \ref{fundamentalestimate} with $c_1 \mathcal{M}^-_{\varepsilon}$ and $c_2\mathcal{M}_{\varepsilon}^+$ instead of $\mathcal{E}_{\varepsilon}$, $u \defas u_{\varepsilon}^\pm$ and $v\defas u_{x,\varepsilon}^\nu$, obtaining, up to \NNN taking $\varepsilon$ small, \EEE functions $w^{\pm}_{\varepsilon}\in W^{2,2}(U)$ such that $w^\pm_{\varepsilon}=u_{x,\varepsilon}^\nu$ on $U \setminus \overline{A}^{\prime}$ and 
\begin{equation}
\label{eq5.17}
    \mathcal{M}_{\varepsilon}^\pm(w^{\pm}_{\varepsilon},U)\leq \mathcal{M}^\pm_{\varepsilon}(u^\pm_{\varepsilon},A^{\prime})+\frac{c_2}{c_1}\mathcal{M}_{\varepsilon}^+(u^\nu_{x,\varepsilon},Q_{\rho}^\nu(x)\setminus \overline{Q}_{(1-\delta)\rho}^\nu(x))+\frac{\omega_{\varepsilon}(u_{\varepsilon},u_{x,\varepsilon}^\nu,A,A^\prime)}{c_1},
\end{equation}
where we also used \eqref{growthcondition} and that $\mathcal{M}^{+}_{\varepsilon}$, having a nonnegative density, defines an increasing set function.
Now we can extend $w^{\pm}_{\varepsilon}$ on all $Q^\nu_\rho(x)$ by setting it equal to $u^\nu_{x,\varepsilon}$ on $Q^\nu_\rho(x) \setminus \overline{U}$. 
\NNN Notice that, since $Q_{(1-\delta)\rho}^\nu(x)\subset U \subset Q_{\rho}^\nu(x)$ and because of \eqref{controlsurface}, we have
\begin{equation*}
\mathcal{M}_{\varepsilon}^\pm(w^{\pm}_{\varepsilon},Q^\nu_{\rho}(x))\leq \mathcal{M}_{\varepsilon}^\pm(w^{\pm}_{\varepsilon},U)+C_{\eta}\rho^{n-1}(1-(1-\delta)^{n-1}).    
\end{equation*}
\EEE
By applying \NNN again \EEE \eqref{controlsurface}, we get \NNN from \eqref{eq5.17} \EEE
\begin{equation}
\label{eq5.18}
    \mathcal{M}_{\varepsilon}^\pm(w^{\pm}_{\varepsilon},Q^{\nu}_\rho(x))\leq \mathcal{M}_{\varepsilon}^\pm(u^\pm_{\varepsilon},A^{\prime})+\NNN \Big(\frac{c_2}{c_1}+1\Big)C_{\eta}\EEE\rho^{n-1}(1-(1-\delta)^{n-1})+\NNN \frac{\omega_{\varepsilon}(u_{\varepsilon},u_{x,\varepsilon}^\nu,A,A^\prime)}{c_1}\EEE.
\end{equation}
Finally, notice that $w^{\pm}_{\varepsilon}$ is admissible for $m_{\mathcal{M}^\pm_{\varepsilon}}(u^\nu_{x,\varepsilon},Q^\nu_{\rho}(x))$ and that $u_{\varepsilon}-u_{x,\varepsilon}^\nu \to 0$ in $L^2(A^\prime \setminus \overline{A})$. 
Passing to the limit as $\varepsilon \to 0$ in \eqref{eq5.18}, using \eqref{limomega} in Proposition \ref{fundamentalestimate}, and \eqref{eq5.16}, we get
\begin{equation}
\limsup\limits_{\varepsilon \to 0}m_{\mathcal{M}^\pm_{\varepsilon}}(u^\nu_{x,\varepsilon},Q^\nu_{\rho}(x))\leq  \sigma^\pm \rho^{n-1}+\NNN \Big(\frac{c_2}{c_1}+1\Big)C_{\eta}\EEE\rho^{n-1}(1-(1-\delta)^{n-1}),  
\end{equation}
and by letting $\delta \to 0$ we obtain \eqref{eq5.15}.
\end{proof}
Finally, Proposition  \ref{propositioncontroldensity} follows from the previous results of this section. 
\begin{proof}[Proof of Proposition \ref{propositioncontroldensity}]
\NNN The proof of \eqref{firstequationprop3.4} relies on standard arguments, for convenience of the reader we include it in the Appendix (Proposition \ref{prop5.1}). The proof of the Borel measurability can be easily adapted from the one of \cite[Proposition 6.2, \emph{Step 1}]{marziani2022gammaconvergence} up to take into account Proposition \ref{prop6.1}. To conclude we just need to prove \eqref{secondequationprop3.4}. Clearly, because of \eqref{growthcondition}, it holds
\begin{equation}
\label{eq6.12}
c_1   m_{\mathcal{M}^-_{\varepsilon}}({u}^\nu_{x,\varepsilon}, Q^\nu_{\rho}(x)) \leq m_{\mathcal{E}_{\varepsilon}}({u}^\nu_{x,\varepsilon}, Q^\nu_{\rho}(x)) \leq c_2 {m}_{\mathcal{M}_{\varepsilon}^+}({u}^\nu_{x,\varepsilon}, Q^\nu_{\rho}(x)).
\end{equation}
Hence, \eqref{secondequationprop3.4} follows by Lemma \ref{lemma5.2} and the definitions in \eqref{f´} and \eqref{f´´}.
\end{proof}

\section{Stochastic homogenization}
\label{sec: stochastichomogenization}
This section is devoted to the proof of Theorem \ref{stochhomformula} and Theorem \ref{stochastichomogenization}.
We state two lemmas that will be used in the proof of one of the main results of this section, namely Proposition \ref{propx0}. 
\begin{lemma}
\label{lemmaA.1bis}
Let $f \in \mathcal{F}$, let $\nu \in \mathbb{S}^{n-1}$, $x, \Tilde{x}\in \mathbb{R}^n$, and $\Tilde{r}>r>4$ be such that
\begin{equation*}
    (i)\:\: Q_{r+2}^\nu(x)\subset \subset Q^\nu_{\Tilde{r}}(\Tilde{x})\quad (ii)\:\: \mathrm{dist}(\Tilde{x},\Pi^\nu(x))<\frac{r}{4},
\end{equation*}
where $\Pi^\nu(x)$ is the hyperplane orthogonal to $\nu$ and passing through $x$.
Then there exists a constant $L>0$ (independent of $\nu,x,\Tilde{x},r,\Tilde{r}$) such that
\begin{equation}
\label{lemmaA.1}
m_{\mathcal{E}}(u_{\Tilde{x}}^\nu,Q_{\Tilde{r}}^\nu(\Tilde{x}))\leq m_{\mathcal{E}}(u_x^\nu,Q_{{r}}^\nu({x}))+L(\vert x-\Tilde{x}\vert+\vert r-\Tilde{r}\vert +1)\Tilde{r}^{n-2}.
\end{equation}
\end{lemma}
\begin{lemma}
\label{lemmaA.2}
Let $f \in \mathcal{F}$, $\alpha \in \big(0,\frac{1}{2}\big)$, and $\nu,\Tilde{\nu}\in \mathbb{S}^{n-1}$ be such that
\begin{equation}
\label{eqA.2}
    \max_{1\leq i \leq n-1}\vert R_{\nu}e_i-R_{\Tilde{\nu}}e_i\vert + \vert \nu - \Tilde{\nu}\vert \leq \frac{\alpha}{\sqrt{n}}.
\end{equation}
Then there exists a constant $c_{\alpha}>0$ (independent of $\nu$, $\Tilde{\nu}$), with $c_{\alpha}\to 0$ as $\alpha \to 0$, such that for every $x \in \mathbb{R}^n$ and every $r>2$ we have
\begin{equation}
\label{LemmaA.2}
m_{\mathcal{E}}(u_{rx}^{\Tilde{\nu}}, Q_{(1+\alpha)r}^{\Tilde{\nu}}(rx))-c_{\alpha}r^{n-1}    \leq m_{\mathcal{E}}(u_{rx}^{{\nu}}, Q_{r}^{{\nu}}(rx)).
\end{equation}
\end{lemma}
The proof of Lemma \ref{lemmaA.1bis} and Lemma \ref{lemmaA.2} can be easily adapted from the ones in \cite[Lemma A.1, Lemma A.2]{marziani2022gammaconvergence} up to minor modifications. Indeed we can observe that because of Theorem \ref{positivesubadditiveprocess} and $\NNN(f_2)\EEE$, $m_{\mathcal{E}}(u^{\nu}_x,Q^{\nu}_{r}(x))\geq 0$ if $r\geq 1$, and thus one can find competitors with energy equal to $m_{\mathcal{E}}(u^{\nu}_x,Q^{\nu}_{r}(x))$ up to an arbitrary small error. Then, the only modifications of the proofs are related to the fact that our energies are not necessarily subadditive set functions. However one can check that an analogue version of \cite{marziani2022gammaconvergence}[(A.3) and (A.11)], as well as the remaining inequalities in the proofs can be obtained using the additivity of the functional $\mathcal{E}$, \eqref{growthcondition} (which implies $\vert \mathcal{E}(\cdot,A)\vert \leq c_2\mathcal{M}^{+}(\cdot,A)$), \eqref{controlsurface} and the fact that $\mathcal{M}^{+}$ defines an increasing set function. \EEE \\
For any positive integer $k$ and for any $a,b \in \mathbb{R}^k$, with $a_i<b_i$ for all $i=1,...,k$, we define the cuboids  
\begin{equation*}
    [a,b)=\prod^{k}_{i=1}[a_i,b_i)=\{x \in \mathbb{R}^k: a_i \leq x_i < b_i \:\: \forall \, i=1,...,k\},
\end{equation*}
and 
\begin{equation*}
\mathcal{R}_k = \{[a,b): a_i<b_i \:\: \forall \, i=1,...,k\}.
\end{equation*}
Then, given $R=\prod^{n-1}_{j=1}[a_j,b_j) \in \mathcal{R}_{n-1}$, for every $\nu \in \mathbb{S}^{n-1}\cap \mathbb{Q}^{n}$, we define the rotated $n$-dimensional cuboid $T_{\nu}R$ as 
\begin{equation}
\label{Tnu}
\begin{split}
T_{\nu}R:= M_{\nu}R_{\nu}(R \times [-c,c)),  \quad \text{where } c:=\max \Big\{\frac{1}{2},\max_{1\leq j\leq n-1}\frac{(b_j-a_j)}{2}\Big\}.
\end{split}
\end{equation}
Notice that the \NNN set $R_{\nu}(R\times [-c,c))$ \EEE adds \NNN to \EEE $\NNN R\EEE$ a thin layer of at least length $1$ along the $\nu$ direction. This is related to the length of the transition layer described by the function $u^{\nu}_{x}$.

\begin{definition}{(Subadditive process)}
\label{subadditivedef}
A \emph{\NNN(bounded) subadditive process} with respect to a group $(\tau_{z})_{z \in \mathbb{Z}^k}$ (resp.\  $(\tau_{z})_{z \in \mathbb{R}^k}$) of $\mathbb{P}$-preserving transformations on $(\Omega, \mathcal{I}, \mathbb{P})$ is a function $\mu \colon \Omega \times \mathcal{R}_k  \to \mathbb{R}$ satisfying:\\
\begin{enumerate}
    \item[$(\mathrm{i})$] (measurability) for any $R \in \mathcal{R}_k$ the function $\omega \to \mu(\omega,R)$ is $\mathcal{I}$ measurable, 
    \item[$(\mathrm{ii})$] (\NNN stationarity) \EEE for any $\omega \in \Omega$, $R \in \mathcal{R}_k$, $z\in \mathbb{Z}^k$ (resp.\  $z\in \mathbb{R}^k$) it holds $\mu(\tau_z \omega,R)=\mu(\omega,z+R)$,
    \item[$(\mathrm{iii})$] (subadditivity) for any $R \in \mathcal{R}_k$ and for any finite family $(R_i)_{i \in I} \subset \mathcal{R}_k$ of pairwise disjoint sets  such that $\bigcup_{i\in I}R_i=R$, it holds
\begin{equation*}
    \mu(\omega,R)\leq \sum_{i \in I}\mu(\omega,R_i)\:\: \text{for every}\:\: \omega \in \Omega,
\end{equation*}
\item[$(\mathrm{iv})$] (uniform boundedness) there exists a constant $c>0$ such that $0\leq\mu(\omega,R) \leq c\,\mathcal{L}^k (R)$ for every $\omega \in \Omega$ and for every $R \in \mathcal{R}_k$.
\end{enumerate}
\end{definition}
For every $\nu \in \mathbb{S}^{n-1}$ we consider the function $\mu_{\nu}\colon \Omega \times \mathcal{R}_{n-1}\to \mathbb{R}$ defined by 
\begin{equation}
\label{subadditiveprocess}
\mu_{\nu}(\omega,R)\defas \begin{cases}
\frac{m_{\mathcal{E}[\omega]}(u_0^\nu,\mathrm{int} (T_{\nu}R))}{M^{n-1}_{\nu}} \:\: \text{if}\:\: \mathcal{L}^1({I}_i)\geq 1\:\: \text{for every}\:\:i \in \{1,...,n-1\}
\\
c_2 C_{\eta}\mathcal{L}^{n-1}(R)\:\: \quad \:\: \text{otherwise},
\end{cases}
\end{equation}
for every $\omega \in \Omega$ and $R=\prod^{n-1}_{i=1}I_i \in \mathcal{R}_{n-1}$.
We notice that because of Theorem \ref{positivesubadditiveprocess}, $(\NNN f_2\EEE)$, and the definition of $T_{\nu}$, it holds $\mu_{\nu}(\omega,R)\geq 0$. Notice that the definition of $\mu_{\nu}$ is slightly different compared to the ones that can be found in the literature (see for example \cite[(5.3)]{cagnetti2017stochastic}, \cite[(7.5)]{marziani2022gammaconvergence}). The main reason, is that $m_{\mathcal{E}[\omega]}(u_0^\nu,\mathrm{int} (T_{\nu}R))$ is in principle not positive if $R$ has one side of length strictly less than 1. Our definition coincides with the one in \cite{cagnetti2017stochastic} and \cite{marziani2022gammaconvergence} on large cuboids. In addition, using $\NNN(f_3)\EEE$ and \eqref{controlsurface} we get
\begin{equation}
\label{subadditivecontrolsurface}
\frac{m_{\mathcal{E}[\omega]}(u_0^\nu,\mathrm{int} (T_{\nu}R))}{M^{n-1}_{\nu}}\leq \frac{{\mathcal{E}}[\omega](u_0^\nu,\mathrm{int} (T_{\nu}R))}{M^{n-1}_{\nu}}\leq c_2 C_{\eta}\mathcal{L}^{n-1}(R).  
\end{equation}

\begin{proposition}
\label{propsubadd}
Let $f$ be a random density which is stationary with respect to a group of $\mathbb{P}$-preserving transformation $(\tau_z)_{z \in \mathbb{Z}^n}$ (resp. $(\tau_z)_{z \in \mathbb{R}^n}$) on $(\Omega,\mathcal{I},\mathbb{P})$. Let $\nu \in \mathbb{S}^{n-1}\cap \mathbb{Q}^n$ and let $\mu_{\nu}$ be as in \eqref{subadditiveprocess}. Then, there exists a group of $\mathbb{P}$-preserving transformations $(\tau_{z^\prime})_{z^\prime \in \mathbb{Z}^{n-1}}$ (resp. $(\tau_{z^\prime})_{z^\prime \in \mathbb{R}^{n-1}}$) such that $\mu_{\nu}$ is a \NNN (bounded) \EEE subadditive process on $(\Omega,\mathcal{I},\mathbb{P})$ with respect to  $(\tau_{z^\prime})_{z^\prime \in \mathbb{Z}^{n-1}}$ (resp. $(\tau_{z^\prime})_{z^\prime \in \mathbb{R}^{n-1}}$). Moreover it holds
\begin{equation}
\label{uniformboundedness}
0\leq \mu_{\nu}(\omega,R)\leq c_2 C_{\eta}\mathcal{L}^{n-1}(R),
\end{equation}
for $\mathbb{P}$-almost every $\omega \in \Omega$ and $R\in \mathcal{R}_{n-1}$.
\end{proposition}
\begin{proof}
\emph{Step 1: Measurability}. \NNN The proof of the measurability of $\omega \mapsto \mu_{\nu}(\omega,R)$ can be easily adapted from the one in \cite[Lemma B.1]{RufSch23} (see also \cite[Lemma A.8]{CaDaScZe21}) up to some minor modifications. 
\\ In particular, for every $l \in \mathbb{N}$ let $f^l$ be defined as
\begin{equation*}
    f^l(\omega,x,u,\xi,\zeta)\defas \begin{cases}
        f(\omega,x,u,\xi,\zeta)\:\: \text{if}\:\: f(\omega,x,u,\xi,\zeta)\geq -l, \\
        -l \quad \quad \quad \quad \quad \: \text{otherwise},
    \end{cases}
\end{equation*}
and let $\mathcal{E}^l$ be the corresponding energy.
\EEE Because of Lemma \ref{appendix1} in the Appendix, it is sufficient to prove that for every \NNN $l \in \mathbb{N}$ \EEE
\begin{equation}
\label{measurabilityclaim}
\begin{split}
\omega \to m_{\mathcal{E}^l[\omega]}(\Tilde{u},A)\:\: \text{is $\mathcal{I}$-measurable, for every}\:\: A \in \mathcal{A}\:\: \text{and}\:\: \Tilde{u}\in W^{2,2}(A).    
\end{split}
\end{equation}
Notice also that
\begin{equation*}
    \inf_{u \in \mathcal{S}(\Tilde{u},A)}\int_{A}\Big(f^l(\omega,x,u,\nabla u,\nabla^2 u)+l \Big) \, \mathrm{d}x= m_{\mathcal{E}^l[\omega]}(\Tilde{u},A)+l\mathcal{L}^n(A)
\end{equation*}
for every $\omega \in \Omega$. Hence, it is not restrictive to assume $f \geq 0$ \NNN in this step. \EEE 
Let $A \in \mathcal{A}$ and let $(A_j)_j \subset \mathcal{A}$ be such that $A_j \uparrow A$. 
Consider the set
\begin{equation*}
    \mathcal{S}_j(\Tilde{u},A)\defas \{u\in W^{2,2}(A)\colon u=\Tilde{u}\:\: \mathcal{L}^n\text{-a.e on}\:\: A \setminus \overline{A}_j\}.
\end{equation*}
Notice that $\mathcal{S}_j(\Tilde{u},A)$ is a closed subset of $W_0^{2,2}(A)+\Tilde{u}$ and thus it defines a complete, separable metric space if equipped with the metric induced by the $\Vert \cdot \Vert_{W^{2,2}(A)}$ norm. 
Since $0\leq \inf_{u \in \mathcal{S}_j(\Tilde{u},A)}\mathcal{E}[\omega](u,A)\leq\mathcal{E}[\omega](\Tilde{u},A)<\infty$ for every $j \in \mathbb{N}$, it holds
\begin{equation}
    \lim\limits_{j \to \infty}\inf_{u \in \mathcal{S}_j(\Tilde{u},A)}\mathcal{E}[\omega](u,A)=m_{\mathcal{E}[\omega]}(\Tilde{u},A).
\end{equation}
Hence, in order to prove \eqref{measurabilityclaim}, it is sufficient to show that
\begin{equation}
\label{measurabilityclaim2}
\omega \to \inf_{u \in \mathcal{S}_j(\Tilde{u},A)}\mathcal{E}[\omega](u,A)\:\: \text{is}\:\: \mathcal{I}\text{-measurable}.
\end{equation}
\NNN We focus for the rest of the proof of this step in showing \eqref{measurabilityclaim2}. The strategy consists in showing first that
\begin{equation}
\label{newmeasurabilityclaim}
(\omega,u)\in \Omega \times \mathcal{S}_j(\tilde{u},A)\mapsto \int_{A}f(\omega,x,u,\nabla u,\nabla^2u)\,\mathrm{d}x \:\: \text{is}\:\:  \mathcal{I}\otimes\mathcal{B}(\mathcal{S}_j(\tilde{u},A))\text{-measurable}, 
\end{equation}
if $(u,\xi,\zeta)\to f(\omega,x,u,\xi,\zeta)$ is continuous for every $\omega \in \Omega$ and $x \in \mathbb{R}^n$. Then, we apply the monotone class theorem (\cite[Chapter 1, Theorem 21]{meyer1966probability}) in order to show \eqref{newmeasurabilityclaim} for general nonnegative random densities. Finally, we deduce \eqref{measurabilityclaim2} by applying the measurable projection theorem (\cite[Theorem III.13 and 33(a)]{meyer1966probability}). \EEE
\NNN Hence, let us assume that the function $(u,\xi,\zeta)\mapsto f(\omega,x,u,\xi,\zeta)$ is continuous for every $\omega \in \Omega$ and $x\in \mathbb{R}^n$. Up to a truncation argument, it is not restrictive to assume that $f$ is also bounded from above.
\EEE The \EEE joint measurability of $(\omega,u)\mapsto \int_{A}f(\omega,u,\nabla u,\nabla^2 u)\, \mathrm{d}x$ \EEE then \EEE follows from the Fubini's theorem, the continuity of $(u,\xi,\zeta)\mapsto f(\omega,x,u,\xi,\zeta)$ and the separability of $\mathcal{S}_j(\tilde{u},A)$. \EEE Indeed, by Fubini's theorem the map $\omega \mapsto \int_{A}f(\omega,x,u,\nabla u,\nabla^2 u)\, \mathrm{d}x$ is $\mathcal{I}$-measurable, while, with the dominated convergence theorem, one can show that the functional $u \in \mathcal{S}_j(\tilde{u},A) \mapsto \int_{A}f(\omega,x,u,\nabla u,\nabla^2 u)\, \mathrm{d}x$ is continuous with respect to the strong $W^{2,2}(A)$ topology, for every $\omega \in \Omega$. Then the joint measurability is just a consequence of these facts and the separability of $\mathcal{S}_j(\tilde{u},A)$ (see for example \cite[Section 11C, Exercise 11.3]{Kechris}). \EEE  We proceed now by showing that \eqref{newmeasurabilityclaim} holds for general nonnegative random densities. Because of $\NNN(f_3)\EEE$, there exists a function $\hat{f}\colon \Omega \times \mathbb{R}^n\times \mathbb{R}\times \mathbb{R}^n \times \mathbb{R}^{n \times n}\to [0,c_2]$ such that $f(\omega,x,u,\xi,\zeta)=\hat{f}(\omega,x,u,\xi,\zeta)(W(u)+\vert \xi \vert^2+\vert \zeta \vert^2)$. Define the sets
\begin{equation*}
\begin{split}
    \mathcal{C}\defas \{ &\hat{f}\colon \Omega \times \mathbb{R}^n\times \mathbb{R}\times \mathbb{R}^n \times \mathbb{R}^{n \times n}\to [0,c_2] \colon \hat{f}(\omega,x,u,\xi,\zeta)=\alpha(\omega)\beta(x)\gamma(u,\xi,\zeta),\:\: \text{with}\\ &\alpha\:\: \text{bounded and}\:\:\mathcal{I}\text{-measurable},\:\: \beta \:\: \text{bounded and}\:\: \mathcal{B}^n\text{-measurable},\:\: \gamma\:\: \text{continuous and bounded}\}  
\end{split}
\end{equation*}
and 
\begin{equation*}
\begin{split}
  \mathcal{R}\defas \{\hat{f}\colon \Omega \times \mathbb{R}^n\times \mathbb{R}\times \mathbb{R}^n \times \mathbb{R}^{n \times n}\to [0,c_2]\colon \eqref{newmeasurabilityclaim}\:\: \text{holds with}\\ f(\omega,x,u,\xi,\zeta)=\hat{f}(\omega,x,u,\xi,\zeta)(W(u)+\vert \xi \vert^2+\vert \zeta \vert^2)\}.    
\end{split}  
\end{equation*}
Notice that $\mathcal{C}$ is stable under multiplication and that the $\sigma$-algebra generated by $\mathcal{C}$ coincides with $\mathcal{I}\otimes \mathcal{B}^n\otimes \mathcal{B}\otimes \mathcal{B}^n \otimes \mathcal{B}^{n \times n}$. Instead, $\mathcal{R}$
satisfies the definition of monotone class (see for example \cite[Definition 4.12]{aliprantis06}).
Because of what we showed before, we have $\mathcal{C}\subset \mathcal{R}$ and so the monotone class theorem implies that $\mathcal{R}$ must contain all the bounded functions $\hat{f}$, proving in this way \eqref{newmeasurabilityclaim} for general nonnegative random densities. Now, we notice that \eqref{newmeasurabilityclaim} implies
\begin{equation}
\label{measurableset}
\Big\{ (\omega,u)\in \Omega \times \mathcal{S}_{j}(\tilde{u},A)\colon \int_{A}f(\omega,x,u,\nabla u, \nabla^2 u)\, \mathrm{d}x<t\Big\}\:\: \text{is}\:\: \mathcal{I}\otimes \mathcal{B}(\mathcal{S}_j(\tilde{u},A))\text{-measurable},   
\end{equation}
for every $t \in \mathbb{R}$. Finally \eqref{measurabilityclaim2}, follows by the measurable projection theorem. Indeed, let us denote with $\pi \colon \Omega \times \mathcal{S}_j(\tilde{u},A)\to \Omega$ the canonical projection of $\Omega \times \mathcal{S}_j(\tilde{u},A)$ onto $\Omega$. Since $(\Omega,\mathcal{I},\mathbb{P})$ is a complete probability space and since $\mathcal{S}_j(\Tilde{u},A)$ is a complete, separable metric space, the measurable projection theorem and \eqref{measurableset} lead to the $\mathcal{I}$-measurability of
$$\pi\Big(\Big\{ (\omega,u)\in \Omega \times \mathcal{S}_{j}(\tilde{u},A)\colon \int_{A}f(\omega,x,u,\nabla u, \nabla^2 u)\, \mathrm{d}x<t\Big\}\Big)=\Big\{\omega \in \Omega \colon \inf_{u \in \mathcal{S}_j(\Tilde{u},A)}\mathcal{E}[\omega](u,A)<t\Big\},$$ for every $t \in \mathbb{R}$, and so to \eqref{measurabilityclaim2}.
\EEE
\\
\emph{Step 2: \NNN Stationarity.} The proof of the \NNN stationarity \EEE is rather standard (see \cite{cagnetti2017stochastic} for example) but we report it for sake of completeness. We show this property only for the more difficult case, which is the one with a discrete group $(\tau_{z^\prime})_{z^\prime \in \mathbb{Z}^{n-1}}$. In view of \eqref{subadditiveprocess}, one has to check only the case in which $R$ has all sides of at least length 1. For every $z^\prime \in \mathbb{Z}^{n-1}$ and $R \in \mathcal{R}_{n-1}$ it holds
\begin{equation*}
    T_{\nu}(R+z^\prime)=T_{\nu}R+z_{\nu}^\prime,
\end{equation*}
where $z_{\nu}^\prime\defas M_{\nu}R_{\nu}(z^\prime,0) \in \mathbb{Z}^n \cap \Pi^\nu$. Thus, we have
\begin{equation*}
\mu_{\nu}(\omega,R+z^\prime)=\frac{1}{M^{n-1}_{\nu}}m_{\mathcal{E}[\omega]}(u_0^\nu,\mathrm{int} (T_{\nu}R)+z_{\nu}^\prime).
\end{equation*}
Now, let $u \in \mathcal{S}(u_0^\nu, T_{\nu}R+z_{\nu}^\prime)$ and define $\Tilde{u}\colon \mathrm{int}(T_{\nu}R)\to \mathbb{R}$ as $\Tilde{u}(y)\defas u(y+z_\nu^\prime)$. Since $z_{\nu}^\prime \in \Pi^\nu$, we have that $u^{\nu}_0(y)=u_0^{\nu}(y+z_{\nu}^\prime)$ for every $y \in \mathbb{R}^n$. Hence, $\Tilde{u}\in \mathcal{S}(u_0^\nu, T_{\nu}R)$. With a change of variable and using the stationarity of $f$ we get
\begin{align}
\label{eq6.11}
\mathcal{E}[\omega](u, \mathrm{int}(T_{\nu}R+z_{\nu}^\prime))&=\int_{T_{\nu}R+z_{\nu}^\prime}f(\omega,x,u,\nabla u,\nabla^2 u)\, \mathrm{d}x=\int_{T_{\nu}R}f(\omega,x+z_{\nu}^\prime,\Tilde{u},\nabla \Tilde{u},\nabla^2 \Tilde{u})\, \mathrm{d}x \nonumber  \\
&=\int_{T_{\nu}R}f(\tau_{z_{\nu}^\prime}\omega,x,\Tilde{u},\nabla \Tilde{u},\nabla^2 \Tilde{u})\, \mathrm{d}x = \mathcal{E}[\tau_{z_{\nu}^\prime}\omega](\Tilde{u},\mathrm{int}(T_{\nu}R)).  \end{align}
As a consequence, it is natural to consider the group of $\mathbb{P}$-preserving transformations $(\tau^\nu_{z^\prime})_{z^\prime \in \mathbb{Z}^{n-1}}$ defined by $(\tau_{z_{\nu}^\prime})_{z^\prime \in \mathbb{Z}^{n-1}}$. Equation \eqref{eq6.11}, combined with the arbitrariness of $u$, gives
\begin{equation*}
    \mu_{\nu}(\omega,R+z^\prime)=\mu_{\nu}(\tau_{z^\prime}^\nu \omega,R),
\end{equation*}
which proves the \NNN stationarity \EEE of $\mu_{\nu}$ with respect to the previously defined group $(\tau^\nu_{z^\prime})_{z^\prime \in \mathbb{Z}^{n-1}}$. \\
\emph{Step 3: Subadditivity.} Let $R \in \mathcal{R}_{n-1}$ and let $R_1,...,R_N \in \mathcal{R}_{n-1}$ be pairwise disjoint and such that $R=\bigcup^{N-1}_{i=1}R_i$. 
For every $i \in \{1,...,N\}$, we take $(u^j_i)_j \subset W^{2,2}(\mathrm{int}(T_{\nu}R_i))$ as minimizing sequences for $m_{\mathcal{E}[\omega]}(u_0^\nu,T_{\nu}R_i)$. 
Then, we define then $u^j \colon T_{\nu}R \to \mathbb{R}$ as
\begin{equation*}
    u^j(x)=\begin{cases}
        u^j_i(x)\:\: \text{if}\:\: x \in T_{\nu}R_i \\
        u^\nu_0(x)\:\: \text{otherwise}.
    \end{cases}
\end{equation*}
Notice that $u^j \in W^{2,2}(\mathrm{int}(T_{\nu}R))$ since, for every $i \in \{1,...,N\}$, $u^j_i$ is equal to $u^\nu_0$ in a neighborhood of $\partial (T_{\nu}R_i)$. In particular, by construction $u^j$ is admissible for $m_{\mathcal{E}[\omega]}(u_0^\nu,T_{\nu}R)$ and it holds
\begin{equation}
\label{eq6.13}
m_{\mathcal{E}[\omega]}(u^{\nu}_0, \mathrm{int}(T_{\nu}R))\leq \mathcal{E}[\omega](u^j,\mathrm{int}(T_{\nu}R))= \sum^N_{i=1}\mathcal{E}[\omega](u^j_i,\mathrm{int}(T_{\nu}R_i))+\mathcal{E}[\omega]\Big(u_0^\nu,\mathrm{int}(T_{\nu}R \setminus (\bigcup^N_{i=1}T_{\nu}R_i))\Big).    
\end{equation}
Now notice that $\mathcal{E}[\omega](u_0^\nu,\mathrm{int}(T_{\nu}R \setminus (\bigcup^N_{i=1}T_{\nu}R_i))=0$, in fact since $c\geq\frac{1}{2}$ and $M_{\nu}\geq 2$ we have $\{y \in T_{\nu}R: \vert y \cdot \nu \vert \leq \frac{1}{2}\}\subset \bigcup^N_{i=1}T_{\nu}R_i$. As a consequence of this, it must hold $u_0^\nu=\overline{u}_0^\nu$ in $T_{\nu}R \setminus (\bigcup^N_{i=1}T_{\nu}R_i)$. Thus, by sending $j \to \infty$ in \eqref{eq6.13} it follows
\begin{equation}
\label{subadditiveinfimumpr}
m_{\mathcal{E}[\omega]}(u_0^\nu,T_{\nu}R)\leq \sum^N_{i=1}m_{\mathcal{E}[\omega]}(u_0^\nu,T_{\nu}R_i).  
\end{equation}
From the last equation and by the definition in \eqref{subadditiveprocess}, we obtain that $\mu_{\nu}$ satisfies the subadditivity inequality when all the $R_i$ have sides of length larger or equal than 1. \\
If at least one $R_i$ has at least one side of length strictly less than 1, we have two possible cases. If $R$ has still all sides of length larger or equal than one, then the subadditivity inequality follows by \eqref{subadditiveinfimumpr} and \eqref{subadditivecontrolsurface}. If $R$ has at least one side of length lesser or equal than 1, then all the $R_i$ must have at least one side of length lesser or equal than 1 as well. Hence, in this case, the subadditive inequality follows from the additivity of the $n-1$ dimensional Lebesgue measure and observing that $\mu_{\nu}(\omega,R_i)=c_2 C_{\eta}\mathcal{L}^{n-1}(R_i)$ for every $i$.
\end{proof}
\begin{proposition}[Homogenized surface integrand for $x=0$]
\label{propx0}
Let $f$ be a random density which is stationary with respect to a group of $\mathbb{P}$-preserving transformations $(\tau_z)_{z \in \mathbb{Z}^n}$ (resp. $(\tau_z)_{z \in \mathbb{R}^n}$ ) on $(\Omega,\mathcal{I},\mathbb{P})$. For every $\omega \in \Omega$ let $m_{\mathcal{E}[\omega]}$ as in \eqref{randominfimumproblem}. 
Then, there exists an $\mathcal{I}\otimes \mathcal{B}(\mathbb{S}^{n-1})$-measurable function $f_{\mathrm{hom}}\colon \Omega \times \mathbb{S}^{n-1}\to [0,\infty)$ such that, given 
\begin{equation}
\label{omegatilde}
    \Tilde{\Omega}\defas \Big\{ \omega \in \Omega \colon \lim\limits_{r \to \infty}\frac{m_{\mathcal{E}[\omega]}(u^\nu_0,Q_r^\nu(0))}{r^{n-1}}=f_{\mathrm{hom}}(\omega,\nu)\:\: \text{for every}\:\: \nu \in \mathbb{S}^{n-1}\Big\},
\end{equation}
we have $\Tilde{\Omega}\in \mathcal{I}$ and $\mathbb{P}(\Tilde{\Omega})=1$. 
Moreover, $\Tilde{\Omega}$ and $f_{\mathrm{hom}}$ are $(\tau_z)_{z \in \mathbb{Z}^n}$ (resp. $(\tau_z)_{z \in \mathbb{R}^n}$) invariant, i.e., $\tau_z(\Tilde{\Omega})=\Tilde{\Omega}$ for every $z \in \mathbb{Z}^n$ (resp. for every $z \in \mathbb{R}^n$) and 
\begin{equation}
\label{tauinvariance}
f_{\mathrm{hom}}(\tau_z \omega, \nu)=f_{\mathrm{hom}}(\omega,\nu),    
\end{equation}
for every $z \in \mathbb{Z}^n$ (resp. $z \in \mathbb{R}^n$), $\omega \in \Tilde{\Omega}$ and $\nu \in \mathbb{S}^{n-1}$. \EEE Eventually, if $(\tau_z)_{z \in \mathbb{Z}^n}$ (resp. $(\tau_z)_{z \in \mathbb{R}^n}$) is ergodic, then $f_{\mathrm{hom}}$ is independent of $\omega$ and given by
\begin{equation}
\label{detdensity}
f_{\mathrm{hom}}(\nu)=\lim\limits_{r \to \infty}\frac{1}{r^{n-1}}\int_{\Omega}m_{\mathcal{E}[\omega]}(u_0^\nu,Q_r^\nu(0))\, \mathrm{d}\mathbb{P}(\omega).
\end{equation}
\end{proposition}

\begin{proof}
The proof follows the same steps as in \cite[Proposition 7.7]{marziani2022gammaconvergence} that we report for convenience of the reader. We divide the proof in three steps. \\
\emph{Step 1: Existence of the limit for $\nu \in \mathbb{S}^{n-1}\cap \mathbb{Q}^n$}. In this step we show that there exists an event $\overline{\Omega} \in \mathcal{I}$, with $\mathbb{P}(\overline{\Omega})=1$, such that for every $\nu \in \mathbb{S}^{n-1}\cap \mathbb{Q}^n$ there exists an $\mathcal{I}$-measurable function $f_{\nu}\colon \Omega \to [0,\infty)$ satisfying
\begin{equation*}
\lim\limits_{r \to \infty}\frac{m_{\mathcal{E}[\omega]}(u^\nu_0,Q_r^\nu(0))}{r^{n-1}}=f_{\nu}(\omega)    
\end{equation*}
for every $\omega \in \overline{\Omega}$. Fix $\nu \in \mathbb{S}^{n-1} \cap \mathbb{Q}^{n-1}$. Because of Proposition \ref{propsubadd}, we know that $\mu_{\nu}$ defines a subadditive process. Thus, we can apply the Subadditive Ergodic Theorem in \cite[Theorem 3.11]{cagnetti2017stochastic} (see also \cite[Theorem 2.7]{Krengel1981}) with $\mu_{\nu}$ instead of $\mu$ and $2{Q^\prime}$ replacing $Q$, finding in this way an event $\Omega_{\nu} \in \mathcal{I}$ with probability 1, and an $\mathcal{I}$-measurable function $f_{\nu}\colon \Omega \to [0,\infty)$ such that
\begin{equation}
\label{eq6.18bis}
f_{\nu}(\omega)=\lim\limits_{r \to \infty}\frac{m_{\mathcal{E}[\omega]}(u^\nu_0,2M_{\nu}Q_r^{\nu}(0))}{(2M_{\nu}r)^{n-1}}\:\: \text{for every} \:\:\omega \in \Omega_{\nu}, 
\end{equation}
where we used that $T_{\nu}(2Q_1^\prime)=2M_{\nu}Q_{1}^\nu(0)$. Now we claim that
\begin{equation}
\label{eq6.18}
f_{\nu}(\omega)=\lim\limits_{r \to \infty}\frac{m_{\mathcal{E}[\omega]}(u^\nu_0,Q_r^\nu(0))}{r^{n-1}},
\end{equation}
for every $\omega \in \Omega_{\nu}$. Let $(r_j)_j$ be a sequence such that $r_j \to \infty$ as $j \to \infty$, and define
\begin{equation*}
    r^-_j\defas 2M_{\nu}\Big(\Big\lfloor\frac{r_j}{2M_{\nu}}\Big\rfloor-1\Big)\quad \text{and} \quad r_j^+\defas 2M_{\nu}\Big(\Big\lfloor\frac{r_j}{2M_{\nu}}\Big\rfloor+2\Big).
\end{equation*}
Now, up to take $j$ sufficiently large we can assume $r_j>4(1+M_{\nu})$ and thus $r^-_j>4$. Furthermore, it holds 
\begin{equation*}
Q_{r^-_j+2}^\nu(0)\subset \subset Q_{r_j}^\nu(0)\subset \subset Q_{r_j+2}^\nu(0)\subset \subset Q_{r^+_j}^\nu(0).  
\end{equation*}
We apply Lemma \ref{lemmaA.1bis} twice: the first time with $x=\Tilde{x}=0$, $r=r^-_j$ and $\Tilde{r}=r_j$, the second time with $x=\Tilde{x}=0$, $r=r_j$, and $\Tilde{r}=r^+_j$, and get the two following estimates
\begin{equation}
\label{eq6.19}
    \frac{m_{\mathcal{E}[\omega]}(u^\nu_0,Q_{r_j}^\nu(0))}{r^{n-1}_j}\leq \frac{m_{\mathcal{E}[\omega]}\Big(u^\nu_0,Q_{r^-_j}^\nu(0)\Big)}{(r^-_j)^{n-1}}+\frac{L(r_j-r^-_j+1)}{r_j},
\end{equation}
\begin{equation}
\label{eq6.20}
    \frac{m_{\mathcal{E}[\omega]}(u^\nu_0,Q_{r_j}^\nu(0))}{r_j^{n-1}}\geq \frac{m_{\mathcal{E}[\omega]}\Big(u^\nu_0,Q_{r^+_j}^\nu(0)\Big)}{(r^+_j)^{n-1}}-\frac{L(r^+_j-r_j+1)}{r_j}.
\end{equation}
Using that $r^+_j-r_j\leq 4M_{\nu}$ and $r_j-r^-_j\leq 4M_{\nu}$, \eqref{eq6.18bis}, passing to the $\limsup$ in \eqref{eq6.19} and to the $\liminf$ in \eqref{eq6.20} as $j \to \infty$, yields 
\begin{equation*}
    \limsup\limits_{j \to \infty}\frac{m_{\mathcal{E}[\omega]}(u^\nu_0,Q_{r_j}^\nu(0))}{r_j^{n-1}}\leq f_{\nu}(\omega)
\end{equation*}
and 
\begin{equation*}
 \liminf\limits_{j \to \infty}\frac{m_{\mathcal{E}[\omega]}(u^\nu_0,Q_{r_j}^\nu(0))}{r_j^{n-1}}\geq f_{\nu}(\omega),    
\end{equation*}
for every $\omega \in \Omega_{\nu}$. Thus \eqref{eq6.18} holds since the sequence $(r_j)_j$ was chosen arbitrarily. Finally, we conclude by setting
\begin{equation}
\label{eq6.22}
 \overline{\Omega}\defas \bigcap_{\nu \in \mathbb{S}^{n-1}\cap \mathbb{Q}^n}\Omega_{\nu}  
\end{equation}
and by noticing $\overline{\Omega}\in \mathcal{I}$ and $\mathbb{P}(\overline{\Omega})=1$, being $\overline{\Omega}$ the countable intesection of probability 1 elements of $\mathcal{I}$.  \\
\emph{Step 2: Existence of the limit for $\nu \in \mathbb{S}^{n-1}\setminus \mathbb{Q}^n$}. In this step we prove that there exists a $\mathcal{I}\otimes \mathcal{B}(\mathbb{S}^{n-1})$-measurable function $f_{\mathrm{hom}}\colon \Omega \times \mathbb{S}^{n-1}\to [0,\infty)$ such that
\begin{equation}
\label{eq6.15}   
\lim\limits_{r \to \infty}\frac{m_{\mathcal{E}[\omega]}(u^\nu_0,Q_r^\nu(0))}{r^{n-1}}=f_{\mathrm{hom}}(\omega,\nu)
\end{equation}
holds for every $\omega \in \overline{\Omega}$ and $\nu \in \mathbb{S}^{n-1}$. \EEE Let $\underline{f}, \overline{f}\colon \Omega \times \mathbb{S}^{n-1}\to [0,\infty)$ be defined by
\begin{equation*}
    \underline{f}(\omega,\nu)\defas \liminf\limits_{r \to \infty}\frac{m_{\mathcal{E}[\omega]}(u_0^\nu,Q^\nu_r(0))}{r^{n-1}}\quad \text{and}\quad \overline{f}(\omega,\nu)\defas \limsup\limits_{r \to \infty}\frac{m_{\mathcal{E}[\omega]}(u_0^\nu,Q^\nu_r(0))}{r^{n-1}}.
\end{equation*}
Notice that using Proposition \ref{propsubadd}, and arguing like in \emph{Step 1}, it can be deduced that, for every $\nu \in \mathbb{S}^{n-1}$, $\underline{f}(\cdot,\nu)$ and $\overline{f}(\cdot,\nu)$, can be written as $\liminf$ and $\limsup$ of sequence of $\mathcal{I}$-measurable functions (see also for example \cite[proof of Theorem 7]{Morfe_2020}). Now we observe that $\hat{\mathbb{S}}_{\pm}^{n-1}\cap \mathbb{Q}^n$ is dense in $\hat{\mathbb{S}}_{\pm}^{n-1}$, and that $\mathbb{S}^{n-1}=\hat{\mathbb{S}}_{-}^{n-1}\cup \hat{\mathbb{S}}_{+}^{n-1}$. Furthermore, because of \emph{Step 1}, it holds $\underline{f}(\omega,\nu)=\overline{f}(\omega,\nu)=f_{\nu}(\omega)$ for every $\omega \in \overline{\Omega}$ and $\nu \in \mathbb{S}^{n-1}\cap \mathbb{Q}^n$. Hence, it is enough to show that, for every $\omega \in \overline{\Omega}$, the restrictions of $\underline{f}(\omega,\cdot)$ and $\overline{f}(\omega,\cdot)$ to $\hat{\mathbb{S}}_{\pm}^{n-1}$ are continuous. Indeed, this implies $\underline{f}(\omega,\nu)=\overline{f}(\omega,\nu)$ for every $\nu \in \mathbb{S}^{n-1}$ and $\omega \in \overline{\Omega}$. In addition, we point out that one can deduce from 
\begin{equation*}
\text{for every}\:\: \nu \in \mathbb{S}^{n-1}, \:\: \omega \to \overline{f}(\omega,\nu)\:\: \text{is} \:\:\mathcal{I}-\text{measurable},    
\end{equation*}
together with
\begin{equation*}
\text{for every}\:\: \omega \in \Omega,\:\:\nu \to \overline{f}(\omega,\nu)  \:\:  \text{is continuous on}\:\: \mathbb{S}_{\pm}^{n-1},
\end{equation*}
that $\overline{f}$ is $\mathcal{I}\otimes \mathcal{B}(\mathbb{S}^{n-1})$-measurable. The same applies to $\underline{f}$. \\ From now on, we focus on showing that $\overline{f}(\omega,\cdot)$ is continuous on $\hat{\mathbb{S}}_+^{n-1}$. The proof considering $\underline{f}(\omega,\cdot)$ and/or $\hat{\mathbb{S}}_-^{n-1}$ is analogous.
Consider the function $f_{\mathrm{hom}}\colon \Omega \times \mathbb{S}^{n-1}\to [0,\infty)$ defined as
\begin{equation*}
    f_{\mathrm{hom}}(\omega,\nu)\defas \begin{cases}
        \overline{f}(\omega,\nu)\:\: \text{if}\:\: \omega \in \overline{\Omega}\\
        c_2\sigma^+ \quad\:\: \text{otherwise}.
    \end{cases}
\end{equation*}
  Let $(\nu_j)_j \subset \hat{\mathbb{S}}_+^{n-1}$ and $\nu \in \hat{\mathbb{S}}_+^{n-1}$ be such that $\nu_j \to \nu$ as $j \to \infty$. For every $\alpha \in (0,\frac{1}{2})$ we can find a $j_{\alpha}\in \mathbb{N}$ such that \eqref{eqA.2} (with $\nu_j$ instead of $\nu$) \EEE holds for every $j \geq j_{\alpha}$. Hence, we can apply Lemma \ref{lemmaA.2} with $x=0$ and $\Tilde{\nu}=\nu_j$ getting
\begin{equation*}
    m_{\mathcal{E}[\omega]}(u^{\nu_j}_0,Q_{(1+\alpha)r}^{\nu_j}(0))-c_{\alpha}r^{n-1}\leq m_{\mathcal{E}[\omega]}(u^{\nu}_0,Q_r^\nu(0))\leq m_{\mathcal{E}[\omega]}(u^{\nu_j}_0,Q_{(1-\alpha)r}^{\nu_j}(0))+c_{\alpha}r^{n-1},
\end{equation*}
where $c_{\alpha}\to 0$ as $\alpha \to 0$. By dividing the above inequality by $r^{n-1}$ and passing to the $\limsup$ as $r \to \infty$ we get
\begin{equation}
\label{eq6.24}
(1+\alpha)^{n-1}\overline{f}(\omega,\nu_j)\leq \overline{f}(\omega,\nu)+c_{\alpha},
\end{equation}
\begin{equation}
\label{eq6.25}
(1-\alpha)^{n-1}\overline{f}(\omega,\nu_j)\geq \overline{f}(\omega,\nu)-c_{\alpha}.
\end{equation}
Now, we can pass to the $\limsup$ as $j \to \infty$ in \eqref{eq6.24} and to the $\liminf$ as $j \to \infty$ in \eqref{eq6.25}, and then let $\alpha \to 0$, obtaining 
\begin{equation*}
    \limsup\limits_{j \to \infty}\overline{f}(\omega,\nu_j)\leq \overline{f}(\omega,\nu)\leq \liminf\limits_{j \to \infty}\overline{f}(\omega,\nu_j),
\end{equation*}
for every $\omega \in \overline{\Omega}$. Let $\Tilde{\Omega}$ be as in \eqref{omegatilde}. From the fact that $\overline{\Omega}\subset \Tilde{\Omega}$, $\mathbb{P}(\overline{\Omega})=1$ and $(\Omega,\mathcal{I},\mathbb{P})$ is a complete probability space, it follows $\Tilde{\Omega}\in \mathcal{I}$ and $\mathbb{P}(\Tilde{\Omega})=1$. \EEE \\
\emph{Step 3: Translation invariance}.  In this step we show that $\Tilde{\Omega}$ and $f_{\mathrm{hom}}$ are $(\tau_z)_{z \in \mathbb{Z}^n}$ (resp. $(\tau_z)_{z \in \mathbb{R}^n}$) invariant. Notice that, by virtue of the group properties of $(\tau_z)_{z \in \mathbb{Z}^n}$ (resp. $(\tau_z)_{z \in \mathbb{R}^n}$), in order to show that $\Omega$ is $(\tau_z)_{z \in \mathbb{Z}^n}$-invariant (resp. $(\tau_z)_{z \in \mathbb{R}^n}$-invariant), it is sufficient to show
\begin{equation*}
    \tau_z(\Tilde{\Omega})\subset \Tilde{\Omega}\:\: \text{for every}\:\:z \in \mathbb{Z}^n\:\: (\text{resp}. \:\: z \in \mathbb{R}^n).
\end{equation*} 
Let $z \in \mathbb{Z}^n$ (resp. $z \in \mathbb{R}^n$), $\omega \in \Tilde{\Omega}$, and $\nu \in \mathbb{S}^{n-1}$ be fixed. Let $r>4$. Because of Theorem \ref{positivesubadditiveprocess} and $\NNN(f_2)\EEE$, we know that $m_{\mathcal{E}[\omega]}(u^\nu_0,Q^\nu_r(0))$ is bounded from below and thus there exists some $u \in \mathcal{S}(u_0^\nu,Q_r^\nu(0))$ satisfying 
\begin{equation}
\label{eq6.26}
\mathcal{E}[\omega](u,Q_r^\nu(0))\leq m_{\mathcal{E}[\omega]}(u^\nu_0,Q^\nu_r(0))+1.
\end{equation}
Setting $\Tilde{u}(z)\defas u(y+z)$ and using the stationarity of $f$, we obtain
\begin{equation*}
    \mathcal{E}[\omega](u,Q_r^\nu(0))=\mathcal{E}[\tau_z\omega](\Tilde{u},Q_r^\nu(-z)).
\end{equation*}
This together with \eqref{eq6.26} and the fact that $\Tilde{u}\in \mathcal{S}(u_{-z}^\nu,Q_r^\nu(-z))$ yields
\begin{equation}
\label{eq6.27}
    m_{\mathcal{E}[\tau_z(\omega)]}(u_{-z}^\nu,Q_r^\nu(-z))\leq m_{\mathcal{E}[\omega]}(u^\nu_0,Q_r^\nu(0))+1.
\end{equation}
We choose $r$, $\Tilde{r}$ such that $\Tilde{r}>r$ and
\begin{equation*}
  Q_{r+2}^\nu(-z)\subset \subset Q_{\Tilde{r}}^\nu(0)\:\: \text{and}\:\: \mathrm{dist}(0,\Pi^\nu(-z))\leq \frac{r}{4}.  
\end{equation*}
Then we apply Lemma \ref{lemmaA.1bis} twice: once with $x=-z$ and $\Tilde{x}=0$ to the minimization problem $m_{\mathcal{E}[\tau_{-z}(\omega)]}$ and once with $x=z$ and $\Tilde{x}=0$ to the minimization problem $m_{\mathcal{E}[\omega]}$. Thus, we get
\begin{equation}
\label{eq6.28}
m_{\mathcal{E}[\tau_z \omega]}(u^\nu_0,Q_{\Tilde{r}}^\nu(0))\leq m_{\mathcal{E}[\tau_z \omega]}(u^\nu_{-z},Q_{{r}}^\nu(-z))+L(\vert z \vert +\vert r-\Tilde{r}\vert +1)(\Tilde{r})^{n-2},
\end{equation}
and
\begin{equation}
\label{eq6.29}
m_{\mathcal{E}[\omega]}(u^\nu_0,Q_{\Tilde{r}}^\nu(0))\leq m_{\mathcal{E}[ \omega]}(u^\nu_{z},Q_{{r}}^\nu(z))+L(\vert z \vert +\vert r-\Tilde{r}\vert +1)(\Tilde{r})^{n-2}.
\end{equation}
Thus combining \eqref{eq6.26}--\eqref{eq6.29} we obtain
\begin{equation}
\label{eq6.30}
\frac{m_{\mathcal{E}[\tau_z \omega]}(u^\nu_0,Q_{\Tilde{r}}^\nu(0))}{\Tilde{r}^{n-1}}\leq \frac{m_{\mathcal{E}[\omega]}(u^\nu_0,Q_{\Tilde{r}}^\nu(0))+1}{r^{n-1}}+\frac{L(\vert z \vert +\vert r-\Tilde{r}\vert +1)}{\Tilde{r}}
\end{equation}
and
\begin{equation}
\label{eq6.31}
\frac{m_{\mathcal{E}[ \omega]}(u^\nu_0,Q_{\Tilde{r}}^\nu(0))}{\Tilde{r}^{n-1}}\leq \frac{m_{\mathcal{E}[\tau_z\omega]}(u^\nu_0,Q_{\Tilde{r}}^\nu(0))+1}{r^{n-1}}+\frac{L(\vert z \vert +\vert r-\Tilde{r}\vert +1)}{\Tilde{r}}.
\end{equation}
Now we take the $\limsup$ as $\Tilde{r}\to \infty$ and the limit as $r \to \infty$ in \eqref{eq6.30} to find
\begin{equation}
\label{eq6.32}
\limsup\limits_{\Tilde{r}\to \infty}\frac{m_{\mathcal{E}[\tau_z \omega]}(u^\nu_0,Q_{\Tilde{r}}^\nu(0))}{\Tilde{r}^{n-1}}\leq f_{\mathrm{hom}}(\omega,\nu).
\end{equation}
Similarly we take in \eqref{eq6.29} the limit as $\Tilde{r}\to \infty$ and the $\liminf$ as $r \to \infty$, obtaining
\begin{equation}
\label{eq6.33}
f_{\mathrm{hom}}(\omega,\nu)\leq \liminf\limits_{r \to \infty} \frac{m_{\mathcal{E}[\tau_z\omega]}(u^\nu_0,Q_{\Tilde{r}}^\nu(0))}{r^{n-1}}.
\end{equation}
Gathering \eqref{eq6.32} and \eqref{eq6.33} we deduce $\tau_{z}(\omega)\in \Tilde{\Omega}$ and that
\begin{equation*}
    f_{\mathrm{hom}}(\tau_z(\omega),\nu)= f_{\mathrm{hom}}(\omega,\nu),
\end{equation*}
for every $z \in \mathbb{Z}^n$ (resp. $z \in \mathbb{R}^n$), $\omega \in \Tilde{\Omega}$ and $\nu \in \mathbb{S}^{n-1}$.  \\
When $(\tau_z)_{z \in \mathbb{Z}^n}$ (resp. $(\tau_z)_{z \in \mathbb{R}^n}$) is ergodic, using the fact that $f_{\mathrm{hom}}$ is $(\tau_z)_{z \in \mathbb{Z}^n}$-invariant (resp. $(\tau_z)_{z \in \mathbb{R}^n}$-invariant), and arguing like at the end of the proof of \cite[Theorem 3.11]{cagnetti2017stochastic} it can be shown that actually $f_{\mathrm{hom}}(\omega,\nu)$ does not depend on $\omega$ for $\mathbb{P}$-a.e $\omega$. Finally, \eqref{detdensity} can be deduced from Proposition \ref{propositioncontroldensity} and the Dominated Convergence Theorem.
\end{proof}
\begin{remark}
\label{dalmasoanalogues}
We observe that arguing like in the second step of the proof of Proposition \ref{propx0} it can be shown that the functions 
\begin{equation*}
\nu \to \liminf\limits_{r \to \infty}\frac{m_{\mathcal{E}[\omega]}(u^\nu_{rx},Q_r^\nu(rx))}{r^{n-1}} \:\: \text{and}\:\:   \nu \to \limsup\limits_{r \to \infty}\frac{m_{\mathcal{E}[\omega]}(u^\nu_{rx},Q_r^\nu(rx))}{r^{n-1}}
\end{equation*}
are continuous on $\hat{\mathbb{S}}_{\pm}^{n-1}$, for every $\omega \in \Omega$ and $x \in \mathbb{R}^n$. 
\end{remark}

\begin{proof}[Proof of Theorem \ref{stochhomformula} and Theorem \ref{stochastichomogenization}]
Arguing like in \cite[Theorem 6.1]{cagnetti2017stochastic} up to minor modifications,  replacing \cite[Theorem 5.1] {cagnetti2017stochastic} with Proposition \ref{propx0} and \cite[Lemma 5.5]{cagnetti2017stochastic} with Remark \ref{dalmasoanalogues}, Theorem \ref{stochhomformula} follows. Indeed the authors' proof relies on general arguments coming from probability theory (e.g.\ the Birkhoff Theorem and Conditional Dominated Convergence Theorem) and a property analogues to \eqref{controlsurface}. In addition, the proof in our case is even simpler since it does not require to estimate the energies added by the presence of jumps.
Finally, Theorem \ref{stochastichomogenization} can be proved by combining Theorem \ref{deterministichomogenization} and Theorem \ref{stochhomformula}.     
\end{proof}

\section*{Acknowledgements} 
The author is thankful to his supervisor Manuel Friedrich for many
helpful discussions and suggestions. This work was supported by the RTG 2339 “Interfaces, Complex Structures, and Singular Limits”
of the German Science Foundation (DFG). The support is gratefully acknowledged. 

\appendix
\section{}
\label{sec: appendix}
In this section we prove the results whose proofs rely on arguments that are standard but included for convenience of the reader, for example the proof of Proposition \ref{prop5.1}. Instead, Lemma \ref{increasingsetfunctionlemma} is used in the proof of Theorem \ref{compintrepr} and it states that the $\Gamma-\liminf\limits$ and the $\Gamma-\limsup\limits$, of the energies \eqref{defenergies}, define increasing set functions on $\mathcal{A}_1$. Indeed, this allows to apply a big part of the localization method arguments to our setting. Lemma \ref{appendix1} states that the infimum problem \eqref{infimumproblem} can be approximated with other infimum problems bounded from below. This fact is used in the proof of Proposition \ref{propsubadd}. \\
Before proving Proposition \ref{prop5.1}, we state and prove two intermediate results (Lemma \ref{inneregularenvelope} and Lemma \ref{subadditivitym0}) that will be useful for applying Theorem \ref{compintrepr} and Proposition \ref{fundamentalestimate} to suitable functions defined on cubes. 
\begin{lemma}
\label{inneregularenvelope}
Let $(\varepsilon_j)_j$ and $\mathcal{E}_0$ as in Theorem \ref{compintrepr}. Then, for every $x \in \mathbb{R}^n$, $\nu \in \mathbb{S}^{n-1}$, and $\rho>0$, if $u_{\varepsilon}\to u$ in $L^2(Q^\nu_{\rho}(x))$ and $u_{\varepsilon}=u^\nu_{x,\varepsilon}$ in a neighborhood of $\partial Q^\nu_{\rho}(x)$, then it holds
\begin{equation}
\label{lipschitzliminf}
    \mathcal{E}_0(u,Q^\nu_{\rho}(x))\leq \liminf\limits_{j \to \infty}\mathcal{E}_{\varepsilon_j}(u_{\varepsilon_j},Q^\nu_{\rho}(x)).  
\end{equation}
\end{lemma}
\begin{proof}
By virtue of Theorem \ref{compintrepr} we have that \eqref{lipschitzliminf} clearly holds if $Q^\nu_{\rho}(x)$ is replaced by $A \in \mathcal{A}_1$. Let $\delta>0$ and $A \in \mathcal{A}_1$ such that $Q^\nu_{\rho}(x)\subset A \subset Q_{(1+\delta)\rho}^\nu(x)$. We extend $u_{\varepsilon}$ and $u$ on the whole $Q_{(1+\delta)\rho}^\nu(x)$ by setting them respectively equal to $u^\nu_{x,\varepsilon}$ and $\overline{u}^\nu_x$ outside $Q^\nu_{\rho}(x)$. Notice that the extension of $u_{\varepsilon}$ still converges to the extension of $u$ in $L^2(Q^{\nu}_{(1+\delta)\rho}(x))$. Using that $\mathcal{E}_0$ defines an increasing set function, \eqref{growthcondition}, \eqref{controlsurface}, the additivity of the set function defined by $\mathcal{E}_{\varepsilon}$, and that $-\mathcal{E}_{\varepsilon}(z,U)\leq -c_1\mathcal{M}^{-}_{\varepsilon}(z,U) \leq c_2 \mathcal{M}_{\varepsilon}^+(z,U)$ for every $z \in L_{\mathrm{loc}}^2(\mathbb{R}^n)$ and $U \in \mathcal{A}$, we have 
\begin{align*}
&\mathcal{E}_0(u,Q^\nu_{\rho}(x))\leq \mathcal{E}_0(u,A) \leq \liminf\limits_{j \to \infty}\mathcal{E}_{\varepsilon_j}(u_{\varepsilon_j},A)=\liminf\limits_{j \to \infty}\mathcal{E}_{\varepsilon_j}(u_{\varepsilon_j},Q_{(1+\delta)\rho}^\nu(x))-\mathcal{E}_{\varepsilon_j}(u^\nu_{x,\varepsilon_j},Q_{(1+\delta)\rho}^\nu(x)\setminus \overline{A}) 
\\ 
&\leq  \liminf\limits_{j \to \infty}\mathcal{E}_{\varepsilon_j}(u_{\varepsilon_j},Q_{\rho}^\nu(x))+\mathcal{E}_{\varepsilon_j}(u^{\nu}_{x.\varepsilon_j},Q_{(1+\delta)\rho}^\nu(x)\setminus \overline{Q}^{\nu}_{\rho}(x))+c_2\mathcal{M}_{\varepsilon_j}^+(u^\nu_{x,\varepsilon_j},Q_{(1+\delta)\rho}^\nu(x)\setminus \overline{Q}_{\rho}^\nu(x))   
\\
&\leq\liminf\limits_{j \to \infty}\mathcal{E}_{\varepsilon_j}(u_{\varepsilon_j},Q_{\rho}^\nu(x))+2c_2C_{\eta}((\rho+\delta)^{n-1}-\rho^{n-1}).
\end{align*}
Hence, \eqref{lipschitzliminf} easily follows by sending $\delta \to 0$ in the last equation. 
\end{proof}
\begin{lemma}
\label{subadditivitym0}
Let $\mathcal{E}_0$ as in Theorem \ref{compintrepr}. Then, for every $x \in \mathbb{R}^n$, $\nu \in \mathbb{S}^{n-1}$ and $\delta$, $\rho>0$ and for every $A \in \mathcal{A}_1$ such that $Q_{\rho}^\nu(x) \subset \subset A \subset \subset Q^\nu_{(1+\delta)\rho}(x)$, it holds
\begin{equation*}
m_{\mathcal{E}_0}(\overline{u}^\nu_x,Q_{(1+\delta)\rho}^\nu(x))\leq  m_{\mathcal{E}_0}(\overline{u}^\nu_x,A)+ c_2C_{\eta}\rho^{n-1}((1+\delta)^{n-1}-1). 
\end{equation*}
\end{lemma}
\begin{proof}
Let $\lambda>0$, $u \in BV(A,\{-1,1\})$ such that $u=\overline{u}^\nu_{x}$ in a neighborhood of $\partial A$ and 
\begin{equation}
\label{eq3.25bis}
    \mathcal{E}_0(u,A)\leq m_{\mathcal{E}_0}(\overline{u}^\nu_x,A)+\lambda.
\end{equation}
Let $A^\prime \in \mathcal{A}_1$ such that $Q^\nu_{(1+\delta)\rho}(x) \subset \subset A^\prime \subset \subset Q^\nu_{(1+\delta+\lambda)\rho}(x)$ and ${A}^\prime \setminus \overline{A}\in \mathcal{A}_1$. Define now $\Tilde{u}\in BV(A^\prime;\{-1,1\})$ as
\begin{equation*}
    \Tilde{u}(y)=\begin{cases}
        u(y)\:\: \text{if}\:\: y \in A \\
        \overline{u}^\nu_x \quad \: \text{if}\:\:y \in A^\prime \setminus A.
    \end{cases}
\end{equation*}
Then, clearly $\Tilde{u}$ is admissible for $m_{\mathcal{E}_0}(\overline{u}^\nu_x,Q_{(1+\delta)\rho}^\nu(x))$ and we have
\begin{equation}
\label{eq3.26}
m_{\mathcal{E}_0}(\overline{u}^\nu_x,Q_{(1+\delta)\rho}^\nu(x))\leq \mathcal{E}_0(\Tilde{u},Q^\nu_{(1+\delta)\rho}(x))\leq \mathcal{E}_0(u,A)+\mathcal{E}_0(\overline{u}^\nu_x,A^\prime\setminus \overline{A}). 
\end{equation}
\NNN Because of \eqref{FonsecaMantegazzalimit} and \eqref{growthcondition} \EEE, we have
\begin{align}
\label{eq3.27} 
\mathcal{E}_0(\overline{u}^\nu_x,A^\prime \setminus \overline{A})&\leq  \mathcal{E}_0(\overline{u}^\nu_x,Q^\nu_{(1+\delta+\lambda)\rho}(x) \setminus \overline{Q}^\nu_{\rho}(x))\leq c_2 \mathcal{M}_0^+(\overline{u}^\nu_x,Q^\nu_{(1+\delta+\lambda)\rho}(x) \setminus \overline{Q}^\nu_{\rho}(x))  \\
& \leq \NNN c_2\sigma^+\rho^{n-1}((1+\delta+\lambda)^{n-1}-1)\leq c_2 C_{\eta}\rho^{n-1}((1+\delta+\lambda)^{n-1}-1)\EEE,
\end{align}
\NNN where we recall that $\sigma^+$ is the density of $\mathcal{M}_0^+$ (see \eqref{FonsecaMantegazzalimit}-\eqref{chermisidensity}) and we used that, by virtue of \eqref{controlsurface} and \cite[Theorem 1.3]{Chermisi}, it can be easily shown $\sigma^+ \leq C_{\eta}$. \EEE
Up to letting $\lambda \to 0$, and combining \eqref{eq3.25bis}--\eqref{eq3.27}, the thesis follows. \EEE
\end{proof}
\begin{proposition}
\label{prop5.1}
Let $(f_\varepsilon)_\varepsilon \subset \mathcal{F}$ and let $f_{0}$ and $(\varepsilon_j)_j$ be as in Theorem \ref{compintrepr}. Then it holds
\begin{equation*}
f_{0}(x,\nu)=f^\prime(x,\nu)=f^{\prime \prime}(x,\nu)    
\end{equation*}
for every $x \in \mathbb{R}^n$ and $\nu \in \mathbb{S}^{n-1}$, where $f^\prime$ and $f^{\prime \prime}$ are as in \eqref{f´} and \eqref{f´´} respectively, but with $\varepsilon$ replaced by $\varepsilon_j$.
\end{proposition}
\begin{proof}
Notice that by definition $f^{\prime \prime}\geq f^{\prime}$, thus it is sufficient to prove $f^\prime \geq f_{0}\geq f^{\prime \prime}$. In the following, for notational convenience, we will still denote with $\varepsilon$ the subsequence $\varepsilon_j$ of Theorem \ref{compintrepr}.\\
\emph{Step 1}: In this step we prove $f_{0}\leq f^\prime$. Fix $x \in \mathbb{R}^n$, $\nu \in \mathbb{S}^{n-1}$ and $\rho>0$. 
Because of $(\NNN f_2\EEE)$, we have $m_{\mathcal{E}_{\varepsilon}}({u}_{x,\varepsilon}^\nu,Q_{\rho}^\nu(x))\geq c_1 m_{\mathcal{M}^-_{\varepsilon}}({u}_{x,\varepsilon}^\nu,Q_{\rho}^\nu(x))\geq c_1 \Tilde{m}_{\mathcal{M}^-_{\varepsilon}}({u}_{x,\varepsilon}^\nu,Q_{\rho}^\nu(x))$, where $\Tilde{m}_{\mathcal{M}^-_{\varepsilon}}({u}_{x,\varepsilon}^\nu,Q_{\rho}^\nu(x))$ is defined as in \eqref{FMinfimumproblem}. Hence, by virtue of Proposition \ref{claimsection6}, we have $m_{\mathcal{E}_{\varepsilon}}({u}_{x,\varepsilon}^\nu,Q_{\rho}^\nu(x))>-1$ for every $\varepsilon$ small enough. Let $\delta>0$ and let $u_\varepsilon$ admissible for $m_{\mathcal{E}_{\varepsilon}}({u}_{x,\varepsilon}^\nu,Q_{\rho}^\nu(x))$ such that
\begin{equation}
\label{eq5.4}
    \mathcal{E}_\varepsilon(u_\varepsilon,Q_{\rho}^\nu(x))\leq m_{\mathcal{E}_{\varepsilon}}({u}^{\nu}_{x,\varepsilon}, Q^\nu_{\rho}(x))+\delta \rho^{n-1}\leq (c_2 C_{\eta}+\delta)\rho^{n-1},
\end{equation}
where in the last inequality we used that ${u}^{\nu}_{x,\varepsilon}$ is also admissible for $m_{\mathcal{E}_{\varepsilon}}({u}^{\nu}_{x,\varepsilon}, Q^\nu_{\rho}(x))$ and \eqref{controlsurface}. Clearly, \eqref{eq5.4} implies
\begin{equation*}
    \sup_{\varepsilon>0} \mathcal{E}_\varepsilon(u_\varepsilon,Q_{\rho}^\nu(x)) < \infty.
\end{equation*}
We can extend $u_\varepsilon$ to $\mathbb{R}^n$ by setting it equal to ${u}^{\nu}_{x,\varepsilon}$ outside $Q^{\nu}_{\rho}(x)$.
Notice that $(u_\varepsilon)_\varepsilon \subset W_{\mathrm{loc}}^{2,2}(\mathbb{R}^n)$ and that, because of $\eqref{growthcondition}$ and \cite[Theorem 1.1]{Chermisi}, for every $A \in \mathcal{A}_1$ such that $Q^\nu_{\rho}(x)\subset \subset A$, there exists a function $u \in L_{\mathrm{loc}}^2(\mathbb{R}^n)\cap BV(A;\{-1,1\})$ such that
\begin{equation*}
    u_\varepsilon \to u \:\: \text{in}\:\: L^2(A)\:\: \text{up to subsequence},
\end{equation*}
and $u=\overline{u}^{\nu}_{x}$ in a neighborhood of $\partial A$. Let $\delta>0$ and $A \in \mathcal{A}_1$ such that $Q^{\nu}_{\rho}(x)\subset \subset A \subset \subset Q^{\nu}_{(1+\delta)\rho}(x)$. By virtue of Lemma \ref{inneregularenvelope} and Lemma \ref{subadditivitym0} we have
\begin{align*}
&m_{\mathcal{E}_0}(\overline{u}^\nu_x,Q_{(1+\delta)\rho}^\nu(x))-c_2C_{\eta}\rho^{n-1}((1+\delta)^{n-1}-1)\leq m_{0}(\overline{u}^\nu_x, A)\leq \mathcal{E}_{0}(u,A) \\
&\leq \mathcal{E}_{0}(u,Q^\nu_{(1+\delta)\rho}(x))\leq \liminf\limits_{\varepsilon \to 0}\mathcal{E}_\varepsilon(u_\varepsilon,Q^\nu_{\rho}(x))+\mathcal{E}_\varepsilon({u}^\nu_{x,\varepsilon},Q^\nu_{(1+\delta)\rho}(x)\setminus \overline{Q}_{\rho}^\nu(x)),      
\end{align*}
that together with \eqref{controlsurface} implies
\begin{equation}
\label{eq5.5}
m_{0}(\overline{u}^\nu_x, Q_{(1+\delta)\rho}^\nu(x))\leq\liminf\limits_{\varepsilon \to 0}m_{\mathcal{E}_{\varepsilon}}({u}^{\nu}_{x,\varepsilon}, Q^\nu_{\rho}(x))+\delta \rho^{n-1}+2c_2C_\eta\delta\rho^{n-1}.     
\end{equation}
By dividing at both sides of $\eqref{eq5.5}$ by $\rho^{n-1}$, sending $\rho \to 0$ and using \eqref{cellformula}, we get
\begin{equation*}
(1+\delta)^{n-1}f_{0}(x,\nu)\leq f^\prime(x,\nu)+\delta+2c_2C_\eta((1+\delta)^{n-1}-1). 
\end{equation*}
Finally sending $\delta \to 0$ in the last equation implies $f_{0}\leq f^\prime$. \\
\emph{Step 2}: In this step we will prove $f^{\prime \prime}\leq f_{0}$. Let $\delta \in (0,1)$ and let $u$ be admissible for $m_{0}(u^\nu_x,Q^{\nu}_{\rho}(x))$ and such that
\begin{equation}
\label{eq5.6}
\mathcal{E}_{0}(u,Q_{\rho}^\nu(x))\leq m_{0}(\overline{u}^\nu_x,Q^\nu_{\rho}(x))+\delta \rho^{n-1}.
\end{equation}
Let $A^{\prime} \in \mathcal{A}_1$ with $Q^\nu_{(1-\delta)\rho}(x)\subset \subset A^{\prime} \subset \subset Q^\nu_{\rho}(x)$.
By virtue of the $\Gamma$-limsup inequality induced by Theorem \ref{compintrepr}, there exists a sequence $(u_\varepsilon)_\varepsilon \subset W^{2,2}(A^{\prime})$ such that $u_\varepsilon \to u$ strongly in $L^2(A^{\prime})$ and 
\begin{equation}
\label{eq5.7}
\limsup\limits_{\varepsilon \to 0}\mathcal{E}_\varepsilon(u_\varepsilon,A^{\prime})\leq \mathcal{E}_{0}(u,A^{\prime}).  
\end{equation}
We extend $u$ to $\mathbb{R}^n$ by setting it equal to $\overline{u}^\nu_x$ outside $A^{\prime}$. Let $0<\rho^{\prime }<\rho$, with $\rho^\prime$ close enough to $\rho$ to ensure
\begin{equation}
\label{eq5.8}
 u=u^\nu_{x}\:\: \text{on} \:\: Q^\nu_{\rho}(x)\setminus \overline{Q}_{\rho^{\prime}}^\nu(x).   
\end{equation}
Up to take $\delta$ small enough, we can also suppose $Q_{\rho^{\prime}}^\nu(x)\subset \subset A^{\prime}$. Let $A$, $A^\prime, U\in \mathcal{A}_1$ such that $Q^\nu_{\rho^{\prime}}(x)\subset \subset A\subset \subset A^{\prime} \subset \subset U \subset \subset Q_{\rho}^{\nu}(x)$ and $B \defas U \setminus \overline{A} \in \mathcal{A}_1$. We apply Proposition \ref{fundamentalestimate} getting a sequence $(\Tilde{u}_\varepsilon)_\varepsilon \subset W^{2,2}(U)$ such that $\Tilde{u}_\varepsilon=u_\varepsilon$ on $A$, $\Tilde{u}_{\varepsilon}={u}^\nu_{x,\varepsilon}$ on $B\setminus \overline{A}^\prime$ and
\begin{equation}
\label{eq5.9bis}
\begin{split}
\limsup\limits_{\varepsilon \to 0}\mathcal{E}_\varepsilon(\Tilde{u}_\varepsilon,U)\leq \limsup\limits_{\varepsilon \to 0}\mathcal{E}_\varepsilon(u_\varepsilon,A^{\prime})+\mathcal{E}_\varepsilon({u}^{\nu}_{x,\varepsilon}, U \setminus \overline{A})\leq\limsup\limits_{\varepsilon \to 0}\mathcal{E}_\varepsilon(u_\varepsilon,A^{\prime})+c_2\mathcal{M}^{+}_\varepsilon({u}^{\nu}_{x,\varepsilon},Q_{\rho}^{\nu}(x)\setminus \overline{Q}^{\nu}_{\rho^{\prime }}(x))\\ \leq \mathcal{E}_0(u,A^{\prime})+c_2C_{\eta}(\rho^{n-1}-(\rho^{\prime })^{n-1})\leq m_{0}(\overline{u}^\nu_x,Q^\nu_{\rho}(x))+\delta\rho^{n-1}+c_2C_{\eta}(\rho^{n-1}-(\rho^{\prime })^{n-1}), \:\: 
\end{split}
\end{equation}
where we used \eqref{controlsurface}, \eqref{eq5.8} and so $\Vert u_\varepsilon - {u}^\nu_{x,\varepsilon}\Vert_{L^2((A^\prime\setminus \overline{A})\cap B)}\to 0$, \eqref{eq5.6} and \eqref{eq5.7}. Now, we extend $\Tilde{u}_{\varepsilon}$ on all $Q_{\rho}^\nu(x)$ by setting it equal to $u^\nu_{x,\varepsilon}$ outside $U$. Thus, we get 
\begin{equation}
    \mathcal{E}_{\varepsilon}(\Tilde{u}_{\varepsilon},Q_{\rho}^\nu(x))\leq \mathcal{E}_{\varepsilon}(\Tilde{u}_{\varepsilon},U)+c_2C_{\eta}(\rho^{n-1}-(\rho^{\prime})^{n-1}),
\end{equation}
where we used $Q_{\rho}^\nu(x)\setminus \overline{U}\subset Q_{\rho}^\nu(x) \setminus \overline{Q}_{\rho^{\prime }}^\nu(x)$, $\NNN(f_3)\EEE$, and \eqref{controlsurface}. In this way, we can notice that since $\Tilde{u}_\varepsilon={u}_{x,\varepsilon}^\nu$ on $Q^\nu_{\rho}(x)\setminus \overline{Q}_{\rho^{\prime }}^\nu(x)$, $\Tilde{u}_\varepsilon$ is admissible for $m_{\mathcal{E}_{\varepsilon}}({u}^\nu_{x,\varepsilon},Q^{\nu}_{\rho}(x))$ and so \eqref{eq5.9bis} implies
\begin{equation}
\label{eq4.11}
\limsup\limits_{\varepsilon \to 0}m_{\mathcal{E}_{\varepsilon}}({u}^\nu_{x,\varepsilon},Q^{\nu}_{\rho}(x))\leq m_{0}(u^\nu_x,Q^\nu_{\rho}(x))+\delta\rho^{n-1}+2c_2C_{\eta}(\rho^{n-1}-(\rho^{\prime})^{n-1}). 
\end{equation}
Finally, $f^{\prime \prime}\leq f_{0}$ follows by sending $\rho^{\prime }\uparrow \rho$ in \eqref{eq4.11} and then arguing like at the end of \emph{Step 1}. 
\end{proof}
\begin{lemma}
\label{increasingsetfunctionlemma}
Let $\mathcal{E}_{\varepsilon}$ be as in \eqref{defenergies}. For every $u \in L_{\mathrm{loc}}^2(\mathbb{R}^n)$ and $A \in \mathcal{A}$, let $\mathcal{E}^{\prime}_0$ and $\mathcal{E}^{\prime \prime}_0$ be as in \eqref{gammaliminfsup}. Then, the set functions $\mathcal{E}^\prime_0(u,\cdot)$ and $\mathcal{E}^{\prime\prime}_0(u,\cdot)$ are increasing on $\mathcal{A}_1$.
\end{lemma}
\begin{proof}
Let $A,B\in \mathcal{A}_1$ with $A \subset B$. Because of \eqref{FonsecaMantegazzalimit} and \eqref{eq4.30}, it is sufficient to show the lemma only for $BV$ functions taking values in $\{-1,1\}$. Let $(u_\varepsilon)_{\varepsilon} \subset L_{\mathrm{loc}}^2(\mathbb{R}^n)$ and $u \in BV(B;\{-1,1\})$ be such that $u_{\varepsilon}\to u$ strongly in $L^2(B)$. Then we have
\begin{align}
\label{C1correction}
 \liminf\limits_{\varepsilon \to 0}\mathcal{E}_{\varepsilon}(u_{\varepsilon},B)&=\liminf\limits_{\varepsilon\to 0}\mathcal{E}_{\varepsilon}(u_{\varepsilon},B)-c_1\mathcal{M}^{-}_{\varepsilon}(u_{\varepsilon},B)+c_1\mathcal{M}^{-}_{\varepsilon}(u_{\varepsilon},B)\nonumber \\
&\geq \liminf\limits_{\varepsilon\to 0}\mathcal{E}_{\varepsilon}(u_{\varepsilon},A)-c_1\mathcal{M}^{-}_{\varepsilon}(u_{\varepsilon},A)+c_1 \mathcal{M}^{-}_{\varepsilon}(u_{\varepsilon},B)  \\
&\geq \liminf\limits_{\varepsilon\to 0}\mathcal{E}_{\varepsilon}(u_{\varepsilon},A)+c_1 \mathcal{M}_0^-(u,B \setminus \overline{A})\geq \liminf\limits_{\varepsilon\to 0}\mathcal{E}_{\varepsilon}(u_{\varepsilon},A)  \nonumber,
\end{align}
where we used that $\mathcal{E}_\varepsilon - c_1\mathcal{M}^{-}_{\varepsilon}$ has positive density and so defines an increasing set function, and that the $\Gamma-\liminf\limits$ inequality in \cite[Theorem 1.3]{Chermisi} holds also for \NNN open and bounded subsets \EEE once \NNN we \EEE assume $u \in BV(B;\{-1,1\})$. Indeed, for the proof of the $\Gamma-\liminf\limits$ inequality in \cite[Theorem 1.3]{Chermisi}, the condition to have $C^1$ boundary is used only in order to apply the compactness result \cite[Theorem 1.1]{Chermisi} and to apply \cite[Theorem 1.2]{Chermisi} in order to bound the integrands of $\mathcal{M}^-_{\varepsilon}(u_{\varepsilon}, B)$ in $L^1(B)$ (and thus the integrands of $\mathcal{M}^-_{\varepsilon}(u_{\varepsilon}, B\setminus \overline{A})$ are bounded in $L^1(B \setminus \overline{A})$ as well). By passing in \eqref{C1correction} to the infimum along all sequences $u_{\varepsilon}$ converging to $u$ strongly in $L^2(B)$, it easily follows $\mathcal{E}_0^\prime(u,B)\geq \mathcal{E}_0^\prime(u,A)$. Similarly, it can be proved $\mathcal{E}_0^{\prime \prime}(u,B)\geq \mathcal{E}_0^{\prime\prime}(u,A)$.   
\end{proof}

\begin{lemma} 
\label{appendix1}
Let $f$ be a random density, $A\in \mathcal{A}$, and $\Tilde{u}\in W^{2,2}(A)$. For every $l \in \mathbb{N}$ define
\begin{equation*}
    f^l(\omega,x,u,\xi,\zeta)\defas \begin{cases}
        f(\omega,x,u,\xi,\zeta)\:\: \text{if}\:\: f(\omega,x,u,\xi,\zeta)\geq -l, \\
        -l \quad \quad \quad \quad \quad \: \text{otherwise},
    \end{cases}
\end{equation*}
and denote by $\mathcal{E}^l$ the corresponding energy. Then, it holds
\begin{equation*}
    \lim\limits_{l \to \infty}m_{\mathcal{E}^l[\omega]}(\Tilde{u},A)=m_{\mathcal{E}[\omega]}(\Tilde{u},A),
\end{equation*}
for every $\omega \in \Omega$.
\end{lemma}
\begin{proof}
Since $f^l\geq f$, it clearly holds
\begin{equation*}
    \liminf\limits_{l \to \infty}m_{\mathcal{E}^l[\omega]}(\Tilde{u},A)\geq m_{\mathcal{E}[\omega]}(\Tilde{u},A).
\end{equation*}
\NNN Before showing the other inequality, we restrict the space of competitors for $m_{\mathcal{E}}(\Tilde{u},A)$ to the ones only having finite $\mathcal{M}^+(\cdot,A)$. Indeed, let $u \in W^{2,2}(A)$ be such that
\begin{equation*}
    \mathcal{M}^+(u,A)=\infty \:\: \text{but}\:\: \infty>\mathcal{E}[\omega](u,A)\geq c_1\mathcal{M}^-(u,A).
\end{equation*}
Then, we would have $\infty=\mathcal{M}^{+}(u,A)-\mathcal{M}^-(u,A)=2q\int_{A}\vert \nabla u \vert^2 \, \mathrm{d}x$, which contradicts $\nabla u \in L^2(A;\mathbb{R}^n)$.
\EEE
Now notice that $f_l \downarrow f$ pointwise and so, because of $\NNN(f_3)\EEE$, $c_2(W(u)+q \vert \nabla u \vert^2 + \vert \nabla^2 u \vert^2)-f^l(\omega,x,u,\nabla u,\nabla^2 u)\uparrow c_2(W(u)+q \vert \nabla u \vert^2 + \vert \nabla^2 u \vert^2)-f(\omega,x,u,\nabla u,\nabla^2 u)$ pointwise as $l \to \infty$.
Thus, by Monotone Convergence Theorem, it holds 
\begin{equation*}
\mathcal{E}[\omega](u,A)  = \lim\limits_{l \to \infty}\mathcal{E}^l[\omega](u,A)\geq \limsup\limits_{l \to \infty}m_{\mathcal{E}^l[\omega]}(\Tilde{u},A),    
\end{equation*}
for every $u \in \mathcal{S}(\Tilde{u},A)$ \NNN such that $\mathcal{M}^+(u,A)<\infty$\EEE, which implies
\begin{equation*}
m_{\mathcal{E}[\omega]}(\Tilde{u},A)\geq  \limsup\limits_{l \to \infty}m_{\mathcal{E}^l[\omega]}(\Tilde{u},A). 
\end{equation*}
\end{proof}

\addcontentsline{toc}{chapter}{Bibliography}
\bibliographystyle{plain}
\bibliography{bibliography.bib}
\end{document}